\documentclass[11pt]{amsart}
\usepackage{graphicx}
\usepackage{amsmath}
\usepackage{amssymb}
\usepackage{amsthm}
\usepackage{mathtools}
\usepackage{tikz-cd}
\usepackage{dsfont}
\usepackage{color}
\usepackage{enumitem}
\usepackage{tensor}
\usepackage{faktor}
\usepackage{todonotes}
\usepackage{mathrsfs}
\usepackage{float}
\usepackage[protrusion=true]{microtype}
\usepackage[]{hyperref}
\hypersetup{
	pdftitle={Arithmetic of Double Torus Quotients and the Distribution of Periodic Torus Orbits},
	pdfauthor={Ilya Khayutin},
	bookmarksnumbered=true,     
	bookmarksopen=true,         
	bookmarksopenlevel=1,       
	colorlinks=true,            
	pdfstartview=Fit,           
	pdfpagemode=UseOutlines,   
	pdfpagelayout=TwoPageRight
}

\makeatletter
\def\paragraph{
	\@startsection{paragraph}{4}
	\z@{.5\linespacing\@plus.7\linespacing}{-.5em}%
	{\normalfont\itshape}}
\makeatother

\DeclareFontFamily{U}{mathx}{\hyphenchar\font45}
\DeclareFontShape{U}{mathx}{m}{n}{
      <5> <6> <7> <8> <9> <10>
      <10.95> <12> <14.4> <17.28> <20.74> <24.88>
      mathx10
      }{}
\DeclareSymbolFont{mathx}{U}{mathx}{m}{n}
\DeclareFontSubstitution{U}{mathx}{m}{n}
\DeclareMathAccent{\widecheck}{0}{mathx}{"71}

\renewcommand{\restriction}{\mathord{\upharpoonright}}
\newcommand{\dif}{\,\mathrm{d}}

\newcommand{\Fu}{{\mathbb{F}_u}}
\newcommand{\Fhat}{{\widehat{\mathbb{F}}}}
\newcommand{\Lhat}{{\widehat{\mathbb{L}}}}
\newcommand{\uhat}{{\widehat{u}}}
\newcommand{\what}{{\widehat{w}}}
\newcommand{\Fbar}{{\overline{\mathbb{F}}}}

\newcommand{\Quhat}{{\mathbb{Q}_\uhat}}

\DeclareMathOperator*{\Id}{Id}

\DeclareMathOperator*{\Ad}{Ad}

\DeclareMathOperator*{\vol}{vol}
\DeclareMathOperator*{\disc}{disc}
\DeclareMathOperator*{\Nr}{Nr}
\DeclareMathOperator*{\Nrd}{Nrd}
\DeclareMathOperator*{\Trd}{Trd}
\DeclareMathOperator*{\Lie}{Lie}

\DeclareMathOperator*{\Spec}{Spec}
\DeclareMathOperator*{\sign}{sign}
\DeclareMathOperator*{\Gal}{Gal}
\DeclareMathOperator*{\Hom}{Hom}
\DeclareMathOperator*{\Mor}{Mor}
\DeclareMathOperator*{\End}{End}
\DeclareMathOperator*{\Aut}{Aut}
\DeclareMathOperator*{\diag}{diag}
\DeclareMathOperator*{\Span}{Span}

\DeclareMathOperator*{\rank}{rank}

\DeclareMathOperator{\Nrml}{N}

\newcommand{\sslash}{\mathbin{/\mkern-6mu/}}

\newcommand{\bsslash}{\mathbin{\backslash \mkern-6mu \backslash}}

\newcommand*{\sdfaktor}[3]{
  \raisebox{-0.5\height}{\ensuremath{#1}}
  \mkern-2mu\backslash\mkern-3mu
  \raisebox{0.5\height}{\ensuremath{#2}}
  \mkern-3mu\slash\mkern-2mu
  \raisebox{-0.5\height}{\ensuremath{#3}}
}

\newcommand*{\dfaktor}[3]{
  \raisebox{-0.25\height}{\ensuremath{#1}}
  \mkern-4mu\bsslash\mkern-5mu
  \raisebox{0.25\height}{\ensuremath{#2}}
  \mkern-5mu\sslash\mkern-4mu
  \raisebox{-0.25\height}{\ensuremath{#3}}
}

\newcommand*{\lfaktor}[2]{
  \raisebox{-0.5\height}{\ensuremath{#1}}
  \mkern-2mu\backslash\mkern-3mu
  \raisebox{0.5\height}{\ensuremath{#2}}
}

\newtheorem{theorem}{Theorem}
\numberwithin{theorem}{section} 
\newtheorem*{theorem*}{Statement of Theorem}
\newtheorem*{theorem-quote}{Theorem}
\newtheorem{lemma}[theorem]{Lemma}
\newtheorem{proposition}[theorem]{Proposition}
\newtheorem{corollary}[theorem]{Corollary}

\theoremstyle{definition}
\newtheorem{definition}[theorem]{Definition}
\newtheorem{example}[theorem]{Example}
\theoremstyle{remark}
\newtheorem{remark}[theorem]{Remark}

\title[Double Torus Quotients]{Arithmetic of Double Torus Quotients and the Distribution of Periodic Torus Orbits}
\author[I. Khayutin]{Ilya Khayutin}
\address{Department of Mathematics, Princeton University, Princeton, NJ 08544, USA}
\address{School of Mathematics, Institute for Advanced Study, Princeton, NJ 08540, USA}

\begin{document}

\setcounter{tocdepth}{1}

\begin{abstract}
We describe new arithmetic invariants for pairs of torus orbits on groups isogenous to an inner form of $\mathbf{PGL}_n$ over a number field. These invariants are constructed by studying the double quotient of a linear algebraic group by a maximal torus.

Using the new invariants we significantly strengthen results towards the equidistribution of packets of periodic torus orbits on higher rank $S$-arithmetic quotients. Packets of periodic torus orbits are natural collections of torus orbits coming from a single adelic torus and are closely related to class groups of number fields. The distribution of these orbits is akin to the distribution of integral points on homogeneous algebraic varieties with a torus stabilizer.

The proof combines geometric invariant theory, Galois actions, local arithmetic estimates and homogeneous dynamics.

\end{abstract}

\maketitle

\bgroup
\hypersetup{linkcolor = blue}
\tableofcontents
\egroup

\section{Introduction}
\subsection{Background}
Let $\mathbf{G}$ be a reductive linear algebraic group and $\mathbf{H}<\mathbf{G}$ a maximal torus, both defined over a number field $\mathbb{F}$. Many interesting affine algebraic varieties arise in the form $\mathbf{V}\coloneqq\faktor{\mathbf{G}}{\mathbf{H}}$. Our main interest is the way that families of integral points of $\mathbf{V}$, defined in a suitable sense,  distribute in $\mathbf{V}(\mathbb{R})$ or in related $S$-arithmetic spaces.
We begin by examining a classical special case. For simplicity of the exposition we assume $\mathbb{F}=\mathbb{Q}$.

\subsection{Points on Sphere}\label{sec:2-sphere}
The $2$-dimensional sphere $\mathbf{S}^{2}$ is the affine variety over $\mathbb{Q}$ defined by the equation
\begin{equation*}
x_1^2+x_2^2+x_3^2=1
\end{equation*}
This variety is equipped with a transitive action of $\mathbf{G}=\mathbf{SO}_3$. The stabilizer of a point is a conjugate of $\mathbf{H}=\mathbf{SO}_{2}<\mathbf{SO}_{3}$. Thus we can identify $\mathbf{S}^{2}\cong\faktor{\mathbf{SO}_{3}}{\mathbf{SO}_{2}}=\faktor{\mathbf{G}}{\mathbf{H}}$. Indeed, $\mathbf{H}={\mathbf{SO}_{2}}$ is a maximal torus of absolute rank $1$ in $\mathbf{G}$. 

Let $d\in\mathbb{N}$, consider the \emph{integral} solutions to the equation
\begin{equation}\label{eq:squares-d}
x_1^2+x_2^3+x_3^2=d
\end{equation}
Denote by $\mathscr{H}_d$ the set of \emph{primitive} integral solutions to \eqref{eq:squares-d}, i.e.\ $\gcd(x_1,x_2,x_3)=1$. Legendre's three squares theorem implies $\mathscr{H}_d\neq\emptyset$ if and only if $d\not\equiv 0,4,7\mod 8$. Accordingly, when writing $d\to\infty$ we consider only those $d$ for which there are primitive integral solutions to \eqref{eq:squares-d}.
If primitive integral solutions exit we denote by $$\widetilde{\mathscr{H}_d}\coloneqq\left\{x/\sqrt{d}\mid x \in \mathscr{H}_d \right\}$$ the radial projection of $\mathscr{H}_d$ to the unit sphere $\mathbf{S}^{2}(\mathbb{R})$. Let $m$ be the unique $\mathbf{SO}_3(\mathbb{R})$-invariant probability measure on the sphere.
We say that the collections $\widetilde{\mathscr{H}_d}$ equidistribute as $d\to\infty$ if for any continuous function $f\colon \mathbf{S}^{2}(\mathbb{R})\to\mathbb{C}$
\begin{equation*}
\frac{1}{|\mathscr{H}_d|}\sum_{x\in\widetilde{\mathscr{H}_d}} f(x)\xrightarrow[d\to\infty]{} \int f(x) \dif m(x)
\end{equation*}
This is weak-$*$ convergence of measures. The question whether equidistribution holds when the limit is taken along all possible $d$'s or along specific subsequences of $d$'s is well studied both using analytic number theory and homogeneous dynamics. 

The equidistribution of $\widetilde{\mathscr{H}_d}$ on the sphere has been settled by Duke \cite{Duke} using automorphic forms and building upon the breakthrough of Iwaniec \cite{Iwaniec}.
Much before the appearance of the harmonic analytic solution, 
Linnik developed in the first half to the 20th century his ``ergodic method''. Linnik was able to prove in the 1950's the equidistribution of integral points for specific sequences of $d$'s \cite{LinnikSphereRussian,LinnikSphereEnglish}. Specifically, he needed the sequence $\left(d_i\right)_{i=1}^\infty$ to satisfy a splitting condition; he assumed  that there is a fixed odd prime $p$ such that $p$ splits in all the field extension $\mathbb{Q}(\sqrt{-d_i})$. Under the hood of Linnik's original proof is the action of the \emph{split} torus $\mathbf{SO}_2(\mathbb{Q}_p)$. 
The $S$-arithmetic quotient
\begin{equation*}
\lfaktor{\Gamma}{G_S}=\lfaktor{\mathbf{SO}_{3}(\mathbb{Z}\left[1/p\right])}{\mathbf{SO}_{3}(\mathbb{R})\times \mathbf{SO}_{3}(\mathbb{Q}_p)} 
\end{equation*} is a compact extension of $\lfaktor{\mathbf{SO}_3(\mathbb{Z})}{\mathbf{S}^2(\mathbb{R})}$.
The periodic orbits of the torus $\mathbf{SO}_2(\mathbb{R})\times\mathbf{SO}_2(\mathbb{Q}_p)$ on $\lfaktor{\Gamma}{G_S}$ 
can be partitioned into finite arithmetic collection called \emph{packets}. For a sequence of discriminants $d_i$ satisfying Linnik's splitting condition, the 
equidistribution of the points $\widetilde{\mathscr{H}_{d_i}}$ follows from the equidistribution of packet of $\mathbf{SO}_2(\mathbb{R})\times\mathbf{SO}_2(\mathbb{Q}_p)$-periodic orbits.

\subsubsection{}
In this paper we make progress on the problem of understanding the distribution of packets of periodic orbits for maximal tori in higher rank groups. Because we focus on groups of type $A_n$ it is useful for us to consider the sphere as a homogeneous space for the projective group of units  of the Hamilton quaternion algebra. 

More generally, for any rational $3$-dimensional quadratic space $\mathcal{Q}=(\mathbb{Q}^3,q)$  there exists a quaternion algebra $\mathbf{B}$, defined over $\mathbb{Q}$, such that $\mathcal{Q}$ is isomorphic to the quadratic space $(\mathbf{B}^0(\mathbb{Q}),\Nrd)$, where $\mathbf{B}^0(\mathbb{Q})$ is the space of traceless quaternions and $\Nrd$ is the reduced norm form. The projective group of units $\mathbf{PB}^\times$ acts isometrically on $\mathbf{B}^0$ by conjugation. This action furnishes an isomorphism of algebraic groups over $\mathbb{Q}$ between the projective group of units $\mathbf{PB}^\times$ and $\mathbf{SO}(q)$. Hence we see that $\mathbf{SO}(q)\simeq \mathbf{PB}^\times$ is an inner form of $\mathbf{PGL}_2$. Our main focus is groups isogenous to an inner form of $\mathbf{PGL}_n$ for $n\geq 2$.

To make the connection to the example of the sphere completely explicit, we mention that
the quadratic form $x_1^2+x_2^2+x_3^2$ corresponds to the reduced norm form of the Hamilton quaternion algebra, which is ramified exactly at $\infty$ and $2$. Similarly, the quadratic form $x_2^2-4x_1^2 x_3^2$ corresponds to the reduce norm form of the split matrix algebra $\mathbf{M}_2$.

\subsection{The Modern Variant of Linnik's Method}
In a series of papers Einsiedler, Lindenstrauss, Michel and
Venkatesh \cite{ELMVPeriodic},\cite{ELMVCubic},\cite{ELMVPGL2} have put forward a modern variant of Linnik's method. Their approach relies heavily of the notion of entropy and measure rigidity for diagonalizable actions. In the modern viewpoint, Linnik's method is based upon the dynamics of a an isotropic torus $H$ on an $S$-arithmetic homogeneous space $\lfaktor{\Gamma}{G}$. In the example of the sphere  $G=\mathbf{PB}^\times(\mathbb{R})\times\mathbf{PB}^\times(\mathbb{Q}_p)$ with $p>2$, where $\mathbf{B}$ is the Hamilton quaternion algebra,  and $\Gamma=\mathbf{PB}^\times\left(\mathbb{Z}\left[\frac{1}{p}\right]\right)$. The torus $H$ is the centralizer in $\mathbf{PB}^\times(\mathbb{Q}_p)$ of a point in $\mathscr{H}_d\subset \mathbf{B}^0(\mathbb{Q})\hookrightarrow  \mathbf{B}^0(\mathbb{Q}_p)$.

The Kolmogorov-Sinai entropy is an invariant of a measure preserving system that takes values in $[0,\infty]$ . For $a\in H$ the entropy $h_\mu(a)$ of an $a$-invariant and ergodic Borel probability measure $\mu$ on $\lfaktor{\Gamma}{G}$ is roughly the exponential decay rate in $N$ of the measure of the set of points whose trajectory under $a$ remains uniformly close to the trajectory of a typical point at the times $1,2,\ldots,N$. 

The modern variant of Linnik's method has two distinguished parts.
The gist of the first part is to use an arithmetic well-separateness result for the orbits of $H$, together with a volume computation, to derive an estimate for the metric entropy of any limit measure with respect to an element $a\in H$. In the second part, this entropy estimate is the input to a measure rigidity theorem which supplies the final result about the possible limit measures. In the case of the modern proof of Linnik's theorem \cite{ELMVPGL2} one establishes first maximal entropy in the limit and then uses the, relatively simple, fact that maximal entropy measures are almost unique.
This general proof scheme has been exploited by Einsiedler, Lindenstrauss, Michel and Venkatesh in \cite{ELMVPeriodic,ELMVCubic} to obtain fundamentally new results in higher rank analogues.

The source of the entropy in the limit in Linnik's original argument for the $2$-sphere is Linnik's basic lemma. Its proof depends on counting representations of a binary integral quadratic form by a ternary one. This counting argument is converted to a result about average separation between orbits using a basic invariant - the Euclidean inner-product of two integral points on the sphere. 

{The main contribution of the work presented here is a generalization of this invariant to higher rank and the analysis of some fundamental properties of the new invariants.} It is using these results that we are able to prove an improved bound on the limit entropy in higher rank spaces.

\subsection{Periodic torus orbits on \texorpdfstring{${\mathbf{PGL}_n(\mathbb{Z})}\backslash{\mathbf{PGL}_n(\mathbb{R})}$}{PGLn(Z)\PGLn(R)}}\label{sec:pgln-periodic}
To illustrate our main theorem
consider first the case that $\mathbf{G}=\mathbf{PGL}_n$ and $H$ is a maximal torus in $\mathbf{G}(\mathbb{R})$. At first reading one can restrict to the case that $H$ is the group of diagonal matrices, but some of the most interesting implications of our theorems are when $H$ is isotropic but not split over the reals.

Following \cite{ELMVPeriodic} we say that an $H$-orbit on 
\begin{equation*}
X\coloneqq\lfaktor{\mathbf{PGL}_n(\mathbb{Z})}{\mathbf{PGL}_n(\mathbb{R})}
\end{equation*} 
is periodic if it supports an $H$-invariant probability measure. Let $g\in\mathbf{PGL}_n(\mathbb{R})$. If   $gHg^{-1}=\mathbf{T}(\mathbb{R})$ where $\mathbf{T}<\mathbf{PGL}_n$ is a maximal torus defined and anisotropic over $\mathbb{Q}$ then $g^{^-1}\Gamma g \cap H$ is a lattice in $H$ and the orbit $\mathbf{PGL}_n(\mathbb{Z})gH$ is periodic.  We say that a periodic orbit $\mathbf{PGL}_n(\mathbb{Z})gH$ is \emph{algebraic} if $gHg^{-1}$ admits the structure of a rational torus as above. For $H$ split, e.g.\ $H$ is the group of diagonal matrices, every periodic orbit is algebraic. This statement is wrong for non-split tori. For example, if $n=2$ then any orbit of the compact torus $\mathbf{SO}_2(\mathbb{R})$ is periodic. 

Assume that $H=\mathbf{H}(\mathbb{R})$ and that $\mathbf{H}<\mathbf{PGL}_n$ is a maximal torus defined over $\mathbb{Q}$. Periodic $H$-orbits tautologically correspond to some $\mathbf{PGL}_n(\mathbb{Z})$ orbits on $\faktor{\mathbf{PGL}_n(\mathbb{R})}{\mathbf{H}(\mathbb{R})}$. Algebraic periodic $H$-orbits on $\lfaktor{\mathbf{PGL}_n(\mathbb{Z})}{\mathbf{PGL}_n(\mathbb{R})}$  necessarily correspond to orbits of  $\mathbf{PGL}_n(\mathbb{Z})$ on \emph{rational points} of the affine variety underlying $\faktor{\mathbf{PGL}_n(\mathbb{R})}{\mathbf{H}(\mathbb{R})}$. 

The algebraic periodic $H$-orbits on $X$ are partitioned into natural finite collections called \emph{packets}, whose formal definition we review in \S \ref{sec:general-setting}. These packets generalize the collections $\mathscr{H}_d$ from \S \ref{sec:2-sphere}. 
To each packet of $H$-orbits we can associate an order in a degree $n$ field extension of $\mathbb{K}/\mathbb{Q}$. If $\mathbb{K}$ has $r_1$ real places and $r_2$ inequivalent complex places then necessarily
$\rank_{\mathbb R} H=r_1+r_2-1$. The set of $H$-orbits in  a single packet is a principal homogeneous space for the Picard group of the associated order, see \cite[Corollary 4.4]{ELMVPeriodic}. The simplest packets are associated to maximal orders, viz.\ rings of integers of degree $n$ number fields. 
To each packet one can associate the \emph{discriminant}. This is an arithmetic property of the packet. In the special case discussed in this section it is just the discriminant of the associated order, which coincides with the discriminant of the number field if the order is maximal.

Each packet supports a canonical $H$-invariant probability measure. It is defined as the uniform average over the $H$-invariant probability measures supported on the individual periodic orbits in the packet.
Suppose we are given a sequence of packets $\left\{Z_i\right\}_{i=1}^\infty$ with discriminants $D_i$ and $H$-invariant measures $\mu_i$. How do the measures $\mu_i$ distribute on $X$ when $D_i\to_{i\to\infty}\infty$?
For $n=2$ the space $\lfaktor{\mathbf{PGL}_2(\mathbb{Z})}{\mathbf{PGL}_2(\mathbb{R})}$ is the unit tangent bundle of the modular curve. If $H$ is the group of diagonal matrices then the periodic $H$-orbits are exactly the non-divergent closed geodesics. When $n=2$ the convergence of periodic measures on $H$-invariant packets to the Haar measure is again a theorem of Duke \cite{Duke}, cf.\ \cite{ELMVPGL2,Skubenko}. For $n=3$ see the recent results of \cite{ELMVCubic}. For higher $n$ our state of knowledge is not as satisfactory. Assuming non-escape of mass, only a lower bound on the metric entropy of a limit measure was known \cite[Theorem 3.1]{ELMVPeriodic} and \emph{only for $H$ split}, see Remark \ref{rem:elmv-bound}. We present a significant improvement for this lower bound for packets satisfying a Galois condition and without assuming that $H$ is split.

For any $a\in H$ and $\nu$ a probability measure on $X$ invariant under $H$, denote by $h_\nu(a)$ the metric entropy of $\nu$ with respect to the action by $a$. Let $h_\mathrm{Haar}(a)$ be the entropy of the Haar probability measure on $X$, which is the unique measure of maximal entropy.
In this setting we have the following special case of our main theorem.
\begin{theorem}
Set $H<  \mathbf{PGL}_n(\mathbb{R})$ to be a fixed non-compact maximal torus.
Let $\left\{Z_i\right\}_{i=1}^\infty$ be a sequence of packets of periodic $H$-orbits with discriminants $D_i$ associated to maximal orders in number fields. Assume $D_i\to_{i\to\infty}\infty$.

Let $\mathbb{K}_i$ be a degree $n$ field extension of $\mathbb{Q}$ associated to $Z_i$. Denote by $\mathbb{L}_i$ the Galois closure of $\mathbb{K}_i$. Assume for all $i$ that $\Gal(\mathbb{L}_i/\mathbb{Q})$ is 2-transitive when considered as a subgroup of $S_n$.

Let $\mu_{i}$ be the canonical $H$-invariant probability measure supported on $Z_i$. 
If $\mu_i$ converges in the weak-$*$ topology to a probability measure $\mu$ on $X$ then for any $a\in H$
\begin{equation*}
h_\mu(a)\geq \frac{h_\mathrm{Haar}(a)}{2(n-1)}
\end{equation*}
\end{theorem}
This theorem should be compared with \cite[Theorem 3.1]{ELMVPeriodic} which for $H$ -- the \emph{split torus} of diagonal matrices and
\begin{equation*}
a=\diag\left(\exp\left(\frac{n-1}{2}\right),\exp\left(\frac{n-3}{2}\right),\ldots,
 \exp\left(-\frac{n-1}{2}\right)\right) 
\end{equation*}
provides the bound $h_\mu(a)\geq \frac{3 h_\mathrm{Haar}(a)}{(n+1)n(n-1)}$ without assuming a Galois condition and without a restriction to maximal orders. 

The strength of our result compared to \cite[Theorem 3.1]{ELMVPeriodic} is that it is valid even when $H$ is isotropic but not split over $\mathbb{R}$ and the significantly improved bound on the entropy under a Galois condition, especially for large values of $n$. 

Measure rigidity results for diagonalizable actions imply in favorable situation that a measure invariant and ergodic under a higher rank torus is algebraic, at least assuming positive entropy.
For $H$ split our entropy bound and the measure rigidity results from \cite{EntropySArithmetic} imply that $\mu$ is not supported on a single periodic measure of a closed subgroup $L<\mathbf{PGL}_n(\mathbb{R})$ all whose simple factors have rank $\leq R$ where $R$ satisfies $(R+2)(R+1)R-1<n$, see Corollary \ref{cor:no-single-periodic}. Moreover, if $n$ is prime then no possible intermediate algebraic measures can arise, cf.\ \cite{LindenstraussWeiss,ELMVPeriodic} and we can deduce that $\mu\geq m/(2(n-1))$ where $m$ is the normalized Haar measure.
Indeed, an improvement in the bound on the asymptotic entropy provides more control on the possible algebraic measures appearing in the ergodic decomposition of $\mu$. Extremely strong bounds on the asymptotic entropy can even imply partial non-escape of mass. These results where pioneered in \cite{ELMVPGL2} and studied further in \cite{EK,KadyrovEscape,KadyrovH,KadyrovPohl,EKP,KKLM} and by a different method in \cite{KhLD}. The required entropy bounds for non-escape of mass seem to be far from what is attainable with existing methods. We stress that the general version of our theorem holds for any inner-form of $\mathbf{PGL}_n$, including inner-forms arising from division algebras, for which there is no possibility of escape of mass.

The entropy of the measure $\mu$ with respect to the action by a semisimple element $a\in\mathbf{PGL}_n(\mathbb{R})$ is related to the dimension of $\mu$ along its sections in some directions transversal to the direction of the
$a$-action; specifically, directions in the stable foliation.  When the centralizer of $a$ is not minimal the measure might have non-trivial sections in many transversal direction which do not affect the entropy. This is a significant problem when trying to bound the entropy with respect to an action by an element which has a big centralizer, as is the case for non-split isotropic tori in our setting. The action of the Galois group of the splitting field on appropriate invariants we construct is an important ingredient in overcoming this problem. The method of proof of \cite[Theorem 3.1]{ELMVPeriodic} is fundamentally limited to elements $a$ whose centralizer is minimal, hence the restriction to split tori in that result. 

To prove this theorem we construct an integral model for the GIT double quotient of $\mathbf{G}$ by a maximal torus and show a bound on the separation of points which are integral in an appropriate sense. A major difficulty is that in rank bigger then $1$ the double quotient space is not an affine space. This is a significant obstacle if we were to try and show maximal entropy in the limit.

A main tool we use is geometric invariant theory. Another novel aspect of our method is that we are able to utilize properties of the Galois groups of splitting fields.

\subsection{Inner forms of type \texorpdfstring{$A_n$}{An} and homogeneous toral sets}\label{sec:general-setting}
We move to describe the setting of our main theorem in its general form. The theorem holds for quotients of inner forms of type $A_n$ by congruence lattices. These can all be described using central simple algebras. We fix for the rest of the paper a number field $\mathbb{F}$. Let $\mathbf{B}$ be the affine variety representing a central simple algebra over $\mathbb{F}$. Let $\mathbf{G}$ be a linear reductive group defined over $\mathbb{F}$ and isogenous to $\mathbf{PGL}_1(\mathbf{B})$, i.e.\ $\mathbf{G}$ has a finite center $\mathbf{Z}$ and $\mathbf{G}^\mathrm{ad}\coloneqq\lfaktor{\mathbf{Z}}{\mathbf{G}}\simeq \mathbf{PGL}_1(\mathbf{B})$. We fix once and for all an isogeny
\begin{equation*}
\mathbf{SL}_1(\mathbf{B})\to\mathbf{G}\to\mathbf{PGL}_1(\mathbf{B})
\end{equation*}
where $\mathbf{SL}_1(\mathbf{B})$ is the group of norm $1$ units in $\mathbf{B}$. This group is the simply-connected cover of $\mathbf{PGL}_1(\mathbf{B})$. Recall that an algebraic group $\mathbf{M}$ defined over $\mathbb{F}$ is anisotropic over $\mathbb{F}$ if it admits no non-trivial characters $\mathbf{M}\to\mathbb
 G_m$ over $\mathbb{F}$. The group $\mathbf{G}$ satisfies this property.

Denote by $\mathbb{A}$ the ring of adeles of $\mathbb{F}$.
Recall that the adelic points of $\mathbf{G}$ are defined as
\begin{equation*}
\mathbf{G}(\mathbb{A})=\left\{\left(g_u\right)_u \mid g_u\in\mathbf{G}(\mathbb{F}_u) \textrm{ and } g_u\in K_u \textrm{ for almost all } u\right\}
\end{equation*}
where the index $u$ runs over all the places of $\mathbb F$ and $K_u$ is a fixed compact-open subgroup of $\mathbf{G}(\mathbb F_u)$ for all $u$ non-archimedean. For almost all $u$ we need $K_u<\mathbf{G}(\mathbb F_u)$ to be a maximal compact subgroup.
The subgroup $\mathbf{G}(\mathbb{F})$ embeds diagonally in $\mathbf{G}(\mathbb{A})$ and it is discrete in the induced topology. The quotient space $\lfaktor{\mathbf{G}(\mathbb{F})}{\mathbf{G}(\mathbb{A})}$ carries a  $\mathbf{G}(\mathbb{A})$-invariant measure that is unique up to scaling. Because $\mathbf{G}$ is anisotropic over $\mathbb{F}$ the invariant measure is finite and we fix it to be a probability measure. 

A rational anisotropic torus in $\mathbf{G}$ can be constructed in the following way. Let $\mathbb{K}/\mathbb{F}$ be a degree $n$ field extension with a ring embedding $\iota\colon\mathbb{K}\hookrightarrow\mathbf{B}(\mathbb{Q})$. There is an affine subvariety $\mathbf{E}<\mathbf{B}$ representing $\iota(\mathbb{K})$, i.e.\  $\iota(\mathbb{K})\otimes_{\mathbb F} R=\mathbf{E}(R)$ for any commutative algebra $R/\mathbb{F}$. The algebraic subgroup $\mathbf{SL}_1(\mathbf{E})<\mathbf{SL}_1(\mathbf{B})$ is a maximal torus and its isogeny image in $\mathbf{G}$ is a maximal torus in $\mathbf{G}$ defined and anistropic over $\mathbb{F}$. Moreover, all maximal tori in $\mathbf{G}$ that are defined and anisotropic over $\mathbb{F}$ arise this way. A torus $\mathbf{T}$ that is anisotropic over $\mathbb{F}$ has the important property that $\lfaktor{\mathbf{T}(\mathbb{F})}{\mathbf{T}(\mathbb{A})}$ carries a unique $\mathbf{T}(\mathbb{A})$-invariant Borel \emph{probability} measure.

Let $\mathbf{T}<\mathbf{G}$ be defined and anisotropic over $\mathbb{F}$ and fix $g_\mathbb{A}\in\mathbf{G}(\mathbb{A})$. A subset of the form
\begin{equation*}
\mathbf{G}(\mathbb{F})\mathbf{T}(\mathbb{A})g_{\mathbb{A}} \subset \lfaktor{\mathbf{G}(\mathbb{F})}{\mathbf{G}(\mathbb{A})}
\end{equation*}
is called a \emph{homogeneous toral set}. These have been introduced in \cite{ELMVCubic}. Because $\mathbf{T}$ is anisotropic the homogeneous toral set carries a unique $g_\mathbb{A}^{-1} \mathbf{T}(\mathbb{A}) g_\mathbb{A}$-invariant Borel \emph{probability} measure. 

To each homogeneous toral set one can attach a volume and a discriminant, both are non-negative real numbers. The volume is a global quantity that measures the size of the homogeneous toral set in geometric terms, it is a generalization of the notion of the length of a closed geodesic. 
The discriminant $D$ is a product of local discriminants $D_u$ for each place $u$ of $\mathbb{F}$. The local discriminant at the place $u$ is closely related to the notion of the local height of the subspace $\Lie(g_u^{-1}\mathbf{T}(\mathbb{F}_u)g_u)$ in $\Lie(\mathbf{G})(\mathbb{F}_u)$.
These quantities have been formally defined in \cite{ELMVCubic} and we review the definitions in \S \ref{sec:discriminant-volume}. Moreover, if the torus $\mathbf{T}$ is associated to an embedding of the number field $\mathbb{K}/\mathbb{F}$ then we will construct in \S \ref{sec:global-order} an order $\mathcal{R}\subset\mathbb{K}$ satisfying $|\disc(\mathcal{R})|=\prod_{u<\infty} D_u$. This order depends both on the torus $\mathbf{T}$ and the distortion $g_\mathbb{A}\in\mathbf{G}(\mathbb{A})$. We say that a homogeneous toral set is of \emph{maximal type} if this order is the ring of integers in $\mathbb{K}$. Homogeneous toral sets of maximal type are exactly the higher rank generalization of packets of periodic geodesics attached to \emph{fundamental} discriminants.

\subsection{The \texorpdfstring{$S$}{S}-arithmetic setting and packets}
Let $S$ be a finite set of places of $\mathbb{F}$ including all the archimedean ones. Denote $\mathbb{F}_S=\prod_{u\in S} \mathbb{F}_u$. Suppose $S$ includes at least one place over which $\mathbf{B}$ is isotropic and that it is large enough so that $\mathbf{G}(\mathbb{A})$ has class number one with respect to $S$. 
The class number 1 requirement amounts to
\begin{equation*}
\mathbf{G}(\mathbb{A})=\mathbf{G}(\mathbb{F}) \cdot  \mathbf{G}(\mathbb{F}_S) \cdot \prod_{u\not\in S} K_u
\end{equation*}
Every finite set of places $S$ can always be increased to a finite set satisfying the class number $1$ assumption. If $\mathbf{G}=\mathbf{SL}_1(\mathbf{B})$ is simply-connected then the class number assumption is satisfied for any set $S$ as above because of strong approximation.

Denote $G_S=\prod_{u\in S} \mathbf{G}(\Fu)$ and $K^S\coloneqq\prod_{u\not\in S} K_u$.
Let $\Gamma<G_S$ be the congruence lattice $$\Gamma\coloneqq\mathbf{G}(\mathbb{F})\cap K^S$$
Under the assumptions above on $S$ we can identify $\lfaktor{\Gamma}{G_S}$ with the quotient of  $\lfaktor{\mathbf{G}(\mathbb{F})}{\mathbf{G}(\mathbb{A})}$ by $K^S$. 

Fix a homogeneous toral set
\begin{equation*}
Y=\mathbf{G}(\mathbb{F})\mathbf{T}(\mathbb{A})g_\mathbb{A} \subset \lfaktor{\mathbf{G}(\mathbb{F})}{\mathbf{G}(\mathbb{A})}
\end{equation*}
and denote by $\mu_Y$ the invariant probability measure supported on this set. The image of $Y$ under the projection map
\begin{equation*}
\lfaktor{\mathbf{G}(\mathbb{F})}{\mathbf{G}(\mathbb{A})} \xrightarrow{/K^S} \lfaktor{\Gamma}{G_S}
\end{equation*}
is called a \emph{packet} of torus orbits. This packet is invariant under the action of $g_S^{-1}\mathbf{T}(\mathbb{F}_S) g_S$, where $g_S$ is the $S$-adic part of $g_\mathbb{A}$. A simple calculation shows that the packet is a finite collection of $g_S^{-1}\mathbf{T}(\mathbb{F}_S) g_S$-orbits parametrized by a finite abelian group arising as a double quotient of $\mathbf{T}(\mathbb{A})$.
The packet supports the push-forward of the measure $\mu_Y$ which is a uniform average over the periodic measures of the $g_S^{-1}\mathbf{T}(\mathbb{F}_S) g_S$-orbits in the packet. 

The special case discussed in \S \ref{sec:pgln-periodic} corresponds to $\mathbb{F}=\mathbb{Q}$, $\mathbf{B}=\mathbf{M}_n$, $\mathbf{G}=\mathbf{PGL}_n$ and $S=\left\{\infty\right\}$. The packets of periodic algebraic $H$-orbits presented in \S \ref{sec:pgln-periodic} are a special case of the packets described in the adelic language. They all satisfy the assumption $g_\infty^{-1}\mathbf{T}(\mathbb{R}) g_\infty=H$, where $g_\infty$ is the archimedean component of $g_\mathbb{A}$, i.e.\ they are all invariant under a fixed non-compact torus in the archimedean place.

\subsection{Main Theorem: Lower Bound on Asymptotic Entropy}\label{sec:main-thm}
To present our main theorem we need also to introduce the notion of the ramified part of the discriminant of a homogeneous toral set. If the local discriminant of a homogeneous toral set at a place $u$ is denoted by $D_u$ then the ramified local discriminant is defined as
\begin{equation*}
D_\mathrm{ram}=\prod_{\mathbf{B}\textrm{ ramifies at } u} D_u
\end{equation*}
Notice that the product is over the places where $\mathbf{B}$ is ramified and \emph{not} over the ramified places of the torus. In particular, if $\mathbf{B}=\mathbf{M}_n$ is split then $D_\mathrm{ram}=1$ for all homogeneous toral sets.

\begin{theorem}\label{thm: entropy-bound}
Suppose we are given a sequence of homogeneous toral sets of maximal type $Y_i={\mathbf{G}(\mathbb{F})}{\mathbf{T}_i(\mathbb{A})g_i}$. Denote by $D_i$ the global discriminant of $Y_i$ and let $D_{\mathrm{ram},i}$ be the ramified part of the discriminant. Assume $D_i\to_{i\to\infty}\infty$ and $D_{\mathrm{ram},i}=D_i^{o(1)}$. 

\emph{Galois Condition:}
Let $\mathbb{L}_i/\mathbb{F}$ be the splitting field of $\mathbf{T}_i$. We assume for all $i$ that $\Gal(\mathbb{L}_i/\mathbb{F})$ is 2-transitive.

\emph{Splitting Condition:}
We fix a place $\uhat\in S$ such that $\mathbf{B}$ is isotropic over $\uhat$. Denote $\Fhat\coloneqq\mathbb{F}_\uhat$.  Fix a torus $\mathbf{H}<\mathbf{G}_{\Fhat}$ defined over $\Fhat$ and such that\footnote{The theorem is vacuously true without the non-compactness assumption for $\mathbf{H}(\Fhat)$.} $\rank_\Fhat\mathbf{H}>0$. Set $H=\mathbf{H}(\Fhat)$.
Assume  for all $i$ that ${g_{i,\uhat}}^{-1} \mathbf{T}_i g_{i,\uhat}=\mathbf{H}$, i.e.\ all the sets $Y_i$ are invariant in the place $\uhat$ under a fixed isotropic torus $H$.

Let $\mu_{i}$ be the probability measure on $\lfaktor{\Gamma}{G_S}$ induced by the invariant probability measure on the homogeneous toral set $Y_i$. If $\mu_i$ converges in the weak-$*$ topology to a probability measure $\mu$ on $\lfaktor{\Gamma}{G_S}$ then for any $a\in H$ we have
\begin{equation*}
h_\mu(a)\geq \frac{h_\mathrm{Haar}(a)}{2(n-1)}
\end{equation*}
Where $h_\mathrm{Haar}(a)$ is the entropy of the Haar probability measure on $\lfaktor{\Gamma}{G_S}$, which is a measure of maximal entropy.
\end{theorem}
\begin{remark}\label{rem:elmv-bound}
Assume $\mathbf{H}$ is split over $\Fhat$ and let $\Phi$ be the set of roots for the torus $\mathbf{H}$.
The weaker bound
\begin{equation}\label{eq:elmv-bound}
h_\mu(a)\geq \frac{1}{2} \min_{\alpha\in\Phi} |\log|\alpha(a)|_\uhat|
\end{equation}
has been proven in \cite{ELMVPeriodic}[Theorem 3.1] for any $a\in H$ but without any Galois condition or assumption on maximality or ramification. The entropy of the Haar measure is
\begin{equation*}
h_\mathrm{Haar}(a)= \frac{1}{2} \sum_{\alpha\in\Phi} |\log|\alpha(a)|_\uhat|
\end{equation*}

The proof of \eqref{eq:elmv-bound} goes by considering the separation of orbits implied by attaching to the Lie algebra of $\mathbf{T}_i$ an integral point with denominator ${D_i}$ in $\left(\bigwedge^{\rank\mathbf{G}}\Lie(\mathbf{G})(\mathbb{F})\right)^{\otimes2}$. This procedure is essentially the \emph{definition} of the discriminant.

The ability to deduce better bounds using just the adjoint representation on the Lie algebra seems to be limited by the fact that for $\rank\mathbf{G}>1$ most of the subspaces of dimension $\rank\mathbf{G}$ in $\Lie(\mathbf{G})(\mathbb{F})$ do not correspond to tori. 

Without the condition on the Galois group, our method gives the same entropy bound as in \eqref{eq:elmv-bound}. 
\emph{A key new feature in our paper is that we are able to use additional information about the Galois groups of the splitting fields of the tori to give substantially more precise results.} Moreover, we are able to treat homogeneous toral sets which are invariant under a non split isotropic torus.

It is instructive to calculate the different available bounds for the case that $\mathbf{G}=\mathbf{PGL}_n$ over $\mathbb{Q}$, $\uhat=\infty$ -- the real place, $\mathbf{H}$ -- the standard torus of diagonal matrices and
\begin{equation*}
 a=\diag\left(\exp\left(\frac{n-1}{2}\right),\exp\left(\frac{n-3}{2}\right),\ldots,
 \exp\left(-\frac{n-1}{2}\right)\right)
\end{equation*}
The roots of $\mathbf{H}$ are parameterized by pairs $(i,j)$ where $1\leq i\neq j\leq n $. The values of the roots for $a$ as above are
\begin{equation*}
\log|\alpha_{i,j}(a)|=j-i
\end{equation*}
In particular
\begin{align*}
& h_\mathrm{Haar}(a)=\frac{(n+1)n(n-1)}{6}\\
& \frac{h_\mathrm{Haar}(a)}{2(n-1)}=\frac{(n+1)n}{12}\\
& \frac{1}{2}\min_{\alpha\in\Phi} |\log|\alpha(a)|_\infty|=\frac{1}{2}
\end{align*}
A main difference between the bounds is that the new bound grows quadratically with $n$, unlike the bound \eqref{eq:elmv-bound} which is constant. In conjunction with the measure rigidity results of \cite{EKL,EntropySArithmetic} this implies a new modest qualitative restriction on possible limit measures. Specifically, if $\mathbb{F}=\mathbb{Q}$ then $\mu$ is not supported on a single periodic orbit of a reductive group all whose simple parts have small absolute rank compared to $\mathbf{G}$, see Corollary \ref{cor:no-single-periodic}.
\end{remark}

\subsection{An Overview of the Proof}
We present a general overview of the proof completely in adelic terms. 
\paragraph{Decay of mass for thin tubes.} Fix an element $a\in H$. Due to essentially upper semi-continuity of entropy, in order to bound the entropy of the limit measure with respect to $a$, it is enough to show that the $\mu_i$ mass of small tubes in $\lfaktor{\mathbf{G}(\mathbb{F})}{\mathbf{G}(\mathbb{A})}$ decays exponentially fast with a given rate when these tubes are conjugated by $a$. Notice that these tubes contract in a single place $\uhat$ only. 

The rate at which the measure of the contracting tubes decreases is our bound on the entropy with respect to $a$. For this idea to apply easily we need all our homogeneous toral sets to be invariant under the action of this fixed $a$. To apply the measure rigidity results of \cite{EKL,EntropySArithmetic} one might need to assume that $\mathbf{H}$ is split over $\Fhat$.

This method for proving a bound on the entropy is well known and has been used for similar means in \cite{ELMVPeriodic,ELMVPGL2}. 
\paragraph{Double torus quotient invariants.} Let $\mathbf{T}$ be a rational torus appearing in the series $\left\{\mathbf{T}_i\right\}_i$. We construct new polynomial invariants $\Psi_\sigma\colon{\mathbf{G}(\mathbb{F})}\to \mathbb{L}$ where $\mathbb{L}/\mathbb{F}$ is the splitting field of $\mathbf{T}$. These invariants depend on $\mathbf{T}$ and we have one such invariant for each $\sigma$ in the absolute Weyl group of $\mathbf{T}$. It is useful to choose an identification of the Weyl group with the symmetric group $S_n$ and consider the invariants as parametrized by permutations. The polynomials $\Psi_\sigma$ are invariant under both the left and the right action of $\mathbf{T}$.

\paragraph{Application of Geometric invariant theory.} Using GIT we are able to prove for two torus cosets $\delta_L\mathbf{T}(\mathbb{A})$, $\delta_R\mathbf{T}(\mathbb{A})$ with $\delta_L,\delta_R\in\mathbf{G}(\mathbb{F})$ that if for all $\sigma$
\begin{equation*}
 \Psi_\sigma({\delta_L}^{-1}\delta_R)=\Psi_\sigma(\mathrm{id})
\end{equation*}
then necessarily $\delta_L\mathbf{T}(\mathbb{A})=\delta_R\mathbf{T}(\mathbb{A})$.

\paragraph{Arithmetic properties of the invariants.}We show that the invariants $\Psi_\sigma$ satisfy an important arithmetic property. One variant of this result is that for all non-archimedean places $u$ of $\mathbb{F}$ where $\mathbf{B}$ is unramified we have that if ${\delta_L}^{-1}\delta_R$ is in an appropriate identity neighborhood in $\mathbf{G}(\Fu)$, i.e.\ satisfy an integrality condition, then 
\begin{equation}\label{eq:idea-integrality}
|{\Nr}_{\mathbb{L}/\mathbb{F}}\Psi_\sigma({\delta_L}^{-1}\delta_R)|_u
<{D_u}^{\left[\mathbb{L}\colon\mathbb{F}\right]}
\end{equation}
Where $D_u$ is the local discriminant of the homogeneous toral set at the place $u$. Moreover, an analogues result holds in the archimedean places. It is here that we use the maximality assumption for the homogeneous toral set.
\paragraph{Algebraic relations between the invariants and Galois orbits thereof.} We exploit two more tools: the algebraic relations between the invariants $\Psi_\sigma$ for different $\sigma$ and the action of $\Gal(\mathbb{L}/\mathbb{F})$ on each $\Psi_\sigma$. We form a new single invariant $\Psi_\mathscr{C}$ by taking the product of all the Galois conjugates of a specific $\Psi_{\sigma_0}$, where $\sigma_0$ can be taken as any permutation that has no fixed points. The resulting invariant takes values in the fixed field $\mathbb{F}$.

Obviously, ${\Nr}_{\mathbb{L}/\mathbb{F}}\Psi_{\sigma_0}$ is an integral power of $\Psi_\mathscr{C}$ and we can deduce from inequality \eqref{eq:idea-integrality} a similar inequality for $\Psi_\mathscr{C}$.

Under the assumption that $\Gal(\mathbb{L}/\mathbb{F})$ is 2-transitive we are able to combine our understanding of both the Galois orbits of all the invariants $\Psi_\sigma$ and the algebraic relations between them to show that the equality 
\begin{equation*}
\Psi_\mathscr{C}({\delta_L}^{-1}\delta_R)=0
\end{equation*} 
is enough to deduce that $\Psi_\sigma({\delta_L}^{-1}\delta_R)=\Psi_\sigma(\mathrm{id})$ for \emph{all} $\sigma\in S_n$. 

\paragraph{Using the invariants to study thin tubes.} At this point we have most of the tools to apply the method described in the fist step. We study the value of $\Psi_\mathscr{C}$ on pairs of translated cosets $\delta_L \mathbf{T}(\mathbb{A})g$, $\delta_R \mathbf{T}(\mathbb{A})g$ that lie in the same thin tube. By choosing our tube carefully in all the places we are able to insure that the integrality conditions required for \eqref{eq:idea-integrality} hold. Hence we have a bound on the absolute value of $\Psi_\mathscr{C}({\delta_L}^{-1}\delta_R)$ in every place $u$ in terms of the local discriminant $D_u$.

Moreover, at the single place $\uhat$ where the tube is contracted we can given an exponentially good bound on the $\uhat$ absolute value of $\Psi_\mathscr{C}({\delta_L}^{-1}\delta_R)$ for pairs of cosets as above. The rate of this exponential bound affects critically the final bound on the entropy. Here we use 2-transitivity of the Galois group again in conjugation with the fact that we chose $\sigma_0$ to be without fixed points to optimize this rate.

We multiply the bounds for all the places of $\mathbb{F}$ to have a bound on the adelic norm of $\Psi_\mathscr{C}({\delta_L}^{-1}\delta_R)$ in terms of the global discriminant and the exponential decay at the single place. By pushing the exponential bound far enough in comparison with the global discriminant term  we can decrease the bound on the norm below 1. But $\Psi_\mathscr{C}({\delta_L}^{-1}\delta_R)\in\mathbb{F}$, viz.\ it is a rational number, hence if it has norm smaller then 1 then it is equal to 0. Combining the previous statements we conclude that a tube which has been contracted long enough must include a single coset at most. This implies immediately a bound on the tube's measure.

\subsection{Organization of the paper.}

In \S \ref{sec:setting}, we review the concept of homogeneous toral sets in reductive linear algebraic groups and discuss their volume and discriminant data.

In \S \ref{sec:git} we construct the double quotient of a reductive linear group by a torus and study elementary properties of the quotient using geometric invariant theory.

In \S \ref{sec:dtq-sln} we present the canonical generators for the double quotient of $\mathbf{SL}_n$ and $\mathbf{PGL}_n$ by the maximal diagonal torus and analyze the algebraic relations between them.

In \S \ref{sec:dtq-pgl2} we review the relation between the canonical generators of $\mathbf{PGL}_2$ and the discriminant inner product used by Linnik and Skubenko.

In \S \ref{sec:dtq-csa} we use the results of \S \ref{sec:dtq-sln} to construct canonical generators for a general rational torus in a groups isogenous to the projective group of units of a central simple algebra. Notably, we study the Galois orbits of the values of canonical generators.

In \S \ref{sec:packets-csa} we prove arithmetic results relating the denominators of values of canonical generators for central simple algebras with the discriminant data of homogeneous toral sets. 

In \S \ref{sec:entropy} we show that values of canonical generators decay exponentially in contracting Bowen balls and prove the entropy lower bound.

In \S \ref{sec:rigidity} we combine the results of \S \ref{sec:entropy} with the measure rigidity theorems of \cite{EKL,EntropySArithmetic} to derive a statement about possible limit measures for a sequence of packets.

\subsection{Acknowledgments}
This paper is part of the author's PhD thesis conducted at the Hebrew University of Jerusalem under the guidance of Prof. E. Lindenstrauss, to whom I am grateful for introducing me to homogeneous dynamics and number theory, and for many helpful discussions and insights. I would like to express my gratitudes to Akshay Venkatesh for encouraging and valuable conversations. I thank Philippe Michel and Menny Aka for helpful comments on a preliminary version of this manuscript. I am grateful to the referees for suggestions that greatly improved the presentation of the paper.

Work on this project began during the MSRI program "Geometric and Arithmetic Aspects of Homogeneous Dynamics" in 2015. It is a pleasure to thank MSRI and the organizers for the program and for the institute's hospitality. 

Last but not least, I would like to thank my wife, Olga Kalantarov Hautin, for her unconditional support during the preparation of this paper.

The author has been supported by the European Research Council [AdG Grant 267259] throughout his PhD studies.

\section{The General Setting}\label{sec:setting}
\subsection{Algebraic Groups and Adeles}
Let $\mathbf{G}$ be a reductive linear algebraic group defined over a number field $\mathbb{F}$. Our objects of study are collections of periodic orbits of a maximal torus in $\mathbf{G}$ on an $S$-arithmetic quotient. Specifically, we study packets of periodic torus orbits as defined in \cite{ELMVCubic}, cf.\ \S \ref{sec:general-setting}.
 
For each place $u$ of $\mathbb{F}$ denote by $\Fu$ the completion of $\mathbb{F}$ at $u$.  We denote by $\mathcal{V}_\mathbb{F}$, $\mathcal{V}_{\mathbb{F},\infty}$ and $\mathcal{V}_{\mathbb{F}, f}$  the set of all places on $\mathbb{F}$,  the set of all archimedean places and the set of all non-archimedean places respectively. For $u$ non-archimedean we write $|\cdot|_u$ to be the canonical absolute value on $\Fu$, i.e.\ $|\varpi|_u=q^{-1}$ where $\varpi$ is a uniformizer for the ring of integers of $\Fu$ and $q$ is the size of the residue field. For $u$ archimedean $|\cdot|_u$ is the standard absolute value on $\mathbb{R}$ or $\mathbb{C}$. 

The choice for a norm in a complex place is slightly non-standard. For the product formula to hold we will consider conjugate complex embeddings as distinct places, although they define the same norm. The upside of this choice is that we will treat real and complex places uniformly.

Let $S$ be a finite set of places of $\mathbb{F}$ including all the infinite ones. Denote by $\mathbb{A}$ the adele ring of $\mathbb{F}$ and let
$\mathbb{A}^S\subset\mathbb{A}$ be the set of adeles having zero $v$-component for each $v\in S$. Set also $\mathbb{F}_S \coloneqq \prod_{u\in S} \mathbb{F}_u$, then $\mathbb{A}=\mathbb{F}_S\times \mathbb{A}^S$. 

We always treat $\mathbb{F}$ as diagonally embedded in $\mathbb{A}$, $\mathbb{A}^S$ and
$\mathbb{F}_S$. Similarly, for an affine algebraic variety $\mathbf{V}$ defined over $\mathbb{F}$, we treat $\mathbf{V}(\mathbb{F})$ as diagonally embedded in $\mathbf{V}(\mathbb{A})$ and $\mathbf{V}(\mathbb{F}_S)$. The spaces $\mathbf{V}(\mathbb{A})$ and $\mathbf{V}(\mathbb{F}_S)$ inherit a 
locally compact Hausdorff topology from the standard topology on $\mathbb{A}$ and $\mathbb{F}_S$.

For a field extension $\mathbb{M}/\mathbb{F}$ and an algebraic variety $\mathbf{V}$ defined over $\mathbb{F}$ we denote by $\mathbf{V}_\mathbb{M}$ the base change of $\mathbf{V}$ to $\mathbb{M}$. In scheme theoretic language $\mathbf{V}_\mathbb{M}\coloneqq \mathbf{V}\times_{\Spec\mathbb{F}} \Spec\mathbb{M}$. As all the main varieties we deal with are affine ones, the base change can be described completely by the transformation of the ring of regular functions $\mathbb{M}[\mathbf{V}_\mathbb{M}]\coloneqq \mathbb{F}[\mathbf{V}]\otimes_\mathbb{F}\mathbb{M}$.

\subsection{\texorpdfstring{$S$}{S}-Arithmetic and Adelic Quotients}\label{sec:adelic-quotient}
Denote $G_S\coloneqq \prod_{u\in S}\mathbf{G}(\mathbb{F}_u)$. We are interested in the locally homogeneous space $X_S\coloneqq\lfaktor{\Gamma}{G_S}$ for a congruence lattice $\Gamma\coloneqq\mathbf{G}(\mathbb{F})\cap K^S$ where $K^S$ is a compact-open subgroup in $\mathbf{G}(\mathbb{A}^S)$  and the intersection takes place in $\mathbf{G}(\mathbb{A}^S)$. 

This space can be naturally identified with the identity component of a factor space of the adelic quotient $X\coloneqq \lfaktor{\mathbf{G}(\mathbb{F})}{\mathbf{G}(\mathbb{A})}$.
We have the obvious right action of $G_S$ on $\faktor{X}{K^S}=\sdfaktor{\mathbf{G}(\mathbb{F})}{\mathbf{G}(\mathbb{A})}{K^S}$. By a theorem of Borel, the finiteness
of class number, the double quotient $\sdfaktor{\mathbf{G}(\mathbb{F})}{\mathbf{G}(\mathbb{A})}{K^SG_S}$ is finite \cite[Theorem 5.1]{BorelFinite} (see also \cite{ConradFinite}). In particular $X_S$ can be identified with the identity component of $\faktor{X}{K^S}$.

\subsection{Discriminant and Volume}\label{sec:discriminant-volume}
Recall from the introduction that a homogeneous toral set is a closed subset of $X$ of the form
\begin{equation*}
Y\coloneqq{\mathbf{G}(\mathbb{F})}{\mathbf{T}(\mathbb{A}) g_\mathbb{A}}\subset
\lfaktor{\mathbf{G}(\mathbb{F})}{\mathbf{G}(\mathbb{A})}
\end{equation*}
where $\mathbf{T}<\mathbf{G}$ a maximal torus defined and anisotropic over $\mathbb F$ and $g_{\mathbb A}\in\mathbf{G}(\mathbb A)$. If $g_\mathbb{A}=(g_u)_{u\in\mathcal{V}_\mathbb{F}}$, then for each place $u$ we can define over $\Fu$ the torus $\mathbf{H}_u={g_u}^{-1} \mathbf{T}_\Fu g_u$. The homogeneous toral set is invariant under the action of the geometric torus $\mathbf{H}_u(\Fu)$ for each $u$.

The homogeneous toral set can be identified with $\lfaktor{\mathbf{T}(\mathbb{F})}{\mathbf{T}(\mathbb{A})}$. Because $\mathbf{T}$ is anisotropic over $\mathbb{F}$, the subgroup $\mathbf{T}(\mathbb{F})$ is a lattice in $\mathbf{T}(\mathbb{A})$. In particular, $\lfaktor{\mathbf{T}(\mathbb{F})}{\mathbf{T}(\mathbb{A})}$ carries a unique probability measure invariant under $\mathbf{T}(\mathbb{A})$.
 
The volume and discriminant of a homogeneous toral set have been defined in \cite[\S\S 4.2--4.3]{ELMVCubic}. We reproduce these definitions here for convenience's sake.

\subsubsection{Volume}\label{sec:volume}
This definition depends on a fixed pre-compact identity neighborhood $\mathcal{K}\subset\mathbf{G}(\mathbb{A})$. Let $\mu_\mathbf{T}$ be the Haar measure on $\mathbf{T}$ such that the volume of
$\lfaktor{\mathbf{T}(\mathbb{F})}{\mathbf{T}(\mathbb{A})}$ is $1$, then one defines
\begin{equation*}
\vol(Y)=\mu_\mathbf{T}\left(t\in\mathbf{T}(\mathbb{A})\mid {g_\mathbb{A}}^{-1} t g_\mathbb{A} \in \mathcal{K} \right)^{-1}
\end{equation*}
The different values of the volume for different sets $\mathcal{K}$ are comparable up to constants depending on these sets only.
It is convenient to choose $\mathcal{K}=\prod_{u\in\mathcal{V}_\mathbb{F}}\mathcal{K}_u$ such that $\prod_{u\not\in S} \mathcal{K}_u=K^S$.

\subsubsection{Discriminant Data}\label{sec:disc}
The discriminant data is composed of a local discriminant $D_u$ for each place $u$ of $\mathbb{F}$. For $u$ non-archimedean $D_u\in q^\mathbb{Z}$ where $q$ is the size of the residue field of $\Fu$. For $u$ archimedean $D_u\in\mathbb{R}_{>0}$.

Let $\mathfrak{g}$ be the Lie algebra of $\mathbf{G}$, $\mathfrak{t}$ the Lie algebra of $\mathbf{T}$ and $\mathfrak{h}_u=\Ad({g_u}^{-1})\mathfrak{t}_{\Fu}$ the Lie algebra of $\mathbf{H}_u={g_u}^{-1}\mathbf{T}_{\Fu}g_u$. Notice that $\mathfrak{h}_u$ is only defined over $\Fu$.

Set $r$ to be the absolute rank of $\mathbf{G}$ and let $\mathbf{V}\coloneqq \left(\bigwedge^r\mathfrak{g}\right)^{\otimes2}$. This is an affine space defined over $\mathbb{F}$.
To define the local discriminant we need to make a choice of a good Goldman-Iwahori norm $\|\cdot\|_u$ on $\mathbf{V}(\Fu)$ for all places $u$. See the discussion at \cite[\S 7]{ELMVCubic} regarding norms. 

These norms must satisfy a compatibility condition which says that for almost all non-archimedean $u$ the unit ball of the norm coincides with the closure of a fixed $\mathcal{O}_\mathbb{F}$ lattice in $\mathbf{V}(\mathbb{F})$.

The subspace $\mathfrak{h}_u(\Fu)<\mathfrak{g}(\Fu)$ defines a point in $\mathbf{V}(\Fu)$ by a variant of the Pl{\"u}cker embedding  
\begin{equation*}
x_{\mathfrak{h}_u}\coloneqq
(f_1\wedge\ldots\wedge f_r)^{\otimes 2}\cdot \left|\det\left(B(f_i,f_j)\right)_{1\leq i,j \leq n}\right|_u^{-1}
\end{equation*}
Where $f_1,\ldots,f_r$ is a base for $\mathfrak{h}_u(\Fu)$ and $B$ is the Killing form.

The local discriminant is then defined to be $D_u=\| x_{\mathfrak{h}_u} \|_u$. The global discriminant is $D\coloneqq \prod_{u\in\mathcal{V}_\mathbb{F}} D_u$.
Notice that the local discriminant is an invariant of $\mathbf{H}_u$.

\begin{remark}
In many classical cases the archimedean discriminant is uniformly bounded for the homogeneous toral sets in question so it doesn't play a role. 

When studying periodic orbits of a fixed split torus $\mathbf{H}$ over $\mathbb{R}$ with $\mathbb{F}=\mathbb{Q}$ as in \cite{ELMVPGL2,ELMVPeriodic} we have $\mathbf{H}_\infty=\mathbf{H}$ and the archimedean discriminant is constant.

When studying the points on the $2$-sphere and modular analogues in \cite{LinnikBook,PointsOnSphere} all the archimedean tori $\mathbf{H}_\infty$ in the associated homogeneous toral sets are conjugate through elements of a compact group. Hence the archimedean discriminant is uniformly bounded.
\end{remark}

\section{Geometric Invariant Theory of a Double Quotient by a Torus}\label{sec:git}
In this part we construct and study the double quotient of a reductive linear algebraic $\mathbf{G}$ group by a torus $\mathbf{T}$. To any two orbits of  $\mathbf{T}$ on $\mathbf{G}$ we are able to attach a point on the double quotient space. 

The gist of the current section is a characterization of all the pairs of $\mathbf{T}$-orbits which have the trivial point attached to them in the double quotient space.

\subsection{Preliminaries}
Let $\mathbf{G}$ be a reductive linear algebraic group defined over a characteristic 0 field $\mathbb{F}$ and let $\Fbar$ be an algebraic closure. Take $\mathbf{T}<\mathbf{G}$ to be a non-trivial torus defined over $\mathbb{F}$. We always denote by $e$ the identity element of groups when written in multiplicative form.

Recall that for a torus $\mathbf{H}$ defined over $\mathbb{F}$ we can form the character group $X^\bullet(\mathbf{H}) = \Hom_{\Fbar}(\mathbf{H},\mathbb{G}_m)$ 
and the cocharacter group
$X_\bullet(\mathbf{H})=\Hom_{\Fbar}(\mathbb{G}_m, \mathbf{H})$. Those are free abelian groups whose rank is equal to the absolute rank of 
$\mathbf{H}$. The character and cocharacter groups
come with a natural perfect pairing
$\left<,\right>\colon X^\bullet(\mathbf{H})\times X_\bullet(\mathbf{H})\to\End(\mathbb{G}_m)\simeq \mathbb{Z}$  defined by the composition of a character
with a cocharacter.  We have an action of the absolute Galois group $\Gal(\Fbar/\mathbb{F})$ on $X^\bullet(H)$, $X_\bullet(H)$ making them Galois modules. In particular ${X^\bullet(H)}^{\Gal(\Fbar/\mathbb{F})}$, ${X_\bullet(H)}^{\Gal(\Fbar/\mathbb{F})}$ are the groups of characters, respectively cocharacters,  defined over $\mathbb{F}$.

\subsection{Synopsis}
For convenience's sake we briefly summarize the results of this chapter. The main tools in proving the following statements are the Geometric
invariant theory of Mumford \cite{GIT} and Kempf's work on the numerical criteria for stability of orbits \cite{instability}.

\begin{enumerate}
\item There exists an affine algebraic variety defined over $\mathbb{F}$, $\dfaktor{\mathbf{T}}{\mathbf{G}}{\mathbf{T}}$, and an $\mathbb{F}$-morphism 
$\pi\colon \mathbf{G} \to \dfaktor{\mathbf{T}}{\mathbf{G}}{\mathbf{T}}$ which is equivariant with respect to the left and right actions of
$\mathbf{T}$ on $\mathbf{G}$.

\item The variety $\dfaktor{\mathbf{T}}{\mathbf{G}}{\mathbf{T}}$ can be realized explicitly by choosing a generating set $\Psi_1,\ldots,
\Psi_m$ in  the ring of regular function on $\mathbf{G}$ which are invariant under both the left and the right action of $\mathbf{T}$.  Using the
generating functions we define a closed morphism $\Psi=\left(\Psi_1,\ldots,\Psi_m\right) \colon  \mathbf{G}  \to \mathbb{A}^m$. The image of $\Psi$ is isomorphic to $\dfaktor{\mathbf{T}}{\mathbf{G}}{\mathbf{T}}$. This isomorphism conjugates $\Psi$ to $\pi$.

\item For each $\lambda\in\mathbf{G}(\mathbb{F})$ such that $\pi(\lambda)=\pi(e)$ we have $\lambda\in\mathbf{T}(\mathbb{F})$.
\end{enumerate}

The last result is used heavily in what follows. We are able to prove that if two torus orbits $\delta_L\mathbf{T}(\mathbb{A})$ and $\delta_R\mathbf{T}(\mathbb{A})$ with $\delta_L,\delta_R\in\mathbf{G}(\mathbb{F})$ from the same homogeneous toral set come too close to each other then $\pi({\delta_L}^{-1}\delta_R) =\pi(e)$. The result quoted above insures us that in this case they must actually be the same orbit. This is exactly the well-separateness property required to prove our entropy inequality.

\subsection{Construction and Zariski Closed Orbits}
To define the affine  algebraic variety $\dfaktor{\mathbf{T}}{\mathbf{G}}{\mathbf{T}}$ we start by looking at the action of the algebraic torus $\mathbf{T} \times \mathbf{T}$ on $\mathbf{G}$ defined by $(t_L,t_R).g=t_L g {t_R}^{-1}$. This is an algebraic action of a reductive algebraic group $\mathbf{T}\times\mathbf{T}$ on an affine algebraic variety $\mathbf{G}$. Let $\tensor[^{\mathbf{T}}]{\mathbb{F} [\mathbf{G}]}{^{\mathbf{T}}}$ be the ring of regular function on $\mathbf{G}$ which are invariant under both the left and the right action of $\mathbf{T}$. By a result of Hilbert this is a finitely generated algebra over $\mathbb{F}$ as $\mathbf{T} \times \mathbf{T}$ is reductive.

The GIT quotient is defined by $\dfaktor{\mathbf{T}}{\mathbf{G}}{\mathbf{T}} \coloneqq \Spec{\tensor[^{\mathbf{T}}]{\mathbb{F} [\mathbf{G}]}{^{\mathbf{T}}} }$. This is an affine algebraic variety defined over $\mathbb{F}$ because $\tensor[^{\mathbf{T}}]{\mathbb{F} [\mathbf{G}]}{^{\mathbf{T}}}$ is finitely generated. This variety comes with a natural $\mathbf{T}\times\mathbf{T}$-equivariant $\mathbb{F}$-morphism $\pi\colon \mathbf{G} \to \dfaktor{\mathbf{T}}{\mathbf{G}}{\mathbf{T}}$ which is induced by the inclusion map $\tensor[^{\mathbf{T}}]{\mathbb{F} [\mathbf{G}]}{^{\mathbf{T}}} \hookrightarrow \mathbb{F} [\mathbf{G}]$.

The following theorem is a classical result of Mumford's theory for GIT quotients in the affine case.
\begin{theorem}[\cite{GIT}]
\label{thm:quotient}
$\dfaktor{\mathbf{T}}{\mathbf{G}}{\mathbf{T}}$ is a universal categorical quotient. In particular the following holds
\begin{enumerate}
\item $\pi$ is surjective.
\item Each fiber of $\pi$ contains a unique Zariski closed orbit.
\item For a field extension $\mathbb{L}/\mathbb{F}$ denote by $\mathbf{G}_{\mathbb{L}}$ and $\left(\dfaktor{\mathbf{T}}{\mathbf{G}}{\mathbf{T}}\right)_{\mathbb{L}}$ the corresponding 
varieties after base change. Then $\mathbf{G}_\mathbb{L} \to \left(\dfaktor{\mathbf{T}}{\mathbf{G}}{\mathbf{T}}\right)_\mathbb{L}$ is a categorical quotient for the
$\mathbf{T}_{\mathbb{L}}\times\mathbf{T}_{\mathbb{L}}$ action.
\end{enumerate}
\end{theorem}
\begin{proof}
See \cite[\S 1.2]{GIT} and \cite[\S 4.4]{invariant}.
\end{proof}

A significant complication is that the $\mathbf{T}\times\mathbf{T}$ action on $\mathbf{G}$ necessarily has non Zariski closed orbits of $\Fbar$-points. As a consequence this is not a geometric quotient and such a quotient does not exist, at least not for the whole variety $\mathbf{G}$.

A main tool for studying the Zariski closure of orbits are cocharacters of the acting group $\mathbf{T}\times\mathbf{T}$, also called multiplicative 1-parameter subgroups. Those are morphisms $\mathbb{G}_m\to\mathbf{T}\times\mathbf{T}$. We cite the following important result of Kempf building upon the work of Mumford and Hilbert (see also \cite[Chapter 2]{GIT}).

\begin{definition}
Fix a separated scheme $\mathbf{X}$.
Let $r\colon \mathbb{G}_m \to \mathbf{X}$ be a morphism. We say that $\lim_{\tau\to 0} r(\tau)$ exists and equals to a point $x$ of $\mathbf{X}$ if $r$ can be extended to a morphism $\mathbb{A}^1\to \mathbf{X}$ and $x$ is the image of the origin.
\end{definition}

\begin{proposition}\cite[Corollary 4.3]{instability}
\label{prop:instability}
Let $\mathbf{H}$ be a connected reductive algebraic group acting on an affine scheme $\mathbf{X}$. Let $\mathbf{S}$ be a closed $\mathbf{H}$-invariant subscheme
of $\mathbf{X}$. Let x be an $\mathbb{F}$-point of $\mathbf{X}$ and assume that $\mathbf{S}$ meets the Zariski closure of the orbit $\mathbf{H}x$. Then there exists a cocharacter $\lambda$ of $\mathbf{H}$ defined over $\mathbb{F}$ such that $\mathbf{S}\cap\overline{\mathbf{H}x}^{\mathrm{Zar}}$ contains the $\mathbb{F}$-point $\lim_{\tau\to 0} \lambda(\tau)x$.
\end{proposition}
\begin{corollary}
\label{cor:no-cocharacters}
\cite[Remark, p.314]{instability}
If $\mathbf{H}$ has no non-trivial cocharacters defined over $\mathbb{F}$ then the $\mathbf{H}$-orbit of any $\mathbb{F}$-point of $X$ is Zariski closed.
\end{corollary}

Notice that although the notion of limit in Proposition \ref{prop:instability} is defined algebraically, because of \cite[Lemma 1.2]{instability}, for an affine algebraic variety $X$ and a local field $\mathbb{F}$ the algebraic definition of the limit coincides with limit in the Hausdorff topology induced from the topology on the field.

\begin{example}
\label{exmpl:non-closed}
Let $\mathbf{B}=\mathbf{T}\ltimes\mathbf{U}$ be a Borel subgroup defined over $\Fbar$ and containing $\mathbf{T}$ and let $n\in \Nrml_\mathbf{G}(\mathbf{T}) (\Fbar)$. The double torus orbit of $x=n u$ for any $u \in \mathbf{U}(\Fbar)$, $u\neq e$, is not Zariski closed.
\end{example}
\begin{proof}
We will show that $n$ is in the Zariski closure of $\mathscr{O} \coloneqq \mathbf{T}(\Fbar) x \mathbf{T}(\Fbar)$ but $n \not\in \mathscr{O}$ because $u\neq e$. This will prove $\overline{\mathscr{O}}^{\mathrm{Zar}}\neq {\mathscr{O}}$. 

Recall that for any Borel subgroup $\mathbf{B}'$ containing the torus $\mathbf{T}$ there is always a cocharacter of $\mathbf{T}$, $\lambda\in X_\bullet(\mathbf{T})$, such that for any character $\alpha\in X^\bullet(\mathbf{T})$ we have $\left<\lambda, \alpha\right> > 0$ if and only if $\alpha$ is a positive
character relative to $\mathbf{B}'$. Let $\lambda\in X_\bullet(\mathbf{T})$ be such cocharacter relative to $\mathbf{B}^-$, i.e.\ $\left<\lambda,\alpha\right> < 0$ if and only if $\alpha$ is a positive character relative $\mathbf{B}$. Now $(n\lambda n^{-1}, \lambda)$ is cocharacter of $\mathbf{T}\times\mathbf{T}$ and
\begin{equation}\label{eq:lambda}
\lim_{\tau\to0} \left(n \lambda(\tau) n^{-1}\right) x {\lambda(\tau)}^{-1}=  n \lim_{\tau\to0} \left[\lambda(\tau) u {\lambda(\tau)}^{-1}\right]
\end{equation}

We claim that $\lim_{\tau\to0} \left[\lambda(\tau) u {\lambda(\tau)}^{-1}\right]=e$ as required. To see that use the fact that $\mathbf{U}$ is
generated by the the 1-parameter unipotent subgroups $\mathbf{U}_\alpha$ for positive characters $\alpha$ and decompose $u$ into a product
of 1-parameter unipotents. For each 1-parameter unipotent $\mathbf{U}_\alpha$ there is an isomorphism $\epsilon_\alpha:\mathbb{G}_a\to
\mathbf{U}_\alpha$ such that $t \epsilon_\alpha(s) t^{-1}=\epsilon_\alpha(\alpha(t) s)$ for all $t\in\mathbf{T}(\Fbar)$ and all $s\in\mathbb{G}_a
(\Fbar)$. In particular for the cocharacter $\lambda$ it holds $\lambda(\tau) \epsilon_\alpha(s) {\lambda(\tau)}^{-1}
=\epsilon_\alpha\left(t^{\left<\lambda,\alpha\right>}s\right)$ which tends to $0$ for all $s$ if and only if $\left<\lambda, \alpha\right> < 0$. Hence the claim follows by choosing $\lambda$ as in \eqref{eq:lambda}.
\end{proof}
 
In the following we are mainly interested in double torus orbits of rational points $g\in \mathbf{G}(\mathbb{F})$ for a torus $\mathbf{T}$ anisotropic $\mathbb{F}$, those points always have Zariski closed double torus orbit.
\begin{proposition}
\label{prop:closed-anisotropic}
Assume $\mathbf{T}$ is anisotropic over $\mathbb{F}$ then the double torus orbit of any $g\in\mathbf{G}(\mathbb{F})$ is Zariski closed.
\end{proposition}
\begin{proof}
By Corollary \ref{cor:no-cocharacters} it is enough to show that $\mathbf{T}\times\mathbf{T}$ has no cocharacter over $\mathbb{F}$, viz.\ that the $\mathbb{F}$-cocharacter group
$X_\bullet(\mathbf{T}\times \mathbf{T})^{\Gal(\Fbar/\mathbb{F})}$ is trivial. This is equivalent to
$X_\bullet(\mathbf{T})^{\Gal(\Fbar/\mathbb{F})}$ being trivial. But for a torus the pairing $\left<,\right>$ between characters and cocharacters is perfect and
Galois equivariant, hence this is the same as $\mathbf{T}$ being anisotropic over $\mathbb{F}$.
\end{proof} 

\subsection{Projections to The Identity}
Our argument is going to use only very rudimentary properties of the general theory of double torus quotients. We need only to understand which rational points have the same image in $\dfaktor{\mathbf{T}}{\mathbf{G}}{\mathbf{T}}$ as the identity. 

Much stronger results can be proved using Galois cohomology about the rational points in the fiber of $\pi$ over any rational point in the double torus quotient. Unfortunately, difficulties in subsequent parts of the argument, specifically in counting integral points in $\dfaktor{\mathbf{T}}{\mathbf{G}}{\mathbf{T}}$, stop us from exploiting such results.

\begin{proposition}\label{prop:identity-fiber}
Assume that $\mathbf{T}$ is anisotropic over $\mathbb{F}$. Then for any $g\in\mathbf{G}(\mathbb{F})$ such that $\pi(g)=\pi(e)$ we have that $g\in\mathbf{T}(\mathbb{F})$.
\end{proposition}
\begin{proof}
Proposition \ref{prop:closed-anisotropic} implies that both $\mathbf{T}g\mathbf{T}$ and $\mathbf{T}e\mathbf{T}$ are Zariski closed. The latter orbit is actually equal to $\mathbf{T}$.

Theorem \ref{thm:quotient} says that there is a unique Zariski closed orbit over $\pi(e)$, hence $\mathbf{T}=\mathbf{T}e\mathbf{T}=\mathbf{T}g\mathbf{T}$. In particular, over an algebraically closed field $\Fbar$ we have $\mathbf{T}(\Fbar)=\mathbf{T}(\Fbar)g\mathbf{T}(\Fbar)$. That is there exist $t_0,t_L,t_R\in\mathbf{T}(\Fbar)$ such that
\begin{equation*}
t_0=t_L g {t_R}^{-1} \implies g= {t_L}^{-1} t_0 t_R
\end{equation*}
Hence $g\in\mathbf{T}(\Fbar)\cap\mathbf{G}(\mathbb{F})=\mathbf{T}(\mathbb{F})$.
\end{proof}

\section{Double Torus Quotient for \texorpdfstring{$\mathbf{SL}_n$}{SLn} and \texorpdfstring{$\mathbf{PGL}_n$}{PGLn}}\label{sec:dtq-sln}
Let $\mathbb{F}$ be a number field. We treat the matrix algebra $\mathbf{M}_n$ and the associated linear algebraic groups $\mathbf{GL}_n$, $\mathbf{PGL}_n$ and $\mathbf{SL}_n$ as defined over $\mathbb{F}$. We denote by $h\colon \mathbf{SL}_n\to\mathbf{PGL}_n$ the standard isogeny between these groups.

\subsection{Canonical Generators for Standard Tori}
Let $\mathbf{Diag}<\mathbf{M}_n$ be the commutative subalgebra of diagonal matrices and denote $\mathbf{A}=\mathbf{GL}_1(\textbf{Diag})$. Define $\mathbf{SA}\coloneqq\mathbf{SL}_1(\mathbf{Diag})$ to be the subgroup of diagonal matrices of determinant 1 and let $\mathbf{PA}\coloneqq\mathbf{PGL}_1(\mathbf{Diag})$ to be the full diagonal torus of $\mathbf{PGL}_n$. The groups $\mathbf{PA}$ and $\mathbf{SA}$ are maximal tori in the corresponding groups split over $\mathbb{F}$ .

We build a specific set of generators for $\tensor[^{\mathbf{PA}}]{\mathbb{F} [\mathbf{PGL}_n]}{^{\mathbf{PA}}}$ and
$\tensor[^{\mathbf{SA}}]{\mathbb{F} [\mathbf{SL}_n]}{^{\mathbf{SA}}}$ which provides us with a description in coordinates of both $\dfaktor{\mathbf{PA}}{\mathbf{PGL}_n}{\mathbf{PA}}$ and $\dfaktor{\mathbf{SA}}{\mathbf{SL}_n}{\mathbf{SA}}$.

\subsubsection{Double Torus Quotient for the Variety of Matrices}
Let $\mathbf{M}_n$ be the affine algebraic variety of $n\times n$ matrices, with $x_{i,j}\in\mathbb{F}
[\mathbf{M}_n]$, $1\leq i,j \leq n$, the function attaching to a matrix its $(i,j)$ entry. In this form
the coordinate ring of $\mathbf{M}_n$ is the polynomial algebra 
$R\coloneqq \mathbb{F}[x_{i,j}]_{1\leq i,j \leq n}$. The algebraic group $\mathbf{SA}\times\mathbf{SA}$ acts on $\mathbf{M}_n$.  The the first copy of $\mathbf{SA}$ in the product acts by the standard left action and the second copy by the standard right action. The closed embedding
$\mathbf{SL}_n\to\mathbf{M}_n$ is obviously equivariant under this action.

Our first step is to describe $R_0\coloneqq \tensor[^{\mathbf{SA}}]{R}{^{\mathbf{SA}}}$, the ring of regular function on
$\dfaktor{\mathbf{SA}}{\mathbf{M}_n}{\mathbf{SA}}$. This ring is generated by those polynomials $P\in R$ which are
invariant under the $\mathbf{SA}\times\mathbf{SA}$ action.

All the following calculations are executed over a characteristic $0$ algebraically closed field.
When an $n\times n$ matrix is multiplied on the left by $\lambda=
\diag{(\lambda_1,\ldots,\lambda_n)}$ and on the right by the inverse of $\mu=\diag{(\mu_1,\ldots,\mu_n)}$ its $(i,j)$ coordinate is multiplied by $\faktor{\lambda_i}{\mu_j}$, i.e. $(\lambda, \mu). x_{i,j}=
\faktor{\lambda_i}{\mu_j} x_{i,j}$. 

 A polynomial in $P\in R$ is of the form
\begin{equation*}
P=
\sum_M {b_{M} \prod_{1\leq i,j\leq n} x_{i,j}^{M_{i,j}}}
\end{equation*}
Where $M$ runs over a finite set of monomials, and $M_{i,j}$ is the power of $x_{i,j}$ in the $M$-monomial. 
We identify $M=\left(M_{i,j}\right)_{1\leq i,j \leq n}$ with a matrix of non-negative integers. Denote
\begin{itemize}
\item $M_i\coloneqq \sum_{1\leq j \leq n}{M_{i,j}}$ - the sum of the columns in row $i$.
\item $M^j\coloneqq \sum_{1\leq i \leq n}{M_{i,j}}$ - the sum of the rows in column $j$.
\end{itemize}

We see that the double torus action on $P$ takes the form
\begin{equation*}
(\lambda,\mu).P=
\sum_M  {b_{M}
\left(\prod_{1\leq i \leq n}{\lambda_i}^{M_i}\right) 
\left(\prod_{1\leq j \leq n}{\mu_j}^{M^j}\right)^{-1}
 \prod_{1\leq i,j\leq n} x_{i,j}^{M_{i,j}}}
\end{equation*}

In particular the $\mathbf{SA}\times\mathbf{SA}$ actions sends monomials to monomials.
If the polynomial $P$ is invariant, then because the polynomial algebra is free and the action conserves monomials we must have that each monomial appearing in $P$ is invariant.
This implies that $R_0$ is generated by invariant monomials. For such an invariant monomial with power matrix $M$ we must have for all diagonal matrices $\lambda,\mu$ with determinant $1$ 
\begin{equation}\label{eq:lambda-mu=1}
\left(\prod_{1\leq i \leq n}{\lambda_i}^{M_i}\right) 
\left(\prod_{1\leq j \leq n}{\mu_j}^{M^j}\right)^{-1}=1
\end{equation}
The only relation between the $\lambda_i$'s is $\prod_{1\leq i \leq n} \lambda_i =1$ and the same hold for the $\mu_j$'s. Equality \eqref{eq:lambda-mu=1} can hold if and only if for all $1\leq i_1,i_2 \leq n$: $M_{i_1}=M_{i_2}$ and for all $1\leq j_1,j_2 \leq n$: $M^{j_1}=M^{j_2}$. A matrix $M$ of non-negative integers
all whose rows sum to the same value and all whose columns sum to the same value is called a semi-magic square. 
Notice that because $\sum_{1\leq i \leq n} M_i =
\sum_{1\leq j \leq n} M^j$ the value must be the same for the rows and the columns.

It is a classical result of D. K{\"o}nig that each semi-magic square is a sum of permutation matrices \cite{Konig} 
(see also the exposition \cite{marriage}).
An immediate consequence is that $R_0$ is generated by the monomials corresponding to permutation matrices. The discussion
above sums to the proof of the following.
\begin{proposition}\label{prop:generators-matrix}
The ring $R_0=\tensor[^{\mathbf{SA}}]{\mathbb{F} [\mathbf{M}_n]}{^{\mathbf{SA}}}$ is generated by the monomials
\begin{equation*}
\Psi^0_\sigma \coloneqq \sign{\sigma}\prod_{1\leq i \leq n} x_{\sigma(i),i}
\end{equation*}
for $\sigma\in S_n$.
\end{proposition}
\begin{remark}
There are many relations between these generators. In particular notice that we have $n!$ generators, but $\dfaktor{\mathbf{SA}}{\mathbf{M}_n}{\mathbf{SA}}$ has a Zariski open subset of dimension $n^2-2(n-1)$. This is the set of stable points - the points that have a Zariski closed orbit and 
trivial stabilizer. Nevertheless, $R_0$ is finitely presented and we can describe all the relations explicitly.

\end{remark}

\subsubsection{Passing from \texorpdfstring{$\mathbf{M}_n$}{Mn} to \texorpdfstring{$\mathbf{SL}_n$}{SLn}}
We recall the definition of the determinant using the Leibniz formula
\begin{equation*}
\mathrm{det}
\coloneqq\sum_{\sigma \in S_n}{\sign{\sigma}\cdot \prod_{i=1}^n{x_{\sigma(i),i}}}
=\sum_{\sigma\in S_n}{ \Psi^0_\sigma } \in R
\end{equation*}
This is a homogeneous polynomial of degree $n$. Denote by $I\coloneqq \left<\mathrm{det}-1\right>$ the principal ideal in $R$ generated by 
$\mathrm{det}-1$. The definition of the special linear group implies  $\mathbb{F} [\mathbf{SL}_n]=\faktor{R}{I}$. 

We are interested in the ring of invariants $\left(\faktor{R}{I}\right)_0\coloneqq \tensor[^{\mathbf{SA}}]{\left(\faktor{R}{I}\right)}{^{\mathbf{SA}}}$. Define 
$I_0\coloneqq I \cap R_0$.
We have an injective homomorphism of $\mathbb{F}$-algebras $\faktor{R_0}{I_0}\hookrightarrow \left(\faktor{R}{I}\right)_0$. The following lemma which appears in the work of Nagata \cite[Lemma 5.1.A]{Nagata} implies that this homomorphism is actually an isomorphism (see also \cite[Lemma 3.5]{Dolgachev}).
\begin{lemma}[Nagata]
Let a linearly reductive algebraic group $\mathbf{H}$ act on the $\mathbb{F}$-algebras $S$ and $S'$ for a field $\mathbb{F}$. Assume everything is defined over $\mathbb{F}$. 
If $\phi\colon S\to S'$ is an 
$\mathbf{H}$-equivariant $\mathbb{F}$-homomorphism \textbf{onto} $S'$, then induced homomorphism $S^\mathbf{H}\to {S'}^\mathbf{H}$ is onto.
\end{lemma}
\begin{corollary}\label{cor:mn-sln}
The homomorphism $\faktor{R_0}{I_0}\hookrightarrow \left(\faktor{R}{I}\right)_0$ is an isomorphism. This allows
us to identify hereon $\faktor{R_0}{I_0} \simeq \left(\faktor{R}{I}\right)_0$. 

In particular $\left(\faktor{R}{I}\right)_0$ is generated by $\left\{ \Psi^0_\sigma \mid \sigma \in S_n\right\}$ and we can describe
$\dfaktor{\mathbf{SA}}{\mathbf{SL}_n}{\mathbf{SA}}$ in coordinate form by the image of $\Psi\colon \mathbf{SL}_n \to \mathbb{A}^{(n!)}$ defined
by $\Psi=\left(\Psi^0_\sigma\right)_{\sigma\in S_n}$ where we have fixed an arbitrary order on $S_n$.
\end{corollary}
\begin{proof}
Applying the lemma to the homomorphism $R\to \faktor{R}{I}$ shows that the homomorphism $R_0\to \left(\faktor{R}{I}\right)_0$ is onto, but this homomorphism factors through $\faktor{R_0}{I_0}$. This shows that $\faktor{R_0}{I_0}\to \left(\faktor{R}{I}\right)_0$ is bijective hence an isomorphism.
\end{proof}

\subsubsection{Passing from \texorpdfstring{$\mathbf{SL}_n$}{SLn} to \texorpdfstring{$\mathbf{PGL}_n$}{PGLn}}
\label{sec:pgl-generators}
The isogeny $h\colon \mathbf{SL}_n \to \mathbf{PGL}_n$ restricts to an isogeny of maximal tori 
\begin{equation*}
h\restriction_{\mathbf{SA}}\colon \mathbf{SA} \to \mathbf{PA}
\end{equation*}
Using $h$ we can let $\mathbf{SA}\times\mathbf{SA}$ act on the affine algebraic variety $\mathbf{PGL}_n$ to form the double quotient $\dfaktor{\mathbf{SA}}{\mathbf{PGL}_n}{\mathbf{SA}}$. Surjectivity of $h\restriction_{\mathbf{SA}}$ implies that 
\begin{equation*}
\dfaktor{\mathbf{SA}}{\mathbf{PGL}_n}{\mathbf{SA}} = \dfaktor{\mathbf{PA}}{\mathbf{PGL}_n}{\mathbf{PA}}
\end{equation*}
From now on we identify those spaces.

Next we use the surjectivity of the full isogeny $h$. The morphism $h$ induces a surjective morphism $\dfaktor{\mathbf{SA}}{\mathbf{SL}_n}{\mathbf{SA}} \to \dfaktor{\mathbf{PA}}{\mathbf{PGL}_n}{\mathbf{PA}}$. This morphism intertwines the $\mathbf{SA}$ actions and it is the unique morphism having this property, as follows from the categorical properties of the quotient.

Because of the duality between affine algebraic varieties over $\mathbb{F}$ and finitely generated reduced  algebras over $\mathbb{F}$ we have an \emph{injective}
algebra homomorphism 
\begin{equation*}
\tensor[^{\mathbf{PA}}]{\mathbb{F} [\mathbf{PGL}_n]}{^{\mathbf{PA}}} \hookrightarrow \tensor[^{\mathbf{SA}}]{\mathbb{F} [\mathbf{SL}_n]}{^{\mathbf{SA}}}
\end{equation*}

In particular, if one can demonstrate a set of elements in $\tensor[^{\mathbf{PA}}]{\mathbb{F} [\mathbf{PGL}_n]}{^{\mathbf{PA}}}$ whose image in $\tensor[^{\mathbf{SA}}]{\mathbb{F} [\mathbf{SL}_n]}{^{\mathbf{SA}}}$ generates the latter ring, then those elements must generate the former ring. We will find such regular functions on $\mathbf{PGL}_n$ that are equal to the $\Psi_\sigma^0$'s when restricted to the image of $\mathbf{SL}_n$.

To represent $\mathbf{PGL}_n$ as an affine algebraic variety we use the adjoint representation which is a closed immersion $\mathbf{PGL}_n\to\mathbf{GL}_n(\mathbf{M}_n)$. Notably, every regular function on $\mathbf{GL}_n(\mathbf{M}_n)$ induces a regular function on $\mathbf{PGL}_n$. Hence, all we need to do is to exhibit regular functions on  $\mathbf{GL}_n(\mathbf{M}_n)$ such that when restricted to the image of $\mathbf{SL}_n$ through the adjoint representation they would be equal to the $\Psi_\sigma^0$'s.

This is easy enough to do if we rewrite $\Psi_\sigma$ in a slightly more intrinsic way. For $1\leq i \leq n$ let $e^0_i=\diag(0,\ldots,0,1,0,\ldots,0)$ be the diagonal matrix with a single $1$ entry in the $(i,i)$ place. Then for $g\in\mathbf{SL}_n(\Fbar)$
\begin{equation*}
\Psi^0_\sigma(g)=\det\left(\sum_{i=1}^n e^0_{\sigma(i)}\,g\,e^0_i\right)=\det\left(\sum_{i=1}^n e^0_{\sigma(i)}\,{\Ad}_g(e^0_i)\right)
\end{equation*} 
The latter expression can be easily extended to a regular function on the variety $\mathbf{GL}_n(\mathbf{M}_n)$ by considering a general linear transformation on $\mathbf{M}_n$ instead of ${\Ad}_g$. Those regular functions would necessarily generate $\tensor[^{\mathbf{PA}}]{\mathbb{F} [\mathbf{PGL}_n]}{^{\mathbf{PA}}}$ as required. We sum up these results in the following proposition.

\begin{proposition}\label{prop:pgl-generators}
The ring $\tensor[^{\mathbf{PA}}]{\mathbb{F} [\mathbf{PGL}_n]}{^{\mathbf{PA}}}$ is generated by the following regular functions, which by abuse of notation we also denote by $\Psi^0_\sigma$
\begin{equation}\label{eq:psi0-def}
\Psi^0_\sigma(g)=\det\left(\sum_{i=1}^n e^0_{\sigma(i)}\,{\Ad}_g(e^0_i)\right)=\det\left(\sum_{i=1}^n e^0_{\sigma(i)}\,g\,e^0_i\right)\cdot\left(\det g\right)^{-1}
\end{equation}
For any $[g]\in\mathbf{PGL}_n(\mathbb{\Fbar})$ and all $\sigma\in S_n$. In expression \eqref{eq:psi0-def} we have treated $g$ as a matrix representing a point in $\mathbf{PGL}_n$.

Notice that this definition of $\Psi^0_\sigma$ for $\mathbf{PGL}_n$ recovers the original functions $\Psi^0_\sigma$ for $\mathbf{SL}_n$ through the homomorphism $\mathbb{F}[\mathbf{PGL}_n]\to\mathbb{F}[\mathbf{SL}_n]$.

Moreover, the homomorphism $\tensor[^{\mathbf{PA}}]{\mathbb{F} [\mathbf{PGL}_n]}{^{\mathbf{PA}}} \hookrightarrow \tensor[^{\mathbf{SA}}]{\mathbb{F} [\mathbf{SL}_n]}{^{\mathbf{SA}}}$ is an isomorphism of $\mathbb{F}$-algebras, i.e.\ $\dfaktor{\mathbf{SA}}{\mathbf{SL}_n}{\mathbf{SA}} \cong \dfaktor{\mathbf{PA}}{\mathbf{PGL}_n}{\mathbf{PA}}$. In particular, the relations between the $\Psi^0_\sigma$'s for $\mathbf{PGL}_n$ are the same as for $\mathbf{SL}_n$.
\end{proposition}

From now on, unless stated otherwise, we write the $(\det g)^{-1}$ factor and its analogues even when treating $\mathbf{SL}_n$. In that case it should be viewed as the constant function $1$. In any case, we can always consider the functions $\Psi^0_\sigma$ as functions on $\mathbf{GL}_n$.

\subsection{Relations between the Canonical Generators}
We present two different results regarding the relations between the canonical generators. The first one is specially adapted to our needs and is restricted to the question of when one canonical generator being equal to zero implies that another generator is also zero.
In the second part we describe explicitly a complete set of relations between the generators.

\begin{proposition}\label{prop:relations-zero}
To each $\sigma\in S_n$ we attach the \emph{roots set}. 
\begin{equation*}
R_\sigma\coloneqq\left\{ \left(\sigma(i),i\right) \mid 1\leq i \leq n ,\, \sigma(i)\neq i\right\}
\end{equation*}
We say that a subset $\mathscr{C}\subseteq S_n$ \emph{has a complete set of roots} if 
\begin{equation*}
\bigcup_{\sigma\in\mathscr{C}} R_\sigma=\left\{ (j,i) \mid 1\leq j,i \leq n,\, i\neq j \right\}
\end{equation*}
 Notice that if $F<S_n$ is a 2-transitive subgroup then for every $\sigma\neq \mathrm{id}$ the subset $\left\{ \tau \sigma \tau^{-1} \mid \tau\in F\right\}$ has a complete set of roots.

If for $g\in\mathbf{GL}_n(\mathbb{F})$ there exists $\sigma\in S_n$ such that $\sigma$ has no fixed points and $\Psi^0_\sigma(g)=0$ then in every subset $\mathscr{C}\subseteq S_n$  which has a complete set of roots one can find $\tau\in\mathscr{C}$ such that $\Psi^0_\tau(g)=0$.
\end{proposition}

The connection  between a root set to the roots of a maximal torus in an algebraic group is presented in \S \ref{sec:entropy}. 

\begin{proof}
We write $g$ in matrix form $(g_{i,j})_{1\leq i,j \leq n}$, viz.\ $x_{i,j}(g)=g_{i,j}$.
Recall that  $\Psi^0_\sigma(g) \coloneqq \sign{\sigma}\prod_{1\leq i \leq n} g_{\sigma(i),i} \cdot (\det g)^{-1}$. Hence if $\Psi^0_\sigma(g)=0$ then there exists $1\leq i_0 \leq n$ such that $g_{\sigma(i_0),i_0}=0$ and $\sigma(i_0)\neq i_0$.

For every $\tau\in S_n$ such that $\tau(i_0)=\sigma(i_0)$ we have $\Psi^0_\tau(g)=0$. In every subset $\mathscr{C}\subseteq S_n$ that has a complete set of roots one can find such $\tau\in \mathscr{C}$ and the proposition follows.
\end{proof}

\subsubsection{Complete Description of the Relations}
The results in the second part of this section are presented mainly for the sake of completeness
as they do not play a direct role in what follows. Nevertheless, the relations between the generators are a key part of our argument albeit only in the form of Proposition \ref{prop:relations-zero} and not in the explicit form presented here. 

\paragraph{Renormalization of the canonical generators.}
We wish first to describe the relations between the canonical generators in $R_0=\tensor[^{\mathbf{SA}}]{\mathbb{F} [\mathbf{M}_n]}{^{\mathbf{SA}}}$. Here we work with a different normalization for generators of $R_0$
\begin{equation}\label{eq:psi1}
\Psi^1_\sigma=\prod_{i=1}^n x_{i,\sigma(i)}=\prod_{1\leq i,j \leq n} x_{i,j}^{P^\sigma_{i,j}}
\end{equation}
Where $P^\sigma\in\mathbf{M}_n(\mathbb{Z})$ is the standard permutation matrix associated to $\sigma\in S_n$.
Notice the lack of the twist by the determinant as currently we just study $R_0$.

We can easily recover the original canonical generators from Proposition \ref{prop:generators-matrix} using
$\Psi^0_\sigma=\sign\sigma\, \Psi^1_{\sigma^{-1}}$.

\paragraph{The composition map.}
Any possible relation can be expressed as 
\begin{equation*}
Q(\Psi^1_\sigma)_{\sigma\in S_n}=0
\end{equation*}
for a polynomial $Q\in\mathbb{F}[y_\sigma]_{\sigma\in S_n}$ with $\{y_\sigma\}_{\sigma\in S_n}$ formal variables. 

We can define the composition map $C_\Psi \colon \mathbb{F}[y_\sigma]_{\sigma\in S_n} \to R=\mathbb{F}[\mathbf{M}_n]$ for any polynomial $Q\in\mathbb{F}[y_\sigma]_{\sigma\in S_n}$ by $C_\Psi(Q)=Q(\Psi^1_\sigma)_{\sigma\in S_n}$ . This composition is an element of $R_0$ and hence and element of $R$. This map is an $\mathbb{F}$-homomorphism between $\mathbb{F}$-algebras. The relations in the ring $R_0$ correspond to $\ker{C_\Psi}$ which is an ideal of $\mathbb{F}[y_\sigma]_{\sigma\in S_n}$.
By Hilbert's Basis Theorem $\ker{C_\Psi}$ is finitely generated, hence $R_0$ is finitely presented. We now turn to describe explicitly a set of 
generating relations.

\paragraph{Reduction to monomials.}
Because all the polynomials $\Psi^1_\sigma$ for $\sigma\in S_n$ are monomials we have that $C_\Psi$ maps monomials to monomials. Each
$Q\in\ker{C_\Psi}$ can be written in the following form
\begin{equation}\label{eq:Q-monomials}
Q=\sum_{P\in\mathrm{Monomials}(R)} 
\sum_{S\in {C_\Psi}^{-1}(P)} b_S S
\end{equation}
Where we have grouped all the monomials appearing in $Q$ by their $C_\Psi$-image in $R$. 

Because $R$ is a free polynomial algebra $C_\Psi(Q)=0$ if and only if
\begin{equation*}
\sum_{S\in {C_\Psi}^{-1}(P)} b_S P=0
\end{equation*}
for each $P$ appearing in the decomposition of $Q$ in \eqref{eq:Q-monomials}. This is the same as
$\sum_{S\in {C_\Psi}^{-1}(P)} b_S=0$ for each such $P$. Hence $\ker{C_\Psi}$ is generated by elements of the form 
$Q=\sum_{S\in {C_\Psi}^{-1}(P)} b_S S$ for a fixed monomial $P\in R$ with $\sum_{S\in {C_\Psi}^{-1}(P)} b_S=0$.

We consider $(b_S)_{S}$ as a vector in the vector space $\mathbb{F}^{{C_\Psi}^{-1}(P)}$. This vector is contained in the subspace of vectors whose sum is zero. This subspace is spanned by vectors with a 
single $+1$ entry and a single $-1$ and the rest of the entries equal to $0$. Thus $\ker C_\Psi$ is generated by elements of the form $S_1-S_2$ for
two monomials in $\mathbb{F}[y_\sigma]_{\sigma\in S_n}$ such that $C_\Psi(S_1)=C_\Psi(S_2)$.

\paragraph{Translation to permutation matrices.}
Let $S_1, S_2\in \mathbb{F}[y_\sigma]_{\sigma\in S_n}$ be two monomials. Write $S_1=\prod_{\sigma\in S_n} y_\sigma^{a_\sigma}$ and $S_2=\prod_{\sigma\in S_n} y_\sigma^{b_\sigma}$. Then $C_\Psi(S_1)=C_\Psi(S_2)$ if and only if
\begin{equation*}
\sum_{\sigma\in S_n} (a_\sigma-b_\sigma)P^\sigma=0
\end{equation*}

Let $c_\sigma=a_\sigma-b_\sigma$ and define $a_\sigma^0\coloneqq \max(c_\sigma,0)$ and $b_\sigma^0\coloneqq \max(-c_\sigma,0)$ for all $\sigma\in S_n$. Let $S_1^0=\prod_{\sigma\in S_n} y_\sigma^{a_\sigma^0}$ and $S_2^0=\prod_{\sigma\in S_n} y_\sigma^{b_\sigma^0}$. Evidently $C_\Psi(S_1^0)=C_\Psi(S_2^0)$.
Moreover, if $\left<S_1^0-S_2^0\right>$ is the ideal in $\mathbb{F}[y_\sigma]_{\sigma\in S_n}$ generated by $S_1^0-S_2^0$ then $S_1-S_2\in\left<S_1^0-S_2^0\right>$.

To sum up, the ideal $\ker C_\Psi$ is generated by elements of the form $S_f^+=\prod_{\sigma\in S_n} y_\sigma^{\max(f(\sigma),0)}$ and $S_f^-=\prod_{\sigma\in S_n} y_\sigma^{\max(-f(\sigma),0)}$ for any function $f\colon S_n\to\mathbb{Z}_{\geq0}$ such that
\begin{equation}\label{eq:relations-orthogonality}
\sum_{\sigma\in S_n} f(\sigma)P^\sigma=0
\end{equation} 

\paragraph{The group algebra of $S_n$}
The group ring $\mathbb{Z}S_n$ consists of all function $S_n\to\mathbb{Z}$, it is an order in the group algebra $\mathbb{Q}S_n$ which consists of all the functions $S_n\to\mathbb{Q}$. All the irreducible representation of $S_n$ in characteristic zero are defined over $\mathbb{Q}$. Moreover, they can  be defined over $\mathbb{Z}$, i.e.\ they are representations of the form $S_n\to\mathbf{GL}_k(\mathbb{Z})$.

If $(\rho, W)$ is a representation of $S_n$ defined over $\mathbb{Q}$ we have the Fourier transform of any $f\in\mathbb{Q}S_n$ defined by
\begin{equation*}
\widehat{f}(\rho)=\sum_{\sigma\in S_n} f(\sigma)\cdot\rho(\sigma)\in{\End}_\mathbb{Q}(W)
\end{equation*}

Let $\rho_0=\rho_\mathrm{St}\oplus\rho_\mathrm{triv}$ be the direct sum of the standard representation of $S_n$ with the trivial representation. This is just the representation $\sigma\to P^\sigma\in\mathbf{GL}_n(\mathbb{Z})$ written as a sum of irreducible representations. Condition \eqref{eq:relations-orthogonality} can be translated to $\widehat{f}(\rho_0)=0$.

Let $\left\{(\rho_i,W_i)\right\}_i$ be the irreducible representation of $S_n$ defined over $\mathbb{Z}$. The Fourier transform induces an isomorphism $\mathbb{Q}S_n \xrightarrow{\sim} \bigoplus_i \End_\mathbb{Q}(W_i)$. By the Plancherel formula for the finite Fourier transform condition \eqref{eq:relations-orthogonality} becomes the statement that $f$ is orthogonal to the representations $\rho_\mathrm{St}$, $\rho_\mathrm{triv}$. 

The discussion above constitutes the proof of the following proposition.
\begin{proposition}\label{prop:relations}
The ring $\tensor[^{\mathbf{SA}}]{\mathbb{F} [\mathbf{M}_n]}{^{\mathbf{SA}}}$ is generated by the polynomials  $\left\{\Psi^1_\sigma\right\}_{\sigma\in S_n}$. A complete list of relations may be computed in the following way.

For each irreducible representation $(W,\rho)$ other then $\rho_\mathrm{St}$ and $\rho_\mathrm{triv}$ we look at $\End_\mathbb{Q}(W)$ as a subspace of $\mathbb{Q}{S_n}$. Let $f^\rho_1,\ldots,f^\rho_{{\dim W}^2}$ be a primitive base for the lattice $\mathbb{Z}S_n\cap \End_\mathbb{Q}(W)$. Then each base vector defines a relation
\begin{equation*}
\prod_{\sigma\in S_n} \left(\Psi^1_\sigma\right)^{\max(f^\rho_i(\sigma),0)}=
\prod_{\sigma\in S_n} \left(\Psi^1_\sigma\right)^{\max(-f^\rho_i(\sigma),0)}
\end{equation*}
\end{proposition}
\begin{remark}
Notice that the dimension of the maximal component of the variety $\dfaktor{\mathbf{SA}}{\mathbf{M}_n}{\mathbf{SA}}$ is $n^2-(n-1)-(n-1)=n^2-2n+2$. We have $|S_n|=n!$ generators for the ring and the amount of relations is $\dim \mathbb{Q}S_n-\dim{\rho_\mathrm{St}}^2-\dim{\rho_\mathrm{triv}}^2=n!-(n-1)^2-1$. The difference between the amount of generators and the amount of relations is $|S_n|-\left[\dim \mathbb{Q}S_n-\dim{\rho_\mathrm{St}}^2-\dim{\rho_\mathrm{triv}}^2\right]=n^2-2n+2$ which is equal to the above mentioned dimension.
\end{remark}
\begin{corollary}
The ring $\tensor[^{\mathbf{SA}}]{\mathbb{F} [\mathbf{SL}_n]}{^{\mathbf{SA}}}$ is generated by the regular functions  $\left\{\Psi^1_\sigma\right\}_{\sigma\in S_n}$ from the Proposition \ref{prop:relations} with the same relations and the additional Leibniz relation
\begin{equation*}
\sum_{\sigma\in S_n} \sign\sigma \Psi^1_\sigma=1
\end{equation*}
\end{corollary}
\begin{proof}
This follows from Proposition \ref{prop:relations} and \ref{cor:mn-sln}.
\end{proof}
\section{Double Torus Quotient for \texorpdfstring{$\mathbf{PGL}_2$}{PGL2}}\label{sec:dtq-pgl2}
In this part we demonstrate the connection between the canonical generators of $\dfaktor{\mathbf{PA}}{\mathbf{PGL}_2}{\mathbf{PA}}$ over $\mathbb{F}=\mathbb{Q}$ and the discriminant inner product of integral binary quadratic forms. This section is not formally required for the development of our results, yet it serves as a motivation to the study of the double quotient of a group by a torus. Moreover, we stress some inherent differences between $\mathbf{PGL}_2$ and higher rank cases.

The action of $\mathbf{PGL}_2$ on the the space of binary quadratic forms and its relation to class groups of quadratic fields is very well known and goes back to Gauss and Dirichlet. This is an extremely rich subject. We present a very concise introduction suited to our needs.

A close variant of this action has been exploited by Linnik to prove his result regarding equidistribution of integral points in the $2$-sphere. A central tool in his method is the study of inner product between two forms.

\subsection{Binary Quadratic Forms and the Adjoint Action}
We denote by $\mathcal{Q}$ the space of binary quadratic forms.
The action of $\mathbf{PGL}_2$ on $\mathcal{Q}$ is induced from the action of $\mathbf{GL}_2$ on $2$-vectors
\begin{equation*}
g.q(x,y)=\frac{1}{\det g}q((x,y)g)
\end{equation*}

Denote the Lie algebra of $\mathbf{PGL}_2$ by $\mathfrak{pgl}_2$. We identify $\mathfrak{pgl}_2$ with $\mathbf{M}_2^0$ the space of $2\times 2$ trace-free matrices. One can define in the following way an isomorphism over $\mathbb{Q}$ between the adjoint representation of $\mathbf{PGL}_2$ on $\mathfrak{pgl}_2$ and the the representation of $\mathbf{PGL}_2$ on the space of binary quadratic forms 
\begin{equation}\label{eq:quad-adj}
ax^2+bxy+cy^2\mapsto
\begin{pmatrix}
b   & -2a \\
2c  & -b
\end{pmatrix}
\end{equation}
Being an isomorphism of representations it intertwines the action of $\mathbf{PGL}_2$ on binary quadratic forms with the adjoint action on the Lie algebra.

The ring of invariants for the adjoint action of $\mathbf{PGL}_2$ on $\mathfrak{pgl}_2$ is generated by a single invariant, the determinant of the matrix. This is an instance of Chevalley's restriction theorem. The pullback of the determinant under the isomorphism of $\mathfrak{pgl}_2$ with $\mathcal{Q}$ is a multiple of the discriminant of a binary quadratic form. Specifically, we have an isomorphism of quadratic spaces
\begin{equation*}
\left(\mathcal{Q},{\disc}\right)\longleftrightarrow \left(\mathfrak{pgl}_2, -{\det} \right)
\end{equation*}

\subsection{Stabilizers of a Binary Quadratic Form}
The stabilizer in $\mathbf{PGL}_2$ of a non-degenerate \emph{rational} binary quadratic form $q$ is the centralizer of a non-trivial regular element in $\mathfrak{pgl}_2(\mathbb{Q})$, hence it is a maximal torus $\mathbf{T}$ defined over $\mathbb{Q}$.

The torus $\mathbf{T}$ is anisotropic over $\mathbb{Q}$ if and only if $q$ is irreducible over $\mathbb{Q}$. Moreover it is split over $\mathbb{R}$ if and only if $q$ is reducible over $\mathbb{R}$, i.e.\ $\disc(q)>0$.

Rational forms $q$ and $q'$ are homothetic if $\exists \alpha\in\mathbb{Q}^\times$ such that $q=\alpha q'$.
The correspondence between homothety classes of non-degenerate rational binary quadratic forms and rational maximal tori is a bijection. We have already seen how to attach a rational maximal torus to a quadratic form in a way which is obviously invariant under homothety. In the other direction, for a maximal rational torus $\mathbf{T}$ its Lie algebra $\mathfrak{t}(\mathbb{Q})$ is a rational subspace of $\mathfrak{g}(\mathbb{Q})$. Using the isomorphism \eqref{eq:quad-adj} this corresponds to a homothety class of rational quadratic forms. It is straight forward to see that these maps are inverse to each other. For example, the standard diagonal torus corresponds to the homothety class of $q_0(x,y)=xy$.

\subsection{Canonical Generators for \texorpdfstring{$\mathbf{PA}$}{PA} in \texorpdfstring{$\mathbf{PGL}_2$}{PGL2}}
Denote by $\mathfrak{a}$ the Lie algebra of the maximal torus $\mathbf{PA}<\mathbf{PGL}_2$. The Weyl group of $\mathbf{PA}$ consists of two elements: the identity element, $+1$, which acts trivially on $\mathfrak{a}$ and the element $-1$ which acts by negation on $\mathfrak{a}$.

From our study in the previous parts we learn that the canonical generators for $\dfaktor{\mathbf{PA}}{\mathbf{PGL}_2}{\mathbf{PA}}$ are $\Psi_{+1}$ and $\Psi_{-1}$. The algebra of invariants for the double torus quotient is generated by those polynomials and they are related by a single relation coming from the Leibniz formula
\begin{equation*}
\Psi_{+1}+\Psi_{-1}=1
\end{equation*}
In particular, the algebra of invariants is $\tensor[^{\mathbf{PA}}]{\mathbb{Q} [\mathbf{PGL}]}{^{\mathbf{PA}}}\cong \mathbb{Q}[\Psi_{+1}]$, i.e.\ the double quotient space is just a one dimensional affine space.

If we look at the definition of the canonical generators at Proposition \ref{prop:generators-matrix} adding the twist by the determinant for $\mathbf{PGL}_2$ we see that
\begin{equation*}
\Psi_{+1}\left[\begin{pmatrix}
a & b\\
c & d
\end{pmatrix}\right]
=\frac{1}{ad-bc}ad
\quad
\Psi_{-1}\left[\begin{pmatrix}
a & b\\
c & d
\end{pmatrix}\right]
=-\frac{1}{ad-bc}bc
\end{equation*}

\subsection{Canonical Generators for a rational torus in \texorpdfstring{$\mathbf{PGL}_2$}{PGL2}}
Let now $\mathbf{T}<\mathbf{PGL}_2$ be any maximal torus defined over $\mathbb{Q}$. Let $\mathbb{L}/\mathbb{Q}$ be the splitting field of $\mathbf{T}$, either $\mathbb{L}=\mathbb{Q}$ and then $\mathbf{T}=\mathbf{PA}$ or $\mathbb{L}$ is a quadratic extension of $\mathbb{Q}$ and then $\mathbf{T}$ is anisotropic over $\mathbb{Q}$. In both cases $\mathbf{T}$ is conjugate to $\mathbf{PA}$ by an element of $\mathbf{PGL}_2(\mathbb{L})$.

We are going now to construct generators for $\tensor[^{\mathbf{T}}]{\mathbb{L} [\mathbf{PGL}_2]}{^{\mathbf{T}}}$ from the functions $\Psi_{+1}$ and $\Psi_{-1}$. This is a simple instance of the construction we carry out in Proposition \ref{prop:generators-splitting}.

As $\mathbf{T}$ is $\mathbb{L}$-split there exists $g\in\mathbf{PGL}_2(\mathbb{L})$ such that $\mathbf{T}_\mathbb{L}=\Ad_{g} \mathbf{PA}_\mathbb{L}$. This equality defines uniquely the coset of $g$ in $\faktor{\mathbf{PGL}_2(\mathbb{L})}{\Nrml_{\mathbf{PGL}_2}(\mathbf{PA})(\mathbb{L})}$, where $\Nrml_{\mathbf{PGL}_2}(\mathbf{PA})$ is the normalizer of the split torus $\mathbf{PA}$ in $\mathbf{PGL}_2$.
We define
\begin{equation*}
\Psi^\mathbf{T}_{\pm 1}\coloneqq\Psi_{\pm 1}\circ {\Ad}_{g^{-1}}
\end{equation*}
One easily sees by direct computation that $\Psi_{+1}^\mathbf{T}$ and $\Psi_{-1}^\mathbf{T}$ are invariant under the left and right action of $\mathbf{T}_\mathbb{L}$ and that the definition does not depend on the choice of $g$ up to permuting the two generators.

More so, the map $P\mapsto P\circ \Ad_{g^{-1}}$ defines an isomorphism of $\mathbb{L}$ algebras $\tensor[^{\mathbf{PA}}]{\mathbb{L} [\mathbf{PGL}_2]}{^{\mathbf{PA}}}\to\tensor[^{\mathbf{T}}]{\mathbb{L} [\mathbf{PGL}_2]}{^{\mathbf{T}}}$. If $\mathbb{L}\neq\mathbb{Q}$ this isomorphism is a priori not defined over $\mathbb{Q}$ as $g\not\in\mathbb{Q}$. Yet it implies that $\Psi_{+1}^\mathbf{T}$ and $\Psi_{-1}^\mathbf{T}$ generate the whole ring of left and right $\mathbf{T}_\mathbb{L}$-invariant regular functions over $\mathbb{L}$.

What we show in the next section is that although $\Psi_{+1}^\mathbf{T}$ and $\Psi_{-1}^\mathbf{T}$ are defined initially over $\mathbb{L}$, for any rational point $\lambda\in\mathbf{PGL}_2(\mathbb{Q})$ one has that $\Psi_{\pm 1}^\mathbf{T}(\lambda)\in\mathbb{Q}$. Moreover, these functions are actually defined over $\mathbb{Q}$.

Unfortunately, this is \emph{no more true} for higher rank groups. Nevertheless, an analogues statement, Proposition \ref{prop:galois-weyl}, regarding the orbits of the canonical generators under the Galois group $\Gal(\mathbb{L}/\mathbb{Q})$ still holds.
 
\subsection{Discriminant Inner Product}
We have a natural invariant for the $\mathbf{PGL}_2$ action on $\mathcal{Q}$ which is the discriminant. Being a quadratic form it is associated to a bilinear form on $\mathcal{Q}$ which is its polarization. We call this bilinear form the discriminant inner product. Notice that it is not positive definite yet it is non-degenerate. For two binary quadratic forms this inner product is evaluated by
\begin{equation*}
\left<ax^2+bxy+cy^2, a'x^2 +b'xy+c'y^2\right>_\mathrm{disc}=bb'-2ac'-2a'c
\end{equation*}

\begin{proposition}
Let $q_\mathbf{T}(x,y)$ be a quadratic form of discriminant $1$ corresponding to the rational torus $\mathbf{T}<\mathbf{PGL}_2$.
Let $\delta \in\mathbf{G}(\mathbb{C})$ then
\begin{equation*}
(\det\delta)^{-1}\left<q_\mathbf{T},\delta. q_\mathbf{T}\right>_{\disc}=\Psi^\mathbf{T}_{+1}(\delta)-\Psi^\mathbf{T}_{-1}(\delta)
\end{equation*}
\end{proposition}
\begin{proof}
First notice that the polynomial map $\delta\mapsto (\det\delta)^{-1}\left<q_\mathbf{T},\delta\cdot q_\mathbf{T}\right>_{\disc}$ is invariant under both the left and the right action of $\mathbf{T}$. 

To see that the polynomials are actually equal, write $\mathbf{T}=g\mathbf{PA}g^{-1}$ with $g\in\mathbf{G}(\mathbb{C})$. Let $q_0(x,y)=xy$ be the discriminant $1$ binary quadratic form of the torus $\mathbf{PA}$. 
We have $q_\mathbf{T}(x,y)=\pm g.q_0(x,y$). Without loss of generality the sign in the equality is positive.

We now reduce the claim to the case $\mathbf{T}=\mathbf{PA}$.
\begin{align*}
&(\det\delta)^{-1}\left<q_\mathbf{T},\delta. q_\mathbf{T}\right>_{\disc}=
(\det g\delta g^{-1})^{-1}\left<q_0,(g^{-1} \delta g). q_0\right>_{\disc}\\
&\Psi^\mathbf{T}_{+1}(\delta)-\Psi^\mathbf{T}_{-1}(\delta)=
\Psi_{+1}(g^{-1} \delta g)-\Psi_{-1}(g^{-1} \delta g)
\end{align*}

For $\mathbf{T}=\mathbf{PA}$ the claims follows from a direct computation of both expressions.
\end{proof}
\begin{corollary}
The polynomials $\Psi^\mathbf{T}_{+1}$ and $\Psi^\mathbf{T}_{-1}$ are defined over $\mathbb{Q}$.
\end{corollary}
\begin{proof}
As $q_\mathbf{T}$ is a quadratic form over $\mathbb{Q}$ and $\left<,\right>_\mathrm{disc}$ is defined over $\mathbb{Q}$ we deduce from the previous proposition that $\Psi^\mathbf{T}_{+1}-\Psi^\mathbf{T}_{-1}$ is defined over $\mathbb{Q}$. 

The claim follows by the relation $\Psi^\mathbf{T}_{+1} + \Psi^\mathbf{T}_{-1}=1$.
\end{proof}
\begin{remark}
In higher rank we do not have an available linear representation having the same properties as the representation of $\mathbf{PGL}_2$ on $\mathcal{Q}$. The natural analogue is the action of $\mathbf{PGL}_n$ on $\bigwedge^{n-1}\mathfrak{pgl}_n$ which is fruitfully exploited in \cite{ELMVPeriodic}. Yet this representation for $n>2$ has plenty of lines whose stabilizer's identity component is not  a maximal torus.

In what follows we have to study the canonical generators attached to rational tori without the aid of a linear representation. This is mainly done using the action of the Galois group of the splitting field on the generators and the algebraic relations between them.

We note that Bhargava \cite{BhaCubic} has discovered a linear representation of a higher rank group such that the identity components of all the line stabilizers are tori, although not maximal ones. See also the generalization by Wood \cite{Wood}.
\end{remark}

\section{Double Torus Quotient for Central Simple Algebras}\label{sec:dtq-csa}
We are now in position to use the results of the previous parts to study double torus quotients in groups isogenous to an inner form of $\mathbf{PGL}_n$, e.g.\ the projective group of units and the group of units of norm 1 in a central simple algebra.

Let $\mathbf{B}$ be a central simple algebra over $\mathbb{F}$ considered as an $\mathbb{F}$-algebra object in the category of affine algebraic varieties.
Denote the reduced norm by $\Nrd\colon\mathbf{B}\to\mathbb{A}^1$. When restricted to the group of units, $\mathbf{GL}_1(\mathbf{B})$, it is an $\mathbb{F}$-character.

Let $\mathbf{G}$ be a reductive linear algebraic defined over $\mathbb{F}$ with fixed $\mathbb{F}$-isogenies
\begin{equation}\label{eq:fixed-isogenies}
\mathbf{SL}_1(\mathbf{B})\to\mathbf{G}\to\mathbf{PGL}_1(\mathbf{B})
\end{equation}
where $\mathbf{PGL}_1(\mathbf{B})$ is the adjoint form and 
$\mathbf{SL}_1(\mathbf{B})$ is the simply-connected cover. The linear algebraic groups
$\mathbf{SL}_1(\mathbf{B})$ and $\mathbf{PGL}_1(\mathbf{B})$ are inner $\mathbb{F}$-forms of $\mathbf{SL}_n$ and $\mathbf{PGL}_n$ respectively. 
All inner $\mathbb{F}$-forms of $\mathbf{PGL}_n$ and $\mathbf{SL}_n$ arise this way. 
See \cite[Chapter III, \S 1.4]{GaloisSerre} for details.

Let $\mathbf{T}<\mathbf{G}$ be a maximal torus defined over $\mathbb{F}$. By $\mathbf{T}$ we mean an algebraic torus together with a fixed embedding into $\mathbf{G}$.
There is a unique maximal commutative subalgebra $\mathbf{E}<\mathbf{B}$ defined over $\mathbb{F}$ such that the isogeny \eqref{eq:fixed-isogenies} restricts to an isogeny of tori
\begin{equation*}
\mathbf{SL}_1(\mathbf{E})\to\mathbf{T}\to\mathbf{PGL}_1(\mathbf{E})
\coloneqq\lfaktor{\mathbf{Z}_{\mathbf{B}^\times}}{\mathbf{E}^\times}
\end{equation*}
Notice that $\mathbf{PGL}_1(\mathbf{E})$ is defined to be the quotient of the units of $\mathbf{E}$ by the center of $\mathbf{B}^\times$.

The variety $\mathbf{E}$ can be reconstructed from $\mathbf{T}$ by identifying the Lie algebra of $\mathbf{B}$ with the traceless quaternions $\mathbf{B}^0<\mathbf{B}$ and defining $\mathbf{E}<\mathbf{B}$ to be the direct sum of $\Lie\mathbf{T}$ and the one-dimensional space spanned by the identity in $\mathbf{B}$.

The ring  
$\mathbb{K}\coloneqq\mathbf{E}(\mathbb{F})$ is an \'{e}tale-algebra of degree $n$ over $\mathbb{F}$. The torus $\mathbf{T}$ is anisotropic over $\mathbb{F}$ if and only
if $\mathbb{K}$ is a field.

\subsection{Double Torus Quotient over a Splitting Field}
We do not have an explicit description of the double torus quotient $\dfaktor{\mathbf{T}}{\mathbf{G}}{\mathbf{T}}$
over $\mathbb{F}$, but as we will see in this section we can use our description of $\dfaktor{\mathbf{SA}}{\mathbf{SL}_n}{\mathbf{SA}}$ to derive an explicit description of $\dfaktor{\mathbf{T}}{\mathbf{G}}{\mathbf{T}}$ after extending the scalars to a splitting field of $\mathbf{T}$. In other words, the double torus quotient varieties for different rational tori $\mathbf{T}$ are all $\mathbb{F}$-forms of the standard double torus quotient $\dfaktor{\mathbf{SA}}{\mathbf{SL}_n}{\mathbf{SA}}$. 

We say that a commutative algebra of degree $n$ over $\mathbb{L}$ is split if it is isomorphic as an algebra to $\mathbb{L}^n$.
Denote by $\mathbb{L}/\mathbb{F}$ the splitting field of $\mathbf{T}$, equivalently $\mathbf{E}$, 
i.e.\ the unique minimal field extension of $\mathbb{F}$ over which $\mathbf{T}$, equivalently $\mathbf{E}$ and $\mathbf{GL}_1(\mathbf{E})$, split.

A splitting field for $\mathbf{E}$ always splits the central simple algebra as well. That is there exists an $\mathbb{L}$-isomorphism of algebras
$\phi\colon\mathbf{B}_\mathbb{L}\to\mathbf{M}_{n,\mathbb{L}}$ (see the discussion in \cite[\S 2.2]{GilleCentral}). The isomorphism $\phi$ sends $\mathbf{E}_\mathbb{L}$ to a split maximal commutative subalgebra in $\mathbf{M}_{n,\mathbb{L}}$. We insure that this subalgebra is the standard diagonal one by replacing $\phi$ with the composition of $\phi$ with an inner automorphism sending $\phi(\mathbf{E}_\mathbb{L})$ to the standard diagonal subalgebra\footnote{One doesn't need to work with the standard diagonal subalgebra here, any  maximal commutative subalgebra of $\mathbf{M}_{n,\mathbb{L}}$ defined and split over $\mathbb{F}$ will do.}, which we denote by $\mathbf{Diag}_\mathbb{L}$.

The isomorphism $\phi$ induces isomorphisms
\begin{align*}
\mathbf{SL}_1(\mathbf{B})_\mathbb{L}&\to\mathbf{SL}_{n,\mathbb{L}}\\
\mathbf{PGL}_1(\mathbf{B})_\mathbb{L}&\to\mathbf{PGL}_{n,\mathbb{L}}
\end{align*}
which we also denote by $\phi$. 
Moreover, base-changing the fixed isogenies \eqref{eq:fixed-isogenies} and composing with these isomorphisms we have fixed isogenies over $\mathbb{L}$
\begin{align}\label{eq:fixed-isogenies-L}
\mathbf{SL}_{n,\mathbb{L}}&\to\mathbf{G}_\mathbb{L}\to\mathbf{PGL}_{n,\mathbb{L}}\\
\mathbf{SA}_{\mathbb{L}}&\to\mathbf{T}_\mathbb{L}\to\mathbf{PA}_{\mathbb{L}}\nonumber
\end{align}
Using this isogeny we consider every regular function on $\mathbf{PGL}_{n,\mathbb{L}}$ as a regular function on $\mathbf{G}_\mathbb{L}$. 

\begin{proposition}
The isogenies \eqref{eq:fixed-isogenies-L} induce \emph{isomorphisms} 
\begin{equation}\label{eq:double-quotient-isogen-iso}
\left(\dfaktor{\mathbf{SA}}{\mathbf{SL}_n}{\mathbf{SA}}\right)_\mathbb{L}
\to
\left(\dfaktor{\mathbf{T}}{\mathbf{G}}{\mathbf{T}}\right)_\mathbb{L}
\to
\left(\dfaktor{\mathbf{PA}}{\mathbf{PGL}_n}{\mathbf{PA}}\right)_\mathbb{L}
\end{equation}
\end{proposition}
\begin{proof}
The properties of the universal categorical quotient in Theorem \ref{thm:quotient} assure us that $\left(\dfaktor{\mathbf{T}}{\mathbf{G}}{\mathbf{T}}\right)_\mathbb{L}\cong\dfaktor{\mathbf{T}_\mathbb{L}}{\mathbf{G}_\mathbb{L}}{\mathbf{T}_\mathbb{L}}$ and the same holds for $\mathbf{PGL}_{n,\mathbb{L}}$ and $\mathbf{SL}_{n,\mathbb{L}}$.

Because the isomorphism in Proposition \ref{prop:pgl-generators} factors through \eqref{eq:double-quotient-isogen-iso} we conclude that all the maps in \eqref{eq:double-quotient-isogen-iso} are isomorphisms.
\end{proof}

We can now transport using $\phi$ our explicit coordinates $\left\{ \Psi^0_\sigma \right\}_{\sigma\in S_n}$ from $\dfaktor{\mathbf{SA}}{\mathbf{SL}_n}{\mathbf{SA}}$ to $\left(\dfaktor{\mathbf{T}}{\mathbf{G}}{\mathbf{T}}\right)_\mathbb{L}$. We wish to rewrite the polynomials $\Psi^0_\sigma$ using intrinsic notions similar to \S \ref{sec:pgl-generators}. For this we bring the classical definition of a complete set of primitive orthogonal idempotents.

\begin{definition} 
An element of a unital ring $e_0\in R$ is called an \emph{idempotent} if $e_0^2=e_0$. The idempotent $e_0$ is \emph{non-trivia}l if $e_0\neq 0,1$ and it is \emph{primitive} if it cannot be written a sum of two non-trivial idempotents. 
	
Two idempotents $e_1,e_2\in R$ are \emph{orthogonal} if $e_1 e_2= e_2 e_1=0$. A set of idempotents $e_1,\ldots,e_n\in R$ is called a \emph{complete set of primitive orthogonal idempotents} if all the $e_i$'s
are mutually orthogonal primitive idempotents and $1=e_1+\ldots+e_n$.

If $R$ is commutative then a complete set of primitive orthogonal idempotents is unique, up to permutation, if it exists.
\end{definition}

\begin{proposition}\label{prop:generators-splitting}
\hfill
\begin{enumerate}
\item
For each $1\leq i \leq n$ let $e_i^0$ be the diagonal matrix with all zero entries except for a single $1$ entry in the $(i,i)$ place.

A complete set of primitive orthogonal idempotents for $\mathbf{Diag}_\mathbb{L}(\mathbb{L})$ is given by
$\left\{e_i^0\right\}_{i=1}^n$. The $\phi^{-1}$ image of this set is a complete set of primitive orthogonal idempotents in $\mathbf{E}_\mathbb{L}(\mathbb{L})$. Write $e_i=\phi^{-1}(e_i^0)$ for those primitive orthogonal
idempotents.

\item We identify the absolute Weyl group with symmetric group on the ordered complete set of primitive orthogonal idempotents and with the standard symmetric group on $\left\{1,\ldots, n\right\}$ in a consistent manner. For all $\sigma\in\mathbf{W}_\mathbf{T}(\mathbb{L})$ we have $e_{\sigma(i)}=\sigma.e_i$. Notice that this identification depends on the order of $e_1,\ldots,e_n$.

\item
The pullback to $\left(\dfaktor{\mathbf{T}}{\mathbf{G}}{\mathbf{T}}\right)_\mathbb{L}$ of $\Psi^0_\sigma$ is 
\begin{equation*}
\Psi_\sigma(g)\coloneqq\left(\Psi^0_\sigma\circ\phi\right)(g)= \Nrd\left(\sum_{i=0}^n (\sigma.e_i) g e_i\right)\cdot \Nrd(g)^{-1}
\end{equation*}
We  call the polynomials $\left\{\Psi_\sigma\right\}_{\sigma\in\mathbf{W}_\mathbf{T}(\mathbb{L})}$ the \textbf{canonical generators} of
the ring of regular functions of $\left(\dfaktor{\mathbf{T}}{\mathbf{G}}{\mathbf{T}}\right)_\mathbb{L}$.
\end{enumerate}
\end{proposition}
\begin{proof}
The first two parts of the proposition are straightforward.
The third part follows immediately from the facts that the reduced norm is exactly the determinant map for the matrix algebra and that $\phi$, being an isomorphism of central simple algebras, is reduced-norm preserving.
\end{proof}

The idempotents $e_1,\ldots,e_n\in E_\mathbb{L}\coloneqq\mathbf{E}_\mathbb{L}(\mathbb{L})$ form a base for $E_\mathbb{L}$ as an $\mathbb{L}$-vector space. One can use other bases for $E_\mathbb{L}$ to express the canonical generators. To do this we will need the reduced trace form. The central simple algebra $\mathbf{B}_\mathbb{L}$ carriers a natural non-degenerate bilinear form, which is the reduced trace form $\Trd$. This form restricts to a non-degenerate bilinear form on $E_\mathbb{L}$ which is just the standard trace form on the \'{e}tale algebra.

\begin{lemma}\label{lem:dual_base}
Let $b_1,\ldots,b_n \in E_\mathbb{L}$ be a basis for $E_\mathbb{L}\coloneqq\mathbf{E}_\mathbb{L}(\mathbb{L})$ as an $\mathbb{L}$-vector space. Let $\widecheck{b}_1,\ldots,\widecheck{b}_n$ the dual basis with respect to the reduced trace, i.e.\ $\Trd(\widecheck{b}_i b_j)=\delta_{i,j}$. For any $\sigma$ in the absolute Weyl group we have
\begin{equation*}
\Psi_\sigma (g)=\Nrd\left(\sum_{i=1}^n \sigma.\widecheck{b}_i \cdot g  \cdot  b_i \right)\cdot \Nrd(g)^{-1}
\end{equation*}
\end{lemma}
\begin{proof}
For a linear endomorphism $\rho\in \End_\mathbb{L}(E_\mathbb{L})$ denote by $\rho^\dagger$ the adjoint operator with respect to the reduced trace.

Notice that $e_1,\ldots,e_n$ is an orthonormal basis with respect to the reduced trace. Hence if $(\rho_{i,j})_{i,j=1}^n$ is the matrix associated to the endomorphism $\rho$ with respect to the basis $e_1,\ldots,e_n$, then the matrix associated with $\rho^\dagger$ is just the transposed matrix $(\rho_{j,i})_{i,j=1}^n$. 

Choose $\rho\in\End_\mathbb{L}(E)$ such that $\rho(e_i)=b_i$ for all $i=1,\ldots,n$, then $\rho^\dagger(\widecheck{b}_i)=e_i$. Write $b_i=\sum_{j=1}^n \rho_{i,j} e_j$. A calculation now gives
\begin{align*}
\sum_{i=1}^n \sigma.\widecheck{b}_i \cdot g  \cdot  b_i &=
\sum_{i=1}^n \sigma.\widecheck{b}_i \cdot g  \cdot \left(\sum_{j=1}^n \rho_{i,j} e_j \right) =
\sum_{j=1}^n \sigma.\left(\sum_{i=1}^n \rho_{i,j} \widecheck{b}_i \right) \cdot g \cdot e_j \\
&= \sum_{j=1}^n \sigma.\rho^\dagger( \widecheck{b}_i) \cdot g \cdot e_j =
\sum_{j=1}^n \sigma.e_j \cdot g \cdot e_j
\end{align*}
This concludes the proof.
\end{proof}

\subsection{Galois Action}
The splitting field $\mathbb{L}$ is always a Galois extension of the base field $\mathbb{F}$. 
In this section we study the action of the Galois group $\mathfrak{G\coloneqq \Gal\left(\mathbb{L}/\mathbb{F}\right)}$
on the map $\phi$ and derive from this the way the Galois group acts on the algebraic numbers $\Psi_\sigma(g)$ for $g\in\mathbf{G}(\mathbb{F})$.

\paragraph{Galois Descent.}
We briefly review basic facts from the theory of Galois descent for
affine algebraic varieties in a form useful to us. This theory is sufficient for our needs as all our varieties are affine. The affine case is straightforward to verify using the analogues statements for the coordinate rings.

Let $\mathbf{X}$ be an affine variety defined over $\mathbb{F}$.
There is a natural action of $\mathfrak{G\coloneqq \Gal\left(\mathbb{L}/\mathbb{F}\right)}$ on the scheme $\mathbf{X}_\mathbb{L}$ by automorphisms. These automorphisms are not defined over $\mathbb{L}$, as they act non-trivially on $\Spec\mathbb{L}$ through the Galois action on $\mathbb{L}/\mathbb{F}$.
By an abuse of notation we denote the automorphism of $\mathbf{X}_\mathbb{L}$ induced by $\sigma\in\mathfrak{G}$ also by $\sigma$.

Given two affine algebraic varieties $\mathbf{X}$, $\mathbf{Y}$ defined over $\mathbb{F}$, we have an induced  action of $\mathfrak{G}$ on $\Mor_\mathbb{L}(\mathbf{X}_\mathbb{L},\mathbf{Y}_\mathbb{L})$ by $f^\sigma\coloneqq\sigma \circ f \circ \sigma^{-1}$. If $f$ is an isomorphism then so is $f^\sigma$ for each $\sigma\in\mathfrak{G}$. This action extends to the arrow category in the following
way. If we have the following morphisms of algebraic varieties over $\mathbb{F}$
\begin{center}
	\begin{tikzcd}
		\mathbf{X}' \arrow{d} & \mathbf{Y}' \arrow{d}\\
		\mathbf{X}  & \mathbf{Y} \\
	\end{tikzcd}
\end{center}

Together with additional horizontal morphisms making the following diagram commute
\begin{center}
	\begin{tikzcd}
		\mathbf{X}'_\mathbb{L} \arrow{d} \arrow{r}{f'} & \mathbf{Y}'_\mathbb{L} \arrow{d}\\
		\mathbf{X}_\mathbb{L}  \arrow{r}{f} & \mathbf{Y}_\mathbb{L} \\
	\end{tikzcd}
\end{center}

Then the following diagram commutes for each $\sigma\in\mathfrak{G}$
\begin{center}
	\begin{tikzcd}
		\mathbf{X}'_\mathbb{L} \arrow{r}{f'^\sigma} \arrow{d} & \mathbf{Y}'_\mathbb{L} \arrow{d}\\
		\mathbf{X}_\mathbb{L} \arrow{r}{f^\sigma} & \mathbf{Y}_\mathbb{L} \\
	\end{tikzcd}
\end{center}

Set $\Aut\left(\mathbf{Y}'_\mathbb{L}\rightarrow\mathbf{Y}_\mathbb{L}\right)$ to be the group of pairs of automorphisms of $\mathbf{Y}'_\mathbb{L}$ and $\mathbf{Y}_\mathbb{L}$ which intertwine with the vertical morphism $\mathbf{Y}'_\mathbb{L}\rightarrow \mathbf{Y}_\mathbb{L}$.
The map $\sigma\mapsto \left(f'^\sigma \circ f'^{-1}, f^\sigma \circ f^{-1}\right) $ is 1-cocycle in $H^1\left(\mathfrak{G}, \Aut\left(\mathbf{Y}'_\mathbb{L}\rightarrow\mathbf{Y}_\mathbb{L}\right)\right)$.

Moreover if we denote $\mathfrak{G}_{f'}\coloneqq\left\{\sigma\in\mathfrak{G} \mid f'^\sigma=f' \right\}$ and $\mathbb{M}/\mathbb{F}$ is the field extension
corresponding to the subgroup $\mathfrak{G}_{f'}$ then Galois descent for affine varieties 
provides the existence of the dotted morphism over $\mathbb{M}$ making the following diagram commute

\begin{center}
	\begin{tikzcd}
		\mathbf{X}'_\mathbb{M} \arrow[dashed]{r}{f'_\mathbb{M}}  & \mathbf{Y}'_\mathbb{M} \\
		\mathbf{X}'_\mathbb{L} \arrow{u} \arrow{r}{f'} & \mathbf{Y}'_\mathbb{L} \arrow{u} \\
	\end{tikzcd}
\end{center}
In particular, if $f'$ is an isomorphism then so is $f'_\mathbb{M}$.

\paragraph{Galois Descent for The Map $\phi$.}
We are going to apply these generalities to the case $\mathbf{X}=\mathbf{B}$, $\mathbf{Y}=\mathbf{M}_n$, $\mathbf{X}'=\mathbf{GL}_1(\mathbf{E})$, $\mathbf{Y}'=
\mathbf{A}$. The morphisms are $\phi$ and $\phi\restriction_{\mathbf{GL}_1(\mathbf{E})_\mathbb{L}}$, the restriction of $\phi$ to $\mathbf{GL}_1(\mathbf{E})_\mathbb{L}$. The schematic image of $\phi\restriction_{\mathbf{GL}_1(\mathbf{E})_\mathbb{L}}$ is exactly $\mathbf{A}_\mathbb{L}$.

The Skolem-Noether theorem implies that for each $\sigma\in \mathfrak{G}$ there exists $n_\sigma\in \mathbf{GL}_n(\mathbb{L})$ such that 
$\phi^\sigma=\Ad_{n_\sigma}\circ\phi$ where $\Ad_{n_\sigma}$ is the conjugation by $n_\sigma$. The equivalence class of $n_\sigma$ in $\mathbf{PGL}_n(\mathbb{L})$ is
uniquely determined by this equality. 
Denote this class by $[n_\sigma]$.

Denote by $\Nrml_{\mathbf{PGL}_n}(\mathbf{PA})$ the normalizer group of the split torus $\mathbf{PA}$ in $\mathbf{PGL}_n$. 
The following holds because the Galois action is well defined in the arrow category
\begin{equation*}
\left(\phi\restriction_{\mathbf{GL}_1(\mathbf{E})_\mathbb{L}}\right)^\sigma={\Ad}_{n_\sigma}\circ\phi\restriction_{\mathbf{GL}_1(\mathbf{E})_\mathbb{L}}\,.
\end{equation*}
A fortiori, 
$\Ad_{n_\sigma}(\mathbf{PA}_\mathbb{L})=\mathbf{PA}_\mathbb{L}$, thus $[n_\sigma]\in \Nrml_{\mathbf{PGL}_n}(\mathbf{PA})(\mathbb{L})$. 

We can now define a natural homomorphism of finite groups $\eta\colon\mathfrak{G}\to\mathbf{W}(\mathbb{L})$, where $\mathbf{W}$ is the Weyl group of $\mathbf{PA}$. This homomorphism is defined by the composition of the following maps
\begin{equation}\label{eq:eta-def}
\mathfrak{G}\xrightarrow{\sigma\mapsto [n_\sigma]} \Nrml_{\mathbf{PGL}_n}(\mathbf{PA})(\mathbb{L}) \rightarrow 
\faktor{\Nrml_{\mathbf{PGL}_n}(\mathbf{PA})}{\mathbf{PA}}(\mathbb{L}) = \mathbf{W}(\mathbb{L})
\end{equation}
Because $\mathbf{PA}$ is split over $\mathbb{L}$ the group $\mathbf{W}(\mathbb{L})$ is equal to the absolute Weyl group of $\mathbf{PA}$. 

Note that a priori, the map 
$\mathfrak{G}\rightarrow\mathbf{W}(\mathbb{L})$ is only a 1-cocycle, but as the torus $\mathbf{PA}$ is already split over $\mathbb{F}$ the action of $\mathfrak{G}$ on $\mathbf{W}(\mathbb{L})$ is trivial, hence $\eta$ is a homomorphism.

Moreover, we have a natural isomorphism between the Weyl group of $\mathbf{PA}$ in $\mathbf{PGL}_n$ and the Weyl group of $\mathbf{T}_\mathbb{L}$ in $\mathbf{G}_\mathbb{L}$, we identify them both by abuse of notation. 
\begin{proposition}\label{prop:galois-weyl}
The homomorphism of finite groups $\eta\colon\mathfrak{G}\to\mathbf{W}(\mathbb{L})$ defined in \eqref{eq:eta-def} is injective.	
\end{proposition}
\begin{proof}
By abuse of notation we identify elements of the Weyl group and the corresponding automorphisms of the torus.
The discussion above implies that for all $\sigma\in\mathfrak{G}$ we have $\phi\restriction_{\mathbf{GL}_1(\mathbf{E})_\mathbb{L}}^\sigma=\eta(\sigma)\circ\phi\restriction_{\mathbf{GL}_1(\mathbf{E})_\mathbb{L}}$.
Therefore, the Galois group stabilizer of $\phi\restriction_{{\mathbf{PGL}_1(\mathbf{E})}_\mathbb{L}}$ is $\mathfrak{G}_{\phi\restriction_{{\mathbf{PGL}_1(\mathbf{E})}_\mathbb{L}}}=\ker\eta$. Galois descent implies that the torus $\mathbf{GL}_1(\mathbf{E})$ is split already over the subfield
which corresponds to $\ker\eta$. Yet $\mathbb{L}$ is the \emph{minimal} splitting field of $\mathbf{GL}_1(\mathbf{E})$, hence $\ker\eta$ is trivial.
\end{proof}

\subsubsection{Galois Action on Canonical Generators}
We see that the Galois action on $\phi$ factors through the Weyl group. This allows us to compute the way in which the canonical generators of a rational point transform under the Galois group. By abuse of notation we  denote both the $\mathbb{L}$-points in the Weyl group of $\mathbf{PA}$ and the $\mathbb{L}$-points in the Weyl group of $\mathbf{T}$ by $W$, those abstract groups are isomorphic through $\phi$. In addition, because of the previous proposition we can treat $\mathfrak{G}$ as a subgroup of $W$. Henceforth we avoid writing explicitly the homomorphism $\eta$.

We turn to describe the action of the Galois group on the canonical generators at a rational point.
\begin{proposition}\label{prop:psi-galois}
Let $g\in\mathbf{G}(\mathbb{F})$ and $\tau\in\mathfrak{G}$, then for every $\sigma\in W$ 
\begin{equation*}
\tau. \Psi_\sigma(g) = \Psi_{\tau\sigma\tau^{-1}}(g)
\end{equation*}
\end{proposition}
\begin{proof}
Let $n_\tau\in \Nrml_{\mathbf{PGL}_n}(\mathbf{PA})(\mathbb{L})$ be a representative of $\tau\in W$ such that $\phi^\tau=\Ad_{n_\tau}\circ \phi$.

Recall that $\Psi^0_\sigma$ is defined over $\mathbb{F}$ and that $g$ is an $\mathbb{F}$-point, hence
\begin{align*}
\tau. \Psi_\sigma(g) &=\tau. \left(\Psi^0_\sigma(\phi(g))\right)
= \Psi^0_\sigma(\tau.(\phi(g)))
=\Psi^0_\sigma(\phi^\tau(g))\\
&=\Psi^0_\sigma(n_\tau \phi(g) n_\tau^{-1})=\Psi_\sigma(\phi^{-1}(n_\tau) g \phi^{-1}(n_\tau)^{-1})\\
&=\Nrd\left(\sum_{i=1}^n (\sigma.e_i) \phi^{-1}(n_\tau) g \phi^{-1}(n_\tau)^{-1} e_i \right) \Nrd(g)^{-1}\\
&=\Nrd\left(\sum_{i=1}^n \left((\sigma\tau^{-1}).e_i\right)  g  (\tau^{-1}.e_i) \right) \Nrd(g)^{-1}\\
&=\Nrd\left(\sum_{i=1}^n  \left((\tau\sigma\tau^{-1}).e_i\right) \cdot g \cdot e_i \right) \Nrd(g)^{-1}
=\Psi_{\tau\sigma\tau^{-1}}(g)
\end{align*}
\end{proof}

Unlike the case of $\mathbf{PGL}_2$ the canonical generators of rational points are not necessarily rational, except for\footnote{This has been confirmed using the SageMath mathematical software \cite{Sage}.} $\Psi_\mathrm{id}$. Nevertheless, we understand rather well their Galois orbits. Specifically, when $\mathfrak{G}\cong S_n$ each $\Psi_\sigma(g)$ belongs to the subfield of $\mathbb{L}$ which corresponds to the centralizer subgroup $Z_{S_n}(\sigma)$ and the Galois orbits of the canonical generators correspond to conjugacy classes in $S_n$. Thus the number of those orbits is the partition number $p(n)$ which for large $n$ is approximately $\frac{1}{4n\sqrt{3}}\exp\left(\pi\sqrt{\frac{2n}{3}}\right)$. Of course, we have in addition numerous algebraic relations between the generators. We now combine those algebraic relations with Proposition \ref{prop:psi-galois}. This a place where 2-transitivity of the Galois group is heavily used.

\begin{corollary}\label{cor:psi-galois-2}
Assume that $\mathfrak{G}$ is 2-transitive. Let $g\in\mathbf{G}(\mathbb{F})$ and let $\sigma\in W$ be an element of the Weyl group without fixed points. If $\Psi_\sigma(g)=0$ then for every $\tau\neq\mathrm{id}$ we have $\Psi_\tau(g)=0$ and $\Psi_\mathrm{id}(g)=1$. Therefore the image of $g$ in 
$\dfaktor{\mathbf{T}}{\mathbf{G}}{\mathbf{T}}$ is $\mathbf{T}e\mathbf{T}$.
\end{corollary}
\begin{proof}
For every $\tau\neq \mathrm{id}$ we look at $\mathscr{C}_\tau=\left\{\upsilon \tau \upsilon^{-1} \mid \upsilon \in \mathfrak{G} \right\}\subset W\cong S_n$. The set $\mathscr{C}_\tau$ has a complete set of roots in the sense of Proposition \ref{prop:relations-zero}, hence by the same proposition there exists $\upsilon\in\mathfrak{G}$ such that
$\Psi_{\upsilon \tau \upsilon^{-1}}(g)=0$. But $\Psi_\tau(g)$ is Galois conjugate to $\Psi_{\upsilon \tau \upsilon^{-1}}(g)$ so $\Psi_\tau(g)=0$ as well.  

To calculate $\Psi_\mathrm{id}(g)$ we use the relation between the canonical generators coming from the Leibniz formula 
\begin{equation*}
\sum_{\tau\in W} \Psi_\tau(g)=1
\end{equation*}
Because $\Psi_\tau(g)=0$ for $\tau\neq \mathrm{id}$ we have $\Psi_\mathrm{id}(g)=1$.

By now we have proved the $\Psi_\tau(g)=\Psi_\tau(e)$ for all $\tau\in W$. Because the $\Psi_\tau$'s generate the whole ring of regular functions on $\left(\dfaktor{\mathbf{T}}{\mathbf{G}}{\mathbf{T}}\right)_\mathbb{L}$ the image of $g$ as an $\mathbb{L}$-point in $\left(\dfaktor{\mathbf{T}}{\mathbf{G}}{\mathbf{T}}\right)_\mathbb{L}$  is equal to the image of $e$. Yet because the $\mathbb{F}$-points of a variety naturally \emph{inject} into the $\mathbb{L}$-points, we conclude that the image of $g$ in $\dfaktor{\mathbf{T}}{\mathbf{G}}{\mathbf{T}}$ is already equal to the image of the identity.
\end{proof}

\section{Local Integral Structure of Double Torus Quotient}\label{sec:packets-csa}
Our goal is to study the integral properties of the canonical generators evaluated at a point associated to two torus orbits coming from the same packet. 

The type of results we are about to prove are the simplest to state when $\mathbb{F}=\mathbb{Q}$. Then if $\delta_L \mathbf{T}(\mathbb{A}) g_\mathbb{A}$, $\delta_R  \mathbf{T}(\mathbb{A}) g_\mathbb{A}$ with $\delta_R,\delta_L\in\mathbf{G}(\mathbb{Q})$ are two representatives of the same homogeneous toral set then $D_f \Psi_\sigma({g_L}^{-1} g_R)$ is integral, where $D_f$ is the product of the local archimedean discriminants of the homogeneous toral set. For a general number field $\mathbb{F}$ we need to replace $D_f$ by a suitable ideal of $\mathcal{O}_\mathbb{F}$.

These results are crucial, for they allow us to apply rudimentary geometry of numbers later on. This is akin to the simple fact that, in rank 1, the discriminant inner product of two binary \emph{integral} quadratic forms associated with two periodic torus orbits is integral.

\subsection{The integral structure}\label{sec:global-order}
Let $\mathbf{G}$ be as in the previous section, viz.\  $\mathbf{G}$ is isogenous to $\mathbf{PGL}_1(\mathbf{B})$ for a degree $n$ central simple algebra $\mathbf{B}$ over $\mathbb{F}$. We  fix once and for all an $\mathcal{O}_\mathbb{F}$-order $\Omega\subset\mathbf{B}(\mathbb{F})$. Let $\Omega_u\subset\mathbf{B}(\mathbb{F}_u)$ be the $u$-adic closure of $\Omega$, then $\Omega_u$ is a maximal order for almost all places $u$. Define also the compact-open subgroups $K_u<\mathbf{G}(\mathbb{F}_u)$ using the order $\Omega$
\begin{equation*}
K_u\coloneqq \mathbf{Z}(\Fu)\mathbf{PGL}_1(\Omega_u)
\end{equation*}
where $\mathbf{Z}$ is the finite center of $\mathbf{G}$.

\subsection{Homogeneous Toral Sets for Central Simple Algebras}
Let $\mathbf{T}<\mathbf{G}$ be a maximal torus define and anisotropic over $\mathbb{F}$ and $g_\mathbb{A}\in\mathbf{G}(\mathbb{A})$.
This data fixes a homogeneous toral set
\begin{equation*}
Y=\mathbf{G}(\mathbb{F})\mathbf{T}(\mathbb{A})g_\mathbb{A}
\end{equation*}
Recall that the isogeny
\begin{equation*}
\mathbf{SL}_1(\mathbf{B})\to\mathbf{G}\to\mathbf{PGL}_1(\mathbf{B})
\end{equation*}
restricts to an isogeny of maximal tori
\begin{equation*}
\mathbf{SL}_1(\mathbf{E})\to\mathbf{T}\to\mathbf{PGL}_1(\mathbf{E})
\end{equation*}
where the closed subvariety $\mathbf{E}<\mathbf{B}$ represents an embedding $\iota\colon \mathbb{K}\hookrightarrow\mathbf{B}(\mathbb{F})$ of a degree $n$ field extension $\mathbb{K}/\mathbb{F}$.

A homogeneous toral set ${\mathbf{G}(\mathbb{F})}{\mathbf{T}(\mathbb{A})g_\mathbb{A}}$, $g_\mathbb{A}=(g_u)_{u\in \mathcal{V}_\mathbb{F}}$, defines for all $u\in \mathcal{V}_\mathbb{F}$ a torus over $\Fu$
\begin{equation*}
\mathbf{H}_u\coloneqq{g_u}^{-1} \mathbf{T}_{\Fu} g_u
\end{equation*}
The homogeneous toral set is invariant under the geometric torus $\mathbf{H}_u(\Fu)$ at each rational place $u$. Denote by $\Psi^\mathbf{T}_\sigma$ the canonical generators $\dfaktor{\mathbf{T}}{\mathbf{G}}{\mathbf{T}}$ over a splitting field and by $\Psi^{\mathbf{H}_u}_\sigma$ the canonical generators for $\dfaktor{\mathbf{H}_u}{\mathbf{G}_{\Fu}}{\mathbf{H}_u}$ over a splitting field. The equality  $\Psi^\mathbf{T}_\sigma(\lambda)=\Psi^{\mathbf{H}_u}_\sigma({g_u}^{-1} \lambda g_u)$ holds in any field extension of $\Fu$ splitting the torus $\mathbf{T}$.

\emph{Throughout this section the homogeneous toral set $Y$ is fixed.}

\subsection{The associated order}
Define
\begin{equation*}
\mathbf{C}_u = g_u^{-1}\mathbf{E}_{\mathbb{F}_u} g_u
\end{equation*}
where $\mathbf{E}_{\mathbb{F}_u}$ is the base change of $\mathbf{E}$ to $\mathbb{F}_u/\mathbb{F}$. In particular, $\mathbf{C}_u$ is defined over $\mathbb{F}_u$ and the isogeny \eqref{eq:fixed-isogenies} induces an isogeny of tori
\begin{equation*}
\mathbf{SL}_1(\mathbf{C}_u)\to  \mathbf{H}_u \to\mathbf{PGL}_1(\mathbf{C}_u)
\end{equation*}
The commutative algebra $\mathbf{C}_u(\mathbb{F}_u)$ is a degree $n$ \'etale-algebra over $\mathbb{F}_u$ satisfying
\begin{equation*}
\mathbf{C}_u(\mathbb{F}_u)\simeq \mathbb{K}_u=\mathbb{K}\otimes \mathbb{F}_u\simeq \prod_{v\mid u} \mathbb{K}_v
\end{equation*}
We now define a local order $\mathcal{R}_u<\mathbf{C}_u(\mathbb{F}_u)$ by
\begin{equation*}
\mathcal{R}_u=\mathbf{C}_u(\mathbb{F}_u)\cap \Omega_u
\end{equation*}
Notice that $g_u$ maps to $\mathbf{PGL}_1(\Omega_u)$ at almost all places $u$, hence $R_u$, as an order in $\mathbb{K}_u$, coincides with the completion of $\Omega\cap \mathbf{E}(\mathbb{F})$ at almost all places. The latter is an order in $\mathbf{E}(\mathbb{F})\simeq \mathbb{K}$, hence $R_u$ coincides with the unique maximal order in $\mathbb{K}_u$ for almost all $u$. We can then take the intersection
\begin{equation*}
\mathcal{R}\coloneqq \cap_{u<\infty} \mathcal{R}_u
\end{equation*}
inside $\mathbb{K}$. This intersection is again an order in $\mathbb{K}$ that differs from the naive version $\Omega\cap \mathbf{E}(\mathbb{F})$ at the places where $g_\mathbb{A}$ is not integral. We will be chiefly interested in homogeneous toral set such that the attached order $\mathcal{R}$ is the unique maximal order in $\mathbb{K}$. We turn to discuss the relation between the discriminant of $\mathcal{R}$ and the discriminant of the homogeneous toral set. 

\subsection{Discriminant Data}\label{sec:discriminant-data}
For the cases  we study we can provide a more concrete description of the discriminant data. It will be useful in studying the arithmetic properties of the canonical generators.

\subsubsection{Local Nonarchimedean Discriminant.}

\begin{definition}\label{defi:realtive-local-discriminant}
The \emph{relative} local discriminant of $\mathcal{R}_u$ is
\begin{equation*}
\mathscr{D}_u\coloneqq\det(\Trd(f_i f_j))_{1\leq i,j \leq n}
\end{equation*}
where $f_1,\ldots,f_n$ is an $\mathcal{O}_\Fu$ basis for $\mathcal{R}_u$. We could have chosen a different base for $\mathcal{R}_u$ and get a different value for the relative local discriminant, but all those values would differ only by a square of a unit \cite[\S 7.6]{CasselsLocalFields}. Hence the relative local discriminant is a well defined class in $\faktor{\mathcal{O}_\Fu}{{\mathcal{O}_\Fu^\times}^2}$. We fix some representative for $\mathscr{D}_u$.
\end{definition}

The local discriminant $D_u$ of $Y$ whose definition appears in \S \ref{sec:disc} is almost equal to the \emph{absolute} discriminant of  $\mathcal{R}_u$, $|\mathscr{D}_u|_u^{-1}$ \cite[Lemma 9.5]{ELMVCubic}. The two notions differ by a constant factor which for almost all places $u$ is equal to $1$. \emph{We take the liberty to ignore this non-consequential subtlety and redefine the local discriminant of a homogeneous toral set to be equal hereon to the absolute discriminant of $\mathcal{R}_u$, i.e.\ $D_u\coloneqq|\mathscr{D}_u|_u^{-1}$}.

\subsubsection{Local Archimedean Discriminant}\label{sec:local-arch-disc}
Recall that in order to calculate the local discriminant at an archimedean place $u$ we had a choice of a norm on $\left(\bigwedge^r\mathfrak{g}(\Fu)\right)^{\otimes 2}$ with $r=n-1$. We restrict ourselves to a norm of a specific type which will be comfortable later on.

Set $\overline{\Fu}$ to be the algebraic closure of $\Fu$, as $\overline{\Fu}\cong\mathbb{C}$ we identify them both and work over $\mathbb{C}$. A norm on $\left(\bigwedge^r\mathfrak{g}(\overline{\Fu})\right)^{\otimes 2}$ naturally restricts to a norm on $\left(\bigwedge^r\mathfrak{g}(\Fu)\right)^{\otimes 2}$.

We have a decomposition of $\mathbb{C}$-vector spaces
\begin{equation*}
\mathbf{B}(\mathbb{C})=\mathbf{B}^0(\mathbb{C})\oplus \Id\mathbb{C}
\end{equation*}
Where $\mathbf{B}^0\coloneqq\ker\Trd$ is the trace zero part of the algebra.

Let $Q$ be a Hermitian inner product on $\mathbf{B}(\mathbb{C})$  such that the spaces $\mathbf{B}^0(\mathbb{C})$ and $\Id\mathbb{C}$ are orthogonal to each other.

This in turn induces an inner product on $\mathfrak{g}(\mathbb{C})$ using the identification of $\mathfrak{g}(\mathbb{C})$ with the trace zero part $\mathbf{B}^0(\mathbb{C})$.

This inner product can be extended to an inner product on $\bigwedge^r\mathfrak{g}(\mathbb{C})$ using the standard determinant construction. This last inner product is a bilinear form on $\bigwedge^r\mathfrak{g}(\mathbb{C})$ so it extends to a linear functional on $\left(\bigwedge^r\mathfrak{g}(\mathbb{C})\right)^{\otimes 2}$. The absolute value of this linear functional will be the required norm.

As before, $\mathfrak{t}$ is the Lie algebra of $\mathbf{T}$.
If we choose a base $f_1,\ldots,f_r$ for $\Ad(g_u^{-1})\mathfrak{t}(\Fu)$ we can extend it to a base $f_0=1,f_1,\ldots,f_{n-1=r}$ of $\mathbf{C}_u(\Fu)$ by identifying  $\Ad(g_u^{-1})\mathfrak{t}(\Fu)$ with the trace zero part of $\mathbf{C}_u(\Fu)$. For this base we can calculate the discriminant in the following way
\begin{align}
D_u&=\det(Q(f_i,f_j))_{1\leq i,j \leq r}\cdot \left| \det(\Trd(f_i f_j))_{1\leq i,j \leq r} \right|^{-1}
\label{eq:arch-disc-lie}\\
&=\det(Q(f_i,f_j))_{0\leq i,j \leq r}\cdot \left| \det(\Trd(f_i f_j))_{0\leq i,j \leq r} \right|^{-1}
\label{eq:arch-disc-alg}
\end{align}
Expression \eqref{eq:arch-disc-lie} involves determinants of $f_1,\ldots,f_r$ and expression \eqref{eq:arch-disc-alg} involves determinant of $f_0,f_1,\ldots,f_r$. The second expression is independent of the base we have chosen for $\mathbf{C}_u(\mathbb{C})$. Hence for $Q$ of the form we are discussing, we can use \eqref{eq:arch-disc-alg} as the definition of the archimedean discriminant with any base we like for $\mathbf{C}_u(\mathbb{C})$.
\subsubsection{Global Discriminant.}
\paragraph{Matrix algebra over $\mathbb{Q}$}
We start with describing the global discriminant in the classical setting of $\mathbf{M}_n$ over $\mathbb{F}=\mathbb{Q}$.

The following description of the global discriminant from \cite{ELMVPeriodic,ELMVCubic} is important for applications.

Let $\uhat\in S$ be a fixed place. Fix as well a split torus $\mathbf{H}^0$ over $\Quhat$.  

We consider only homogeneous toral sets for which $\mathbf{H}_\uhat=\mathbf{H}^0$ and $\mathbf{H}_u=\mathbf{T}_u$ for all $u\neq\uhat$. Notably, $\mathbf{T}$ is assumed to be split over $\Quhat$. We call such homogeneous toral sets \emph{simple}.

Given $\mathbf{T}$ there are no more than finitely many simple homogeneous toral sets corresponding to $\mathbf{T}$. They are parametrized by the absolute Weyl group of $\mathbf{T}$. 

Notice that the general definition of the local discriminant depends only on the tori $\mathbf{H}_u$. In our concrete case all these tori either come from $\mathbf{T}$ or are all equal to $\mathbf{H}^0$ at the place $\uhat$. Hence the local discriminant at the place $\uhat$ does not depend on the simple homogeneous toral set. Its contribution to the global discriminant is immaterial.

On the other hand at all places $u\neq\uhat$ the local discriminant is a property of $\mathbf{T}$. Accordingly, the global discriminant of such a homogeneous toral set essentially depends only on $\mathbf{T}$.

As expected, the global discriminant in this case is actually an invariant of a global object. To such a homogeneous toral set one can attach a global order $\mathcal{R}\coloneqq \mathbb{K} \cap \Omega$, where $\mathbb{K}=\mathbf{E}(\mathbb{Q})$ and $\mathbf{E}$ is the maximal commutative subalgebra over $\mathbf{T}$. Recall that $\mathbb{K}$ is a degree $n$ field extension of $\mathbb{Q}$ and $\mathcal{R}$ is an order in this number field.

It is shown in \cite[\S 6.3]{ELMVCubic} that the discriminant of the order $\mathcal{R}$ coincides up to a constant factor with the global discriminant of the homogeneous toral set. This equality of discriminants follows from the fact that the discriminant of an order can be calculated locally as well\footnote{The place $\uhat$ splits in $\mathbb{K}$ so the local discriminant of the order $\mathcal{R}$ at $\uhat$ is trivial.}.

\paragraph{Relative non-archimedean discriminant.}
We return to our general setting of a central simple algebra $\mathbf{B}$ over a number field $\mathbb{F}$ with $Y$ a general homogeneous toral set.

We begin with the definition of the \emph{absolute} non-archimedean discriminant.
\begin{definition}
We define the \emph{absolute non-archimedean discriminant} of the homogeneous toral set to be $D_f=\prod_{u\in \mathcal{V}_{f}} D_u$ where the product is taken across all the non-archimedean places.
\end{definition}

There is a notion of a \emph{relative} global non-archimedean discriminant which will be useful to us. This notion refines the global discriminant of simple homogeneous toral sets for $\mathbf{M}_n$ over $\mathbb{Q}$.

\begin{definition}\label{defi:realtive-global-discriminant}
The \emph{relative non-archimedean discriminant} of the homogeneous toral set is the unique \emph{ideal} $\mathscr{D}$ of $\mathcal{O}_\mathbb{F}$ such that for any $u\in\mathcal{V}_{\mathbb{F},f}$ the closure of $\mathscr{D}$ in $\mathcal{O}_\Fu$ is $\mathscr{D}_u \mathcal{O}_\Fu$.
\end{definition}

One can see using the relation between the norm of an ideal and non-archimedean valuations that
\begin{equation*}
{\Nr}_{\mathbb{F}/\mathbb{Q}}\mathscr{D}=\prod_{u\in\mathcal{V}_{\mathbb{F},f}} |\mathscr{D}_u|_u^{-1}=\prod_{u\in\mathcal{V}_{\mathbb{F},f}} D_u=D_f
\end{equation*}
In particular, if $\mathbb{F}$=$\mathbb{Q}$ then $\mathscr{D}=D_f \mathbb{Z}$.

\subsection{Denominators of Canonical Generators}
Recall the we have fixed a homogeneous toral set  $Y\coloneqq{\mathbf{G}(\mathbb{F})}{\mathbf{T}(\mathbb{A})g_\mathbb{A}}$.

Let $\delta_L\mathbf{T}(\mathbb{A})g_\mathbb{A}$, $\delta_R\mathbf{T}(\mathbb{A})g_\mathbb{A}$ with $\delta_L, \delta_R\in\mathbf{G}(\mathbb{F})$ be two representatives of the homogeneous toral set $Y$. We are interested in what can be said about $$\lambda\coloneqq {\delta_L}^{-1}\delta_R$$ when $\delta_R\mathbf{T}(\mathbb{A})g_\mathbb{A}$ belongs to a small neighborhood around $\delta_L\mathbf{T}(\mathbb{A})g_\mathbb{A}$. Specifically, given a  finite place $u$, we look at pairs $(\delta_L,\delta_R)\in\mathbf{G}(\mathbb{F})^{\times 2}$ such that 
\begin{equation*}
\delta_R\mathbf{T}(\Fu)g_u \subset \delta_L\mathbf{T}(\Fu)g_u K_u
\end{equation*}
Or equivalently
\begin{align*}
\lambda={\delta_L}^{-1}\delta_R
&\in \mathbf{G}(\mathbb{F})\cap \prod_{u\in\mathcal{V}_{\mathbb{F},f}}\mathbf{T}(\Fu) g_u K_u {g_u}^{-1} \mathbf{T}(\Fu)\\
&=\mathbf{G}(\mathbb{F})\cap \prod_{u\in\mathcal{V}_{\mathbb{F},f}} g_u \mathbf{H}_u(\Fu) K_u \mathbf{H}_u(\Fu) {g_u}^{-1} 
\end{align*}

We would like to compute a denominator for $\Psi^\mathbf{T}_\sigma(\lambda)\in\mathbb{L}$ in each field extension $\mathbb{L}_w$, when $w$ is a non-archimedean place of $\mathbb{L}$. By a denominator we mean a $w$-adic integer $d_w\in\mathcal{O}_{\mathbb{L}_w}$ such that $d_w \Psi^\mathbf{T}_\sigma(\lambda) \in\mathcal{O}_{\mathbb{L}_w}$. Examining the rank 1 case one observes that a natural candidate for this denominator is the local discriminant of the torus. Actually, it will be the local \emph{relative} discriminant.

Recall that as $\lambda$ is an $\mathbb{F}$-point and the polynomials $\Psi^\mathbf{T}_\sigma$ for $\sigma\in W$ are defined over $\mathbb{L}$ we have that $\Psi^\mathbf{T}_\sigma(\lambda)\in\mathbb{L}$. For a place $w$ of $\mathbb{L}$ extending $u$ we fix a commutative diagram of embeddings
\begin{center}
\begin{tikzcd}
                                                                 & \mathbb{L}  \arrow[hookrightarrow]{dr}  & \\
\mathbb{F} \arrow[hookrightarrow]{ur} \arrow[hookrightarrow]{dr} &                                         &  \mathbb{L}_w\\
                                                                 & \Fu \arrow[hookrightarrow]{ur} & \\
\end{tikzcd}
\end{center}
This induces the following commutative diagram of natural injections for points in $\mathbf{G}$
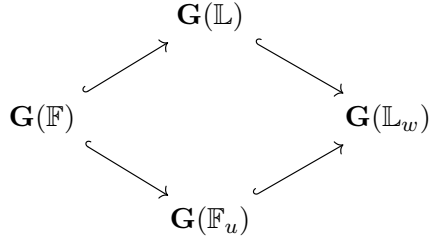
\begin{center}
\begin{tikzcd}
                                                                             & \mathbf{G}(\mathbb{L})  \arrow[hookrightarrow]{dr}  & \\
\mathbf{G}(\mathbb{F}) \arrow[hookrightarrow]{ur} \arrow[hookrightarrow]{dr} &     &  \mathbf{G}(\mathbb{L}_w)\\
                                                                             & \mathbf{G}(\Fu) \arrow[hookrightarrow]{ur} & \\
\end{tikzcd}
\end{center}
We can write in $\mathbf{G}(\Fu)$: ${g_u}^{-1}\lambda g_u=h^\mathrm{L}_u k_u h^\mathrm{R}_u$ with $k_u\in K_u$ and $h^\mathrm{L}_u,h^\mathrm{R}_u \in \mathbf{H}_u(\Fu)$. 
These elementary facts imply that in $\mathbb{L}_w$ one has 
\begin{equation*}
\Psi^\mathbf{T}_\sigma(\lambda)=\Psi^{\mathbf{H}_u}_\sigma({g_u}^{-1}\lambda g_u)=
\Psi^{\mathbf{H}_u}_\sigma(h^\mathrm{L}_u k_u h^\mathrm{R}_u)=\Psi^{\mathbf{H}_u}_\sigma(k_u)
\end{equation*}
The last equality holds because  $\Psi^{\mathbf{H}_u}_\sigma$ is invariant under the left and right actions of ${\mathbf{H}_u}$. Hence for $u$ non-archimedean it is enough to consider expressions of the form $\Psi^{\mathbf{H}_u}_\sigma(k_u)$. 

\begin{proposition}\label{prop:denom}
Let $u$ be a non-archimedean place of $\mathbb{F}$ and let $w$ be a place of $\mathbb{L}$ above $u$. Assume that $\mathbf{B}$ is \emph{unramified} over $u$, equivalently $\mathbf{G}$ is $u$-split.  Assume that $\mathcal{R}_u$ is the \emph{maximal order} of $C_u\coloneqq \mathbf{C}_u(\Fu)$. Then for all $\sigma\in W$ and for all $k_u\in K_u$ we have
\begin{equation*}
\mathscr{D}_u \Psi^{\mathbf{H}_u}_\sigma(k_u) \in \mathcal{O}_{\mathbb{L}_w}
\end{equation*}
\end{proposition}
\begin{proof} 
First notice that it is enough to prove the claim for $\mathbf{G}=\mathbf{PGL}_1(\mathbf{B})$ and $K_u=\mathbf{PGL}_1(\Omega_u)$. The general case would follow from the fact that the isogeny sends $K_u$ to a subset of $\mathbf{PGL}_1(\Omega_u)$.

\paragraph{Sketch.}
We begin by providing a sketch of the proof for $\mathbf{PGL}_1(\mathbf{B})$. The idea is simple but the proof is cluttered with technical details. 
\begin{enumerate}
\item Without loss of generality we can assume that $\Omega_u$ is a maximal order.
\item\label{sketch1} Over a non-archimedean place $u$ of $\mathbb{F}$ where $\mathbf{B}$ is unramified the central simple algebra $\mathbf{B}(\Fu)$ is isomorphic to a matrix algebra over $\Fu$.

\item The isomorphism from item \ref{sketch1} can be chosen so that the image of $\Omega_u$ is exactly the matrix ring over $\mathcal{O}_{\Fu}$.

\item The basic idea is to find a \emph{nice} element $g\in\mathbf{GL}_n(\mathbb{L}_w)$ conjugating the standard diagonal commutative subalgebra of the matrix algebra to $C_w\coloneqq \mathbf{C}_u(\mathbb{L}_w)$ and then transform $\Psi^{\mathbf{H}_u}_\sigma(k_u)$ to an expression involving primitive orthogonal idempotents over the standard diagonal torus, $e_1^0,\ldots,e_n^0$. 

\item The expression for the canonical generator becomes 
\begin{equation*}
\Psi^{\mathbf{H}_u}_\sigma(k_u)=\Nrd\left(\sum_{i=1}^n \sigma.e_i^0 \, \left(g^{-1} k_u g \right) \, e_i^0\right) \cdot \Nrd({k_u}^{-1})
\end{equation*}
\item Both $e_i^0$ and $\sigma.e_i^0$ are in $\Omega_w$ which is the matrix ring over $\mathcal{O}_{\mathbb{L}_w}$. The same holds for $k_u\in K_u\subset \Omega_u\subset \Omega_w$.

\item \emph{The heart of the argument:} We are able to choose $g$ so that both $g^{-1}$ and $\Delta g$ are in the order $\Omega_w$, where $\Delta^{-1}\in C_u$ generates the different ideal of $\mathcal{R}_u$ which is principal. The norm of $\Delta$ contributes the $\mathscr{D}_u$ factor to the end result.
\end{enumerate}

The choice of a nice $g$ as above is done in the following way

\begin{enumerate}
\item In matrix algebras the embedding of rational tori with a given $\mathcal{R}_u$ are in bijection with homothety classes of bases for proper fractional ideals of $\mathcal{R}_u$.

\item Since $\mathcal{R}_u$ is maximal it is a PID, and all the embeddings correspond to bases of the trivial fractional ideal. This translates directly to the fact that we can choose $g$ as above such that the rows of $g^{-1}$ are a base for $\mathcal{R}_u$. In particular, $g^{-1}$ is integral.

\item The rows of $\tensor[^t]{{g}}{}$, equivalently the columns of $g$, will be a dual base with respect to the reduced trace to the rows of $g^{-1}$; hence they span the inverse different. To transform the matrix $g$ into an integral one, we need to multiply it by the generator of the different ideal.
\end{enumerate}

\paragraph{Structure of the Algebra over $\Fu$.}
The order $\Omega_u$ is contained in a maximal order $\Omega^{\max}_u\supseteq\Omega_u$. One has that $K_u$ is a subset of $\mathbf{PGL}_1(\Omega^{\max}_u)$. 

The order $\mathcal{R}_u$ is contained in the order $\Omega^{\max}_u\cap C_u$, where $C_u=\mathbf{C}_u(\Fu)$. As $\mathcal{R}_u$ was assumed to be a maximal order this is actually an equality \footnote{This is not the significant use of the maximality of $\mathcal{R}_u$.}. We see that the assumptions of the theorem would hold with $\Omega_u$ replaced by $\Omega^{\max}_u$ and that its conclusion for $\Omega_u$ follows from that for $\Omega^{\max}_u$. We assume without loss of generality that $\Omega_u=\Omega^{\max}_u$.

As $B_u$ is unramified it is isomorphic to the matrix algebra $\mathbf{M}_n(\Fu)$. Moreover, by \cite[Theorem 3.5 and Theorem 3.8]{AuslanderGoldman} we can choose the isomorphism $\iota: B_u \xrightarrow{\sim}\mathbf{M}_n(\Fu)$ so that $\iota(\Omega_u)=\mathbf{M}_n(\mathcal{O}_{\Fu})$ (see also \cite[Theorem 17.3]{Reiner}). By abuse of notation we will identify $B_u$ with $\mathbf{M}_n(\Fu)$ and $\Omega_u$ with its image.

Let $B_w \coloneqq B_u\otimes_{\Fu} \mathbb{L}_w \cong \mathbf{M}_n(\mathbb{L}_w)$ and let $\Omega_w\coloneqq \mathbf{M}_n(\mathcal{O}_{\mathbb{L}_w})$. The ring $\Omega_w$ is a maximal order in $B_w$ and $\Omega_w \cap B_u=\Omega_u$. 

Let $\mathrm{Diag}_w$ be the diagonal commutative algebra in $\mathbf{M}_n(\mathbb{L}_w)$ and set $A_w\coloneqq\mathbf{PGL}_1(\mathrm{Diag}_w)$.  Write $T_w\coloneqq \mathbf{T}(\mathbb{L}_w)$. The tori $T_w$ and $A_w$ are both split, hence there exists $g_0\in \mathbf{G}(\mathbb{L}_w)$ such that $T_w =g_0 A_w {g_0}^{-1}$. It follows that $C_w=g_0 \mathrm{Diag}_w {g_0}^{-1}$, where $C_w\coloneqq \mathbf{C}(\mathbb{L}_w)$. Any element of $g_0 \, \Nrml_{\mathbf{GL}_n(\mathbb{L}_w)}(\mathrm{Diag}_w)$ will also conjugate those two commutative algebras to each other. Our
aim now is to find a conjugating element whose matrix rows are a base for a proper fractional ideal $\mathcal{I}$ of $\mathcal{R}_u$ in a suitable sense.

\paragraph{Fractional Ideals.}
The explicit definition of a proper fractional ideal of $\mathcal{R}_u$ is as follows. A full $\mathcal{O}_{\Fu}$-lattice $\mathcal{J}\subseteq C_u$ is an $\mathcal{R}_u$ fractional ideal if $\mathcal{R}_u \mathcal{J}\subseteq \mathcal{J}$. The ring associated to a $\mathcal{J}$ is 
\begin{equation}
\mathcal{O}_\mathcal{J} \coloneqq \left\{x\in C_u \mid x\mathcal{J} \subseteq \mathcal{J}  \right\}
\end{equation}
The fractional ideal is proper for $\mathcal{R}_u$ if $\mathcal{O}_\mathcal{J}=\mathcal{R}_u$.

The \'{e}tale-algebra $C_u$ is a product of fields\footnote{When $\mathbf{H}_u$ is an extension of scalars of an anisotropic $\mathbb{F}$ torus then $\mathbb{K}\coloneqq\mathbf{C}(\mathbb{F})$ is a field. Weak approximation for number fields implies $C_u\cong \prod_{v|u} \mathbb{K}_v$ where $v$ runs over the places of $\mathbb{K}$ over $u$.} $C_u\cong\prod_{v\in \mathcal{W}_u} \mathbb{K}_v$.
The fields $\mathbb{K}_v$ for\footnote{We can consider $\mathcal{W}_u$ as a set of equivalence classes of valuation on finite field extensions of $\Fu$.} $v\in\mathcal{W}_u$ are finite extensions of $\Fu$ and $\sum_{v\in\mathcal{W}_u} \left[\mathbb{K}_v \colon \Fu\right]=n$.

Let $1_v\in C_u$ be the identity element of $\mathbb{K}_v$. Obviously $e=\sum_{v\in\mathcal{W}_u} 1_v$. 
By the maximality assumption for $\mathcal{R}_u$ we have $\mathcal{R}_u\cong\prod_{v\in\mathcal{W}_u} \mathcal{O}_{\mathbb{K}_v}$. We can use this decomposition of $\mathcal{R}_u$ to decompose any $\mathcal{R}_u$ fractional ideal $\mathcal{J}$ into a sum of $\mathcal{O}_{\mathbb{K}_v}$ fractional ideals $\mathcal{J} \cong \oplus_{v\in\mathcal{W}_u} \mathcal{J}_v$, where $\mathcal{J}_v\coloneqq 1_v \mathcal{J}$. It holds that $\mathcal{J}$ is $\mathcal{R}_u$-proper if and only if $\mathcal{J}_v$ is $\mathcal{O}_{\mathbb{K}_v}$ proper for all $v\in\mathcal{W}_u$.

If $\mathcal{J}$ is proper then each $\mathcal{J}_v$ is a proper fractional ideal of the Discrete Valuation Ring\footnote{This is one of the significant uses of the maximality assumption.} $\mathcal{O}_{\mathbb{K}_v}$, hence it is principal and invertible. We conclude that $\mathcal{J}$ itself is principal and invertible. One can choose $a_\mathcal{J}\in C_u^\times$ such that $\mathcal{J}= a_\mathcal{J} \mathcal{R}_u$.

For each fractional ideal $\mathcal{J}$ one can define the dual fractional ideal 
\begin{equation}
\widecheck{\mathcal{J}}=\left\{x\in C_u \mid \Trd(x \mathcal{J}) \subseteq \mathcal{O}_{\Fu} \right\}
\end{equation}
The dual is proper if $\mathcal{J}$ is. In particular, the local inverse different $\widecheck{\mathcal{R}_u}$ is proper and principal\footnote{Second significant use of the maximality assumption.}. We write $\widecheck{\mathcal{R}_u}=\Delta^{-1}\mathcal{R}_u$ for $\Delta\in\mathcal{R}_u^\times$. Essentially by definition, $\Delta$ can be chosen so that $\Nrd(\Delta)=\mathscr{D}_u$.

\paragraph{Associated Fractional Ideal}
We follow the proof of \cite[Corollary 4.4]{ELMVPeriodic} to show that a conjugating element can be chosen so its rows are in a suitable manner a base for a proper fractional ideal. 

We start by constructing an $\Fu$-vector space isomorphism $\jmath\colon C_u\to\mathbb{F}_u^n$. For each $v\in\mathcal{W}_u$ we have a linear endomorphism of $\mathbb{F}_u^n$ given by $a\mapsto 1_v\cdot a$, where the right hand side is the multiplication of a matrix by a vector. Because $1_v\neq 0$ for each $v$, the kernel of this endomorphism is a proper linear subspace of $\mathbb{F}_u^n$. The union of a finite collection of proper subspaces never exhausts a vector space over characteristic 0, hence there exists $a\in\mathbb{F}_u^n$ such that $1_v\cdot a\neq 0$ for all $v\in\mathcal{W}_u$. Let $\jmath(x)\coloneqq x\cdot a$. Notice that $\jmath$ intertwines the action of $C_u$ on itself by matrix-matrix multiplication with the action of $C_u$ on $\mathbb{F}_u^n$ by matrix-vector multiplication.

To see that $\jmath$ is an isomorphism we prove it has a trivial kernel. For each $0\neq x\in C_u$ one has $x=\sum_{v\in\mathcal{W}_u} x 1_v$ and $x 1_v$ can be identified with the $\mathbb{K}_v$ component of $x$. In particular, there exists $v_0$, depending on $x$, such that $x 1_{v_0}\neq 0$. Let $y_{v_0}\in C_u$ be the inverse of $x 1_{v_0}$ in $\mathbb{K}_{v_0}$, then $y_{v_0} x =1_{v_0}$. If $0\neq x\in\ker\jmath$ then $x\cdot a=0$ and $0=y_{v_0} x\cdot a =1_{v_0} \cdot a$ which is a contradiction.

The fractional ideal we seek is the $\mathcal{O}_{\Fu}$-lattice $\mathcal{I}\coloneqq\jmath^{-1}(\mathcal{O}_{\Fu}^n)$. Denote the canonical base elements for $\mathbb{F}_u^n$  by $b_i=\left(0,\ldots,0,1,0,\ldots,0\right)$ with $1$ in the $i$ place for $1\leq i \leq n$. Let $\beta_i\coloneqq \jmath^{-1}(b_i)$; this is a base of $\mathcal{I}$.
Notice that for any $x=(x_{i,j})_{1\leq i,j \leq n} \in C_u \subset \mathbf{M}_n(\Fu)$ we have 
\begin{equation*}
x\beta_j= \jmath^{-1}(x b_j)
=\sum_{i=1}^n x_{i,j} \jmath^{-1}(b_i)
=\sum_{i=1}^n x_{i,j} \beta_i
\end{equation*}
This implies immediately that $\mathcal{O}_\mathcal{I}=\mathbf{M}_n(\mathcal{O}_{\Fu}) \cap C_u=\mathcal{R}_u$, i.e.\
$\mathcal{I}$ is a proper $\mathcal{R}_u$ fractional ideal.

\paragraph{The Conjugating Matrix.}
Fix $e_1,\ldots,e_n$ -- an ordered complete set of primitive orthogonal idempotents of $C_w$. 
We embed $\mathcal{I}$ back into $C_w$ by defining an $\mathbb{L}_w$-\emph{algebra} isomorphism $\iota\colon C_w \to\mathbb{L}_w^n$,  given by 
\begin{equation*}
\iota(x)\coloneqq \left(\Trd(x e_1), \Trd(x e_2),\ldots,\Trd(x e_n)\right)
\end{equation*} 
For every $x=(x_{i,j})_{1\leq i,j \leq n} \in C_w$ we can write
\begin{equation*}
\iota(x)\iota(\beta_j)=\iota(x\beta_j)=\sum_{i=1}^n x_{i,j} \iota(\beta_i)
\end{equation*}
Let $M_{\iota(x)}\in \mathbf{GL}_n(\mathbb{L}_w)$ be the \emph{diagonal} matrix in the canonical basis corresponding to the linear transformation  $a\mapsto \iota(x)a$. The equality above implies for any $x\in C_w$ that 
$g_\mathcal{I} M_{\iota(x)}{g_\mathcal{I}}^{-1}=x$, where $g_\mathcal{I}\in  \mathbf{GL}_n(\mathbb{L}_w)$ is the base change matrix defined by $g_\mathcal{I}\, \iota(\beta_i)=b_i$. We conclude that $g_\mathcal{I} \mathrm{Diag}_w {g_\mathcal{I}}^{-1} = C_w$.
 
The ideal $\mathcal{I}$ is principal and invertible because it is proper. Write $\mathcal{I}=a \mathcal{R}_u$ for some $a \in C_u^\times$. The elements $a\beta_1,\ldots,a\beta_n$ form a basis for $\mathcal{R}_u\cong\prod_{v\in\mathcal{W}_u} \mathcal{O}_{\mathbb{K}_v}$, hence they are all integral. As $\iota$ is an algebra isomorphism, we deduce that $r_i\coloneqq\iota(a\beta_i)$ is integral for all $i$, viz.\ $r_i\in \mathcal{O}_{\mathbb{L}_w}^n$. 

Define $g=g_\mathcal{I} M_{\iota(a^{-1})}$ then $g \mathrm{Diag}_w {g}^{-1}=C_w$ and ${g}^{-1}\, b_i=\iota(\alpha\beta_i)=r_i\in \mathcal{O}_{\mathbb{L}_w}^n$. We deduce that ${g}^{-1}\in \Omega_w$.

\paragraph{Integrality of $\Delta g$}
We have shown that ${g}^{-1}$ is integral, the last critical step will be to show the $\Delta g$ is also in $\mathbf{M}_n(\mathcal{O}_{\mathbb{L}_w})$. 

Because of the way we constructed the isomorphism $\iota$, the reduced trace form on $C_w$ is pushed forward by $\iota$ to the regular euclidean form on $\mathbb{L}_w^n$. In particular $\tensor[^t]{{g}}{}$ is the adjoint of $g$ with respect to this bilinear form. Hence if $\widecheck{r_1},\ldots,\widecheck{r_n}$ is the dual base of $r_1,\ldots,r_n$, then $\tensor[^t]{{g}}{}\,b_i=\widecheck{r_i}$.

Notice that $\iota^{-1}(\widecheck{r_1}),\ldots,\iota^{-1}(\widecheck{r_n})$ is the dual base of $a\beta_1,\ldots,a\beta_n$. Thus it is an $\mathcal{O}_{\Fu}$ base of $\widecheck{\mathcal{R}_u}=\Delta^{-1}\mathcal{R}_u$. In particular, $\iota(\Delta) \widecheck{r_i} \in \mathcal{O}_{\mathbb{L}_w}^n$ for all $i$. We see that $M_{\iota(\Delta)}\tensor[^t]{{g}}{} \in \mathbf{M}_n(\mathcal{O}_{\mathbb{L}_w})$, so by applying the transpose we have $\mathcal{O}_{\mathbb{L}_w} \ni g M_{\iota(\Delta)}=g \left({g}^{-1} \Delta g\right)=\Delta g$.

\paragraph{Concluding the Proof.}
Let $e_1^0=g^{-1}e_1 g,\ldots,e_n^0=g^{-1} e_n g$. This is a complete set of primitive idempotents for $\mathrm{Diag}_w$. We can always choose a representative for $k_u$ in $\mathbf{GL}_1(\Omega_u)$ which we now use. Next we compute
\begin{align*}
\Psi_\sigma(k_u) \mathscr{D}_u
&=\Nrd\left(\sum_{i=1}^n \sigma.e_i\, k_u\, e_i\right)\Nrd({k_u}^{-1})\Nrd(\Delta)\\
&=\Nrd\left(\sum_{i=1}^n \sigma.e_i^0\, \left(g^{-1} k_u g\right)\, e_i^0\right)\Nrd({k_u}^{-1})\Nrd(g^{-1}\Delta g)\\
&=\Nrd\left(\sum_{i=1}^n \sigma.e_i^0\, \left(g^{-1} k_u \Delta g\right)\, e_i^0\right)\Nrd({k_u}^{-1})
\end{align*}
In the last line we have used the fact that $g^{-1}\Delta g\in \mathrm{Diag}_w$.

Because $\sum_{i=1}^n \sigma.e_i^0\, \left(g^{-1} k_u \Delta g\right)\, e_i^0$ is in the order $\Omega_w$ it has integral reduced norm. In addition, ${k_u}^{-1}\in \mathbf{GL}_1(\Omega_u)$ and it also has an integral reduced norm. This concludes the proof.
\end{proof}

Unfortunately, the proof above works only if $\mathbf{B}$ is unramified over $u$. Although it would be plausible that the conclusion of the proof holds under weaker assumptions, as is the case in rank 1, we can only produce the following result. This result for the ramified places is significantly weaker but it needs only to be applied for finitely many places of $\mathbb{F}$.

\begin{proposition}\label{prop:denom-weak}
Let $u$ be a non-archimedean place of $\mathbb{F}$ and let $w$ be a place of $\mathbb{L}$ above $u$.  Assume that $\mathcal{R}_u$ is the \emph{maximal order}. Then for all $\sigma\in W$ and for all $k_u\in K_u$ we have

\begin{equation*}
{\mathscr{D}_u}^{1+n/2}\, \Psi_\sigma(k_u) \in \mathcal{O}_{\mathbb{L}_w}
\end{equation*}
\end{proposition}
\begin{proof}
We carry notation from the proof of the previous theorem. Let $r_1,\ldots,r_n$ be an $\mathcal{O}_{\Fu}$-base for $\mathcal{R}_u$, We shall call it the \emph{integral} base. The dual base with respect to the reduced trace, $\widecheck{r_1},\ldots,\widecheck{r_n}$, is an $\mathcal{O}_{\Fu}$-base to $\widecheck{\mathcal{R}_u}=\Delta^{-1}\mathcal{R}_u$. In particular for all $i$ we have $\Delta \widecheck{r_i} \in \mathcal{R}_u$.

Using lemma \ref{lem:dual_base} we can write
\begin{align*}
\Psi_\sigma(k_u)&=\Nrd(\sum_{i=1}^n \sigma.r_i\, k_u \, \widecheck{r_i})\Nrd({k_u}^{-1})\\
&=\Nrd(\sum_{i=1}^n \sigma.r_i\, k_u \, (\Delta\widecheck{r_i}))\Nrd({k_u}^{-1}){\mathscr{D}_u}^{-1}\\
\end{align*}
We would like to transform $\sum_{i=1}^n \sigma.r_i\, k_u \, (\Delta\widecheck{r_i})$ to an element of $\Omega_w$. Both $k_u$ and $\Delta \widecheck{r_i}$ for all $i$ are in $\Omega_u$, but $\sigma.r_i$ is not necessarily so.

For each $1 \leq j \leq n$ we would like to write down in $C_w$ the element $\sigma.r_j$ in the integral base $r_1,\ldots,r_n$. To do this we need to find out the matrix corresponding to the linear transformation defined by the Weyl element $\sigma$ in the $r_1,\ldots,r_n$ base.

In $C_w$ we have another natural base. This is the base of the orthogonal primitive idempotents $e_1,\ldots,e_n$. We shall call this the \emph{orthogonal} base. Let $M$ be the matrix in the orthogonal base of the transformation sending $e_i$ to $r_i$. Because all the $r_i$ are integral we have $M\in\mathbf{M}_n(\mathcal{O}_{\mathbb{L}_w})$. By definition of the discriminant, we have $\det M=\sqrt{\mathscr{D}_u}$.
 
The matrix corresponding to $\sigma$ in the orthogonal base is just the standard permutation matrix $P_\sigma$. Hence the matrix corresponding to $\sigma$ in the \emph{integral} base is $Q_\sigma\coloneqq M P_\sigma M^{-1}=M P_\sigma M^{\mathrm{adj}} \cdot (\det M)^{-1}$. Where $M^{\mathrm{adj}} \in  \mathbf{M}_n(\mathcal{O}_{\mathbb{L}_w})$ is the adjugate matrix of $M$, which is constructed from the minors of $M$.
 
We see that $\sqrt{\mathscr{D}_u} Q_\sigma=M P_\sigma M^{\mathrm{adj}}$ is an integral matrix. This implies that for all $j$ we have 
\begin{equation*}
\sqrt{\mathscr{D}_u}\, \sigma.r_j\in {\Span}_{\mathcal{O}_{\mathbb{L}_w}}{\left(r_1,\ldots,r_n\right)}=\Omega_w\cap C_w\subset \Omega_w
\end{equation*}
 
We have thus proved that $\sum_{i=1}^n \sqrt{\mathscr{D}_u}\sigma.r_i\, k_u \, (\Delta\widecheck{r_i})\in \Omega_w$. This concludes the proof if we notice that the determinant of the matrix $\diag(\sqrt{\mathscr{D}_u},\ldots,\sqrt{\mathscr{D}_u})$ is exactly ${\mathscr{D}_u}^{n/2}$.
\end{proof}

\subsubsection{Archimedean Denominators}
We will also need an archimedean analogue of Proposition \ref{prop:denom}. The proof boils down to the following trivial linear algebra lemma.

\begin{lemma}\label{lem:Gram}
Let $V$ be an $n$-dimensional vector space over $\mathbb{C}$ endowed with an inner product $Q\colon V\times V\to \mathbb{C}$. 

Let $\epsilon_1,\ldots,\epsilon_n$ be a base of $V$.
We form the Gram matrix of $\epsilon_1,\ldots,\epsilon_n\in V$ with respect to the inner product $Q$
\begin{equation*}
\mathrm{Gr}\coloneqq \left(Q\left( \epsilon_i,\epsilon_j \right)\right)_{1\leq,i,j \leq n}
\end{equation*}

Then there exists $U\in \mathbf{U}_n(\mathbb{C})$ and a real diagonal matrix $S$ such that 
\begin{equation}
\mathrm{Gr}=U S^2 U^{-1}
\end{equation}
Evidently, $\det\mathrm{Gr}=\left(\det S\right)^2$.

Moreover, the base $\left\{U S^{-1}\epsilon_i\right\}_{1\leq i \leq n}$ is an orthonormal base with respect to $Q$.
\end{lemma}
\begin{proof}
The matrix $\mathrm{Gr}$ is Hermitian, hence there exists a  unitary diagonalization $\mathrm{Gr}=U \mathcal{D} U^{-1}$. The existence of a square root of $\mathcal{D}$ follows from the positive-definiteness of $Q$. 
\end{proof}

\begin{proposition}\label{prop:denom-arch}
Let $u$ be an archimedean place of $\mathbb{F}$ and let $w$ be a place of $\mathbb{L}$ above $u$. Denote by $|\cdot|_w$ the canonical absolute value on $\mathbb{L}$ associated to $w$.

Let $B_u\subset\mathbf{G}(\Fu)$ be a pre-compact set. Then for any $g\in B_u$ and for all $\sigma\in W$
\begin{equation*}
|\Psi^{\mathbf{H}_u}_\sigma(g)|_w \ll_{B_u} D_u
\end{equation*}
\end{proposition}
\begin{proof}
We have that either $\mathbb{L}_w\cong\Fu$ or $\mathbb{L}_w\cong\overline{\Fu}$, where $\overline{\Fu}\cong\mathbb{C}$ is the algebraic closure of $\Fu$. In any case we can identify $\overline{\mathbb{L}_w}\cong \overline{\Fu} \cong \mathbb{C}$ where all the the absolute values are compatible. Notice that the embedding of $B_u$ in $\mathbf{G}(\mathbb{C})$ is pre-compact as well. We work over $\mathbb{C}$. 

In order to define the archimedean discriminant we had to fix an inner product $Q_u$ on $\mathbf{B}(\Fu)$. We assume that this inner product is the restriction of an Hermitian inner product $Q$ on $\mathbf{B}(\mathbb{C})$ chosen as in \ref{sec:local-arch-disc}.

Let $C\coloneqq \mathbf{C}_u(\mathbb{C})$ be the maximal commutative subalgebra corresponding to $\mathbf{H}_u$ and let $e_1,\ldots,e_n$ a complete set of primitive orthogonal idempotents in $C$.

For any $g\in\mathbf{G}(\mathbb{C})$ we can write
\begin{equation*}
\Psi^{\mathbf{H}_u}_\sigma(g)=\Nrd\left(\sum_{i=1}^n e_{\sigma(i)} \cdot \Ad(g)e_i \right)
\end{equation*}

We define a utility function $\psi\colon \mathbf{G}(\mathbb{C})\times \mathbf{B}(\mathbb{C})^{n^2} \times \mathbf{B}(\mathbb{C})^{n^2}\to\mathbb{C}$
\begin{equation*}
\psi(g, b^\mathrm{L}_{i,j}, b^\mathrm{R}_{i,j})_{1\leq i,j \leq n} \coloneqq \Nrd\left(\sum_{i,j=1}^n b^\mathrm{L}_{i,j} \cdot \Ad(g) b^\mathrm{R}_{i,j} \right)
\end{equation*}

This is a a continuous function so it is bounded on pre-compact sets. Let $K_\mathbb{C}\subset\mathbf{B}(\mathbb{C})$ be the norm 1 ball with respect to\footnote{This is a decent analogue for some purposes of the compact subgroup of the non-archimedean case, hence the notation.} $Q$. Our goal
is to show ${D_u}^{-1} \Psi^{\mathbf{H}_u}_\sigma(g)\in \psi(B_u \times (K_\mathbb{C})^{n^2} \times ( K_\mathbb{C})^{n^2})$ which will prove the claim.

To calculate $D_u$ we can use any base $f_1,\ldots,f_n$ of $C$ and apply equation \eqref{eq:arch-disc-alg}
\begin{equation*}
D_u=\det\left( Q(f_i,f_j)\right)_{1\leq i,j \leq n} \cdot
\left|\det\left(\Trd(f_i f_j)\right)_{1\leq i,j \leq n}\right|_w^{-1}
\end{equation*}

We can compute the discriminant using the complete set of primitive orthogonal idempotents $e_1,\ldots,e_n$. 
\begin{equation*}
D_u=\det\left( Q(e_i,e_j)\right)_{1\leq i,j \leq n}
\end{equation*}
This is just the determinant of the Gram matrix $\mathrm{Gr}$ of $e_1,\ldots,e_n$.

By Lemma \ref{lem:Gram} we can write $\mathrm{Gr}=U S^2 U^{-1}$ with $U$ Hermitian and $S$ diagonal with $\left(\det S\right)^2=D_u$. In addition, the $Q$-orthonormal vectors $f_i\coloneqq US^{-1}e_i$, $1\leq i \leq n$, all belong to $K_\mathbb{C}$.

Write $S^{-1}=\diag(s_1,\ldots,s_n)$ and let $\varsigma\coloneqq \sum_{i=1}^n s_i e_i \in C_u$. We have
\begin{align*}
{D_u}^{-1} \Psi^{\mathbf{H}_u}_\sigma(g)
&=\Nrd(\sigma.\varsigma)\Nrd\left(\sum_{i=1}^n e_{\sigma(i)} \cdot \Ad(g)e_i \right) \Nrd(\Ad(g)\varsigma)\\
&=\Nrd\left(\sum_{i=1}^n s_{\sigma(i)}e_{\sigma(i)} \cdot \Ad(g)s_i e_i \right)
\end{align*}
\end{proof}

Notice that $s_i e_i=U^{-1} f_i$, hence each of the vectors $s_i e_i$ can be expressed as a linear combination of $f_1,\ldots,f_n$ with coefficients equal to the entries of $U^{-1}$. As $U^{-1}$ is unitarian, all its coefficients are smaller then $1$ in absolute value. This concludes the proof.

\section{Lower Bound for Asymptotic Entropy} \label{sec:entropy}
We continue with the notations of the previous sections. We now introduce the $S$-arithmetic setting which is well suited for entropy arguments.
Fix a finite set of places of $\mathbb{F}$, denoted by $S$, such that the following holds
\begin{enumerate}
\item All the archimedean places are included in $S$.
\item There is a place in $S$ over which $\mathbf{B}$ is isotropic. This condition implies that $\mathbf{G}(\mathbb{F}_S)$ is not compact.
\item The set $S$ is large enough so that 
\begin{equation*}
\mathbf{G}(\mathbb{A})=\mathbf{G}(\mathbb{F}) \cdot  \mathbf{G}(\mathbb{F}_S) \cdot \prod_{u\not\in S} K_u
\end{equation*}
That is the algebraic group $\mathbf{G}$ has class number 1 with respect to $S$ and the integral structure induced by $\Omega$.

Let $\mathbb{F}_\infty$ to be the product of the completions of $\mathbb{F}$ at all archimedean places.
Recall \cite[Theorem 5.1]{BorelFinite} that $\sdfaktor{\mathbf{G}(\mathbb{F})}{\mathbf{G}(\mathbb{A})}{\mathbf{G}(\mathbb{F}_\infty) \prod_{u<\infty} K_u}$ is finite. Let $\delta_1,\ldots,\delta_k\in \mathbf{G}(\mathbb{A})$ be a list of representatives for the double cosets. If for all $1\leq i \leq k$ the set $S$ contains all the places $u$ such that $\delta_i\not \in K_u$, then the class number $1$ assumption is trivially satisfied for $S$. In particular, $S$ can always be enlarged to satisfy the class number $1$ assumption.
\end{enumerate}
One can deduce results analogous to ours for smaller sets $S$ not fulfilling the class number $1$ assumption by taking quotients. Otherwise one can deduce all our results directly avoiding this assumption with the downside of obfuscating many statements.

Let $K^S=\prod_{u\not\in S} K_u$ and define $\Gamma= K^S \cap \mathbf{G}(\mathbb{F})$ where the intersection is taken with respect to the \emph{diagonal} embedding of $\mathbf{G}(\mathbb{F})$. Then $\Gamma$ is a lattice in the $S$-arithmetic product  $G_S \coloneqq \prod_{u\in S} \mathbf{G}(\Fu)$, again with respect to the diagonal embedding. 

Using the fact that $\mathbf{G}$ has class number one with respect to $S$, we can identify 
\begin{equation}\label{eq:adelic-to-S}
\sdfaktor{\mathbf{G}(\mathbb{F})}{\mathbf{G}(\mathbb{A})}{K^S}\xrightarrow{\sim} \lfaktor{\Gamma}{G_S}
\end{equation}

\subsection{Limit Entropy of Well Separated Orbits}
We begin by defining the Bowen ball.

\begin{definition}\label{defi:bowen-ball}
Let $a\in G_S$ be a semisimple element. 
Let $B \subset G_S$ be an identity neighborhood. For any $s<t\in\mathbb{Z}$ Define
\begin{equation*}
B^{(s,t)}\coloneqq a^{-s} B a^{s} \cap a^{-t} B a^t
\end{equation*}

If $B= \prod_{u\in S} B_u$ is a box then the $u\in S$ component of $B^{(s,t)}$ is $B_u^{(s,t)}\coloneqq a_u^{-s} B_u a_u^{s} \cap a_u^{-t} B_u a_u^t$.

We call $B^{(s,t)}$ a \emph{Bowen ball} and for $x\in\lfaktor{\Gamma}{G_S}$ we call $x\,B^{(s,t)}$ a \emph{Bowen Ball around $x$}.
\end{definition}

For a semisimple element $a\in G_S$ and $B\subset G_S$ chosen small enough we can use the Lie algebra of $G_S$ to conceptually understand the behavior of the set $B^{(-t,t)}$ for large $t>0$. The Lie algebra of $G_S$ decomposes into eigenspaces of $\Ad(a)$. All the eigenspaces on which $\Ad(a)$ acts with an eigenvalue which is not of absolute value $1$ will be either contracted or expanded under the adjoint action of $a$. In particular, the preimage of $B^{(-t,t)}$ in the Lie algebra will contract in all these eigenspaces. Hence for large $t$ the set $B^{(-t,t)}$ will be a thin tube. 

The reason we care about Bowen balls is that an exponentially decaying bound on the average measure of Bowen balls implies an entropy bound for a probability limit measure. The following proposition is an $S$-arithmetic analogue of \cite[Proposition 3.2]{ELMVPeriodic}. Its proof is the same as of the original proposition.

\begin{proposition}\label{prop:bowen-entropy}
Fix a semisimple element $a\in G_S$. 
Suppose that $\{\mu_i\}_{i=1}^\infty$ is a sequence of $a$-invariant probability measures on $\lfaktor{\Gamma}{G_S}$ converging to a probability measure $\mu$ in the weak-$*$ topology.

Assume that for some fixed $\eta>0$ we have a sequence of integers $\tau_i\to_{i\to\infty}\infty$ such that for any compact subset $F\subset \lfaktor{\Gamma}{G_S}$ there exists an identity neighborhood $B\subset G_S$ such that
\begin{align*}
\mu_i\times \mu_i &\left\{(x,y)\in F\times F \mid y\in x B^{(-\tau_i,\tau_i)} \right\}\\
&=\int_F \mu_i\left(x B^{(-\tau_i,\tau_i)} \cap F\right)\dif\mu_i(x) \ll_F \exp(-2\eta \tau_i)
\end{align*}
Then the metric entropy of the $a$-action with respect to the measure $\mu$ satisfies $h_\mu(a)\geq \eta$.
\end{proposition}

\subsection{Double Torus Quotient of a Bowen Ball}
In this section we study canonical generators of points in a Bowen ball. We will show that the canonical generators are bounded in terms of the size of the Bowen ball.

Fix a place $\uhat\in S$ and denote $\Fhat\coloneqq \mathbb{F}_\uhat$. Let $\mathbf{H}<\mathbf{G}_\Fhat$ be a maximal torus defined over $\Fhat$. 
Let $\Lhat/\Fhat$ be the splitting field of $\mathbf{H}$ and denote by $\widehat{w}$ the unique place of $\Lhat$ with canonical absolute value $|\cdot|_\what$. For $f\in\Fhat$ we have $|f|_\what=|f|_\uhat^d$ for some fixed $d\in\mathbb{N}$.
Evidently, $\mathbf{B}$ is split over $\Lhat$.

We fix $a\in \mathbf{H}(\Fhat)$, all the Bowen balls we discuss are with respect to this fixed $a$.
 
Denote by $\Psi^{\mathbf{H}}_\sigma$ the canonical generators of $\left(\dfaktor{\mathbf{H}}{\mathbf{G}_\Fhat}{\mathbf{H}}\right)_{\Lhat}$. Notice that $\mathbf{H}$ is split over $\Lhat$ so the generators are defined over $\Lhat$.
Our task is to bound $\Psi^\mathbf{H}_\sigma(g)$ when $g\in B_\uhat^{(s,t)}$.

There exists a maximal commutative subalgebra $\mathbf{C}<\mathbf{B}$ such that $\mathbf{H}$ is isogenous to $\mathbf{PGL}_1(\mathbf{C})$ under the fixed isogeny $\mathbf{G}\to\mathbf{PGL}_1(\mathbf{B})$.
Let $e^\mathbf{H}_1,\ldots,e^\mathbf{H}_n$ be a complete set of primitive orthogonal idempotents
for the split commutative algebra $\mathbf{C}(\Lhat)$. As usual we identify the absolute Weyl group of $\mathbf{H}$ with the symmetric group on $e^\mathbf{H}_1,\ldots,e^\mathbf{H}_n$.

\begin{definition}[Coordinates relative to $\mathbf{H}$]
As the algebra $\mathbf{B}$ is split over $\Lhat$ there exists an isomorphism over $\Lhat$
\begin{equation}\label{eq:xsi-def}
\xi\colon\mathbf{B}_\Lhat\to\mathbf{M}_{n,\Lhat}
\end{equation}
such that $\xi(\mathbf{C}_\Lhat)=\textbf{Diag}_\Lhat$. Choose $\xi$ so it is compatible with the fixed orderings of the primitive orthogonal idempotents in $\mathbf{C}_\Lhat (\Lhat)$ and $\textbf{Diag}_\Lhat(\Lhat)$. 

Recall that $x_{i,j}$ are the usual coordinate functions on the matrix algebra. We define 
\begin{equation*}
x^\mathbf{H}_{i,j}=x_{i,j}\circ\xi
\end{equation*}
The functions $x^\mathbf{H}_{i,j}$ are regular functions generating $\Lhat[\mathbf{B}_\Lhat]$. In particular, they are continuous with respect to the Hausdorff topology on $\mathbf{B}_\Lhat(\Lhat)$.
\end{definition}

\begin{lemma}\label{lem:generators-xx}
Let $\sigma$ be an element of the absolute Weyl group of $\mathbf{H}$.
When considered as a function of $\mathbf{G}_\Fhat$, the canonical generator $\Psi^\mathbf{H}_\sigma$  can be written as
\begin{equation*}
\Psi^\mathbf{H}_\sigma=\sign\sigma\, (\det)^{-1} \prod_{i=1}^n x^\mathbf{H}_{\sigma(i),i}
\end{equation*}
\end{lemma}
\begin{proof}
Notice that this is well defined once an isogeny $\mathbf{G}\to\mathbf{PGL}_1(\mathbf{B})$ is fixed.
The claim follows directly from Propositions \ref{prop:generators-matrix} and \ref{prop:generators-splitting}.
\end{proof}

We need to explicate the relation between roots of $\mathbf{H}$ and the set of roots $R_\sigma$ associated to $\sigma\in W_\mathbf{H}$.
The group $W_\mathbf{H}$ is the absolute Weyl group of $\mathbf{H}$ which is identified with the symmetric group over the primitive orthogonal idempotents. Recall the definition of the set of roots of $\sigma\in W_\mathbf{H}\cong S_n$
\begin{equation*}
R_\sigma\coloneqq \left\{(\sigma(i), i) \mid 1 \leq i \leq n,\,\sigma(i)\neq i \right\}
\end{equation*}

The roots of $\mathbf{H}$ are the non-trivial weights of the adjoint action of $\mathbf{H}$ on the Lie algebra $\Lie(\mathbf{G})$. Because $\mathbf{H}$ splits over $\Lhat$ the roots can be defined over $\Lhat$.
The roots $\alpha_{i,j}\colon \mathbf{H}_\Lhat \to \mathbb{G}_m$ are indexed by pairs $\left\{(i,j)\mid 1\leq i,j \leq n,\, i\neq j \right\}$. They can be defined by 
\begin{equation*}
\alpha_{i,j}\coloneqq x^\mathbf{H}_{i,i}/x^\mathbf{H}_{j,j}
\end{equation*}
Hereafter when discussing elements of $R_\sigma$ we identify the pair of integers $(\sigma(i),i)\in R_\sigma$ with the root $\alpha_{\sigma(i),i}$.

\begin{proposition}\label{prop:generators-bound}
Let $|\cdot|_\what$ be the standard absolute value associated to $\what$ on $\Lhat$. Fix $\sigma\in W_\mathbf{H}$. Let $R_\sigma$ be the set of roots associated to $\sigma$ as in Proposition \ref{prop:relations-zero}. Assume $B_\uhat\subseteq \mathbf{G}(\Fhat)$ is pre-compact. The following holds for any $g\in B_\uhat^{(s,t)}$
\begin{equation*}
\left|\Psi^\mathbf{H}_\sigma(g)\right|_\what \ll_{\mathbf{H},B_\uhat} \prod_{\alpha\in R_\sigma} \min\left(|\alpha(a^s)|_\what,|\alpha(a^t)|_\what\right)
\end{equation*}
\end{proposition}
\begin{proof}
Our starting point is Lemma \ref{lem:generators-xx}. We raise the expression for $\Psi^\mathbf{H}_\sigma$ to the $n$'th power
\begin{equation*}
\left(\Psi^\mathbf{H}_\sigma\right)^n=\left(\sign\sigma\right)^n  \prod_{i=1}^n \frac{\left(x^\mathbf{H}_{\sigma(i),i}\right)^n}{\det}
\end{equation*}
We claim that for each $1\leq i,j \leq n$ the function 
\begin{equation*}
\widetilde{x^\mathbf{H}_{i,j}}\coloneqq {\left(x^\mathbf{H}_{\sigma(i),i}\right)^n}\cdot(\det)^{-1}
\end{equation*}
is a regular function on $\mathbf{G}_\Lhat$. Recall that regular functions on $\mathbf{PGL}_1(\mathbf{B})$ extend to regular functions on $\mathbf{G}$ by the fixed isogeny $\mathbf{G}\to\mathbf{PGL}_1(\mathbf{B})$. 
Using the isomorphism $\xi$ defined in \eqref{eq:xsi-def} we can reduce this claim to the analogues claim for $\mathbf{PGL}_n$. 

We look closer at the coordinate rings. The ring of regular functions on $\mathbf{GL}_{n,\Lhat}$ is a quotient of $\Lhat[x_{i,j},(\det)^{-1}]_{1\leq i,j \leq n}$ and the regular functions on $\mathbf{PGL}_{n,\Lhat}$ are the elements of 
degree zero in $\Lhat[\mathbf{GL}_n]$, where the degree is defined the same as for regular polynomials except that the degree of $(\det)^{-1}$ is $-n$. Our functions $\widetilde{x^\mathbf{H}_{i,j}}$ are degree zero, hence they are regular.

From the definition of $\alpha_{i,j}$ and $x^\mathbf{H}_{i,j}$ we see that for any $h\in H$
\begin{equation*}
x^\mathbf{H}_{i,j} \circ \Ad(h)=\alpha_{i,j}(h)\cdot x^\mathbf{H}_{i,j}
\end{equation*}
We define $\alpha_{i,i}$ to be the trivial character for all $1\leq i \leq n$, viz.\ $x^\mathbf{H}_{i,i}$ is invariant under conjugation by $h$.

The following holds for the normalized $n$'th power functions
\begin{equation*}
\widetilde{x^\mathbf{H}_{i,j}} \circ \Ad(h)={\alpha_{i,j}(h)}^n\cdot \widetilde{x^\mathbf{H}_{i,j}}\,.
\end{equation*}

As the functions $\widetilde{x^\mathbf{H}_{i,j}}$ are regular ones they are continuous with respect to the Hausdorff topology on $\mathbf{G}(\Lhat)$, in particular $|\widetilde{x^\mathbf{H}_{i,j}}|_\uhat$ is bounded on the pre-compact set $B_\uhat$. Therefore for $g\in \Ad(h)(B_\uhat)$
\begin{equation*}
|\widetilde{x^\mathbf{H}_{i,j}}(g)|_\what\leq |\alpha_{i,j}(h)|_\what^n\cdot \sup_{g\in B_\uhat} {|\widetilde{x^\mathbf{H}_{i,j}}(g)|_\what}
\end{equation*}

We conclude that for $g\in B_\uhat^{(s,t)}$ the following estimate holds
\begin{equation}\label{eq:xtilde-estimate}
|\widetilde{x^\mathbf{H}_{i,j}}(g)|_\what \ll_{\mathbf{H},B_\uhat} \min\left(|\alpha_{i,j}(a^s)|_\what^n, |\alpha_{i,j}(a^t)|_\what^n \right)
\end{equation}
Notice that for $i=j$ the character $\alpha_{i,i}$ is trivial and the bound does not depend on $s$ and $t$.

The proposition follows from multiplying the estimates \eqref{eq:xtilde-estimate} for all $(\sigma(i),i)$, $1\leq i\leq n$, and taking the positive $n$'th root.
\end{proof}

\subsection{Small Bowen Balls} 
Recall that $\mathbf{H}<\mathbf{G}_\Fhat$ is a maximal torus defined over $\Fhat=\mathbb{F}_\uhat$. Assume henceforth that $\mathbf{H}$ is isotropic over $\Fhat$, viz.\ $\rank_\Fhat \mathbf{H}>0$. 

Let $H\coloneqq \mathbf{H}(\Fhat)$. We consider $H$ as embedded in $G_S$ in the natural way. The isotropy assumption is necessary for $H$ to have non-trivial dynamics on $\lfaktor{\Gamma}{G_S}$.

The main result of this section is Theorem \ref{thm:small-bowen} below. The entropy lower bound presented in Theorem \ref{thm: entropy-bound} follows from Theorem \ref{thm:small-bowen} using well known methods going back essentially to Linnik \cite{LinnikBook}. In their modern form they have been developed in \cite{ELMVPeriodic,ELMVPGL2,PointsOnSphere}. To formulate Theorem \ref{thm:small-bowen} we need a few definitions regarding homogeneous toral sets.

\begin{definition}
A homogeneous toral set $\lfaktor{\mathbf{G}(\mathbb{F})}{\mathbf{T}(\mathbb{A})g_\mathbb{A}}$, $g_\mathbb{A}=(g_u)_{u\in\mathcal{V}_{\mathbb{F}}}$, such that $\mathbf{H}_\uhat\coloneqq {g_\uhat}^{-1} T_\Fhat g_\uhat=\mathbf{H}$ is called an \emph{$H$-invariant} homogeneous toral set.
The projection of $Y$ to $\lfaktor{\Gamma}{G_S}$ using the map from \eqref{eq:adelic-to-S} is called an \emph{$H$-invariant} packet. Notice that $\uhat$ belongs to $S$ by assumption so that the group $H$ has a well-defined action on $\lfaktor{\Gamma}{G_S}$.
\end{definition}

\begin{definition}
Let $Y=\lfaktor{\mathbf{G}(\mathbb{F})}{\mathbf{T}(\mathbb{A})g_\mathbb{A}}$ be a homogeneous toral with local discriminants $\left(D_u\right)_{u\in\mathcal{V}_\mathbb{F}}$. Let $\mathcal{V}_\mathrm{ram}$ be the set of \emph{finite} $\mathbb{F}$-places where $\mathbf{B}$ is ramified. We define the ramified discriminant to be
\begin{equation*}
D_\mathrm{ram}\coloneqq\prod_{u\in\mathcal{V}_\mathrm{ram}} D_u
\end{equation*}
\end{definition}

\begin{definition}
Let $Y=\lfaktor{\mathbf{G}(\mathbb{F})}{\mathbf{T}(\mathbb{A})g_\mathbb{A}}$ be a homogeneous toral set. We say the $Y$ is a homogeneous toral set \emph{of maximal type} if for any finite $u$ the order $\mathcal{R}_u\coloneqq\Omega_u\cap C_u(\Fu)$ from \S \ref{sec:global-order} is maximal in the \'{e}tale algebra $\mathbf{C}_u(\mathbb{F}_u)$. Equivalently, the homogeneous toral set is of maximal type if the associated global order $\mathcal{R}$ is maximal in $\mathbb{K}$.
\end{definition}

\subsubsection{}
For $a\in H$ and an $a$-invariant probability measure $\mu$ on $\lfaktor{\Gamma}{G_S}$ we have denoted by $h_\mu(a)$ the Kolmogorov-Sinai entropy of the measure $\mu$ with respect to the $\mathbb{Z}$ action by $a$. Let $m$ be the Haar measure on $\lfaktor{\Gamma}{G_S}$ and let $\mathfrak{g}_{-}$ be the subspace of the Lie algebra $\Lie(G(\Fhat))$ which is spanned over $\Lhat$ by the eigenvalues of $\Ad_{a}$ of modulus $<1$. Then
\begin{equation*}
h_\mathrm{Haar}(a)\coloneqq h_m(a)=-\log|\det{\Ad}_a\restriction_{\mathfrak{g}_-}|_\uhat
=\frac{1}{2f_\what}\sum_{1\leq i\neq j \leq n} |\log{|\alpha_{i,j}(a)|_\what}|
\end{equation*}
Where the absolute value on the outside each $\log{|\alpha_{i,j}(a)|_\uhat}$ is the regular absolute value on $\mathbb{R}$ and $|x|_\what=|x|_\uhat^{f_\what}$ for each $x\in\Fhat$, i.e.\ $f_\what$ is the inertia degree of $\Lhat/\Fhat$.

\begin{theorem}\label{thm:small-bowen}
Let $Y={\mathbf{G}(\mathbb{F})}{\mathbf{T}(\mathbb{A})g_\mathbb{A}}$ be an $H$-invariant homogeneous toral set of maximal type with discriminant $D$ and ramified local discriminant $D_\mathrm{ram}$. Let $\mathbb{L}$ be the splitting field of $\mathbf{T}$. Assume $\mathfrak{G}\coloneqq\Gal(\mathbb{L}/\mathbb{F})$ is 2-transitive.

Let $a\in H$ and fix $B=\prod_{u\in S} B_u$ such that for all $\uhat\neq u\in S$ non-archimedean the set $B_u$ is contained in $K_u$.

There exists a constant $\kappa$ depending only on $B$, $\mathbb{F}$ and $\mathbf{H}$ such that if $\tau>0$ satisfies
\begin{equation}\label{eq:tau-ineq}
2\frac{h_\mathrm{Haar}(a)}{n-1}\tau> \log\left(D D_\mathrm{ram}^{n/2}\right)+\kappa
\end{equation}
Then  $\mathbf{G}(\mathbb{F}) \cap \mathbf{T}(\mathbb{A}) g_\mathbb{A} \left(B_u^{(-\tau,\tau)} \times K^S\right) {g_\mathbb{A}}^{-1}  \mathbf{T}(\mathbb{A})$ is contained in $\mathbf{T}(\mathbb{F})$.
\end{theorem}
\begin{proof}
For now let $\kappa$ be a free variable and let 
\begin{equation*}
\lambda\in\mathbf{G}(\mathbb{F}) \cap \mathbf{T}(\mathbb{A}) g_\mathbb{A} \left(B_u^{(-\tau,\tau)} \times K^S\right) {g_\mathbb{A}}^{-1}  \mathbf{T}(\mathbb{A})
\end{equation*}
Denote by  $W$ be the absolute Weyl group of $\mathbf{T}$. 

\paragraph{Galois orbits.} Recall that by Proposition \ref{prop:galois-weyl} we can consider $\mathfrak{G}$ as a subgroup of $W$.

Fix $\sigma\in W$. By Proposition \ref{prop:psi-galois} the Galois orbit of $ \Psi^\mathbf{T}_\sigma(\lambda)\in\mathcal{O}_{\mathbb{L}}$ is
\begin{equation*}
\left\{\Psi^\mathbf{T}_{\tau \sigma \tau^{-1}}(\lambda) \mid \tau\in\mathfrak{G} \right\}
\end{equation*}

Denote $\mathscr{C}\coloneqq\left\{\tau\sigma\tau^{-1}\mid \tau\in\mathfrak{G} \right\}$ and let $|\mathscr{C}|$ be the cardinality of $\mathscr{C}$. We now define
\begin{equation*}
\Psi_\mathscr{C}(\lambda)\coloneqq \prod_{\omega\in\mathscr{C}} \Psi^\mathbf{T}_\omega(\lambda)\in\mathbb{F}
\end{equation*}

\paragraph{Archimedean bound for $\Psi_\mathscr{C}(\lambda)$.}
For each archimedean place $u\neq \uhat$ of $\mathbb{F}$ let $m$ be the number of places $w|u$ of $\mathbb{L}$. Recall that we treat each pair of conjugate complex embedding as genuinely different places. Applying Proposition \ref{prop:denom-arch} and using that $B_u^{(-\tau,\tau)}= B_u$ for all archimedean $u\neq\uhat$ we have
\begin{equation*}
|\Psi_\mathscr{C}(\lambda)|_u^m=\prod_{w|u}|\Psi_\mathscr{C}(\lambda)|_w\leq \left(\kappa_0 D_u\right)^{|\mathscr{C}|m}
\end{equation*}
The constant $\kappa_0>0$ depends on $B$ only. Set $\kappa_1\coloneqq\max(\kappa_0,1)$.

Taking the positive real root of order $m$ we have the necessary bound
\begin{equation}\label{eq:bound-arch}
|\Psi_\mathscr{C}(\lambda)|_u\leq \kappa_1^{|\mathscr{C}|} D_u^{|\mathscr{C}|}
\end{equation}

\paragraph{Nonarchimedean bound for $\Psi_\mathscr{C}(\lambda)$.}
Recall that if $u\in S$ then  $B_u\subseteq K_u$.
For all non-archimedean $u\neq\uhat$ we have
\begin{equation*}
\lambda\in \mathbf{T}(\Fu) g_u K_u {g_u}^{-1}  \mathbf{T}(\Fu)
\end{equation*} 
For each non-archimedean place $u\neq\uhat$ of $\mathbb{F}$ where $\mathbf{B}$ is ramified we apply Proposition \ref{prop:denom} 
\begin{align*}
|\Psi_\mathscr{C}(\lambda)|_u^{\left[\mathbb{L}\colon\mathbb{F}\right]}
&=|{\Nr}_{\mathbb{L}/\mathbb{F}}\Psi_\mathscr{C}(\lambda)|_u=\prod_{w|u}|\Psi_\mathscr{C}(\lambda)|_w\\
&\leq \prod_{w|u}|\mathscr{D}_u|_w^{-|\mathscr{C}|} =|{\Nr}_{\mathbb{L}/\mathbb{F}}\mathscr{D}_u|_u^{-|\mathscr{C}|}
=D_u^{|\mathscr{C}|\left[\mathbb{L}\colon\mathbb{F}\right]}
\end{align*}

If $\mathbf{B}$ is ramified at $u\neq\uhat$ we apply the weaker Proposition \ref{prop:denom-weak} and have similarly
\begin{equation*}
|\Psi_\mathscr{C}(\lambda)|_u^{\left[\mathbb{L}\colon\mathbb{F}\right]}
\leq D_u^{(1+n/2)|\mathscr{C}|\left[\mathbb{L}\colon\mathbb{F}\right]}
\end{equation*}

Taking the positive real root of order $\left[\mathbb{L}\colon\mathbb{F}\right]$ we have the necessary bounds
\begin{equation} \label{eq:bound-nonarch}
|\Psi_\mathscr{C}(\lambda)|_u\leq  
\begin{cases}
D_u^{|\mathscr{C}|} & \textrm{ If $\mathbf{B}$ is unramified at $u$} \\
D_u^{(1+n/2)|\mathscr{C}|} & \textrm{ Otherwise}
\end{cases}
\end{equation}

\paragraph{Exponential bound.}
It will be now important to choose $\sigma$ that has \emph{no fixed points}. We are going to calculate the exponential bound on
$\left| \Psi_\mathscr{C}(\lambda)\right|_\uhat$.

Let $\Lhat/\Fhat$ be the splitting field of $\mathbf{H}$ and denote by $\what$ place of $\Lhat$.
Using Proposition \ref{prop:generators-bound} we deduce that
\begin{equation}\label{eq:exp-prebound-split}
\left|\Psi_\mathscr{C}(\lambda)\right|_\what \leq \kappa_2^{|\mathscr{C}|} \exp\left(-\tau\sum_{\omega\in\mathscr{C}}\sum_{\alpha\in R_\omega} \left|\log |\alpha(a)|_\what \right| \right)
\end{equation}
Where $k_2>0$ depends only on $B$. Notice that the outer absolute value of $|\log |\alpha(a)|_\what|$ is the usual absolute value on $\mathbb{R}$.

As $\Psi_\mathscr{C}(\lambda)$ belongs to $\mathbb{F}$ it can be considered as an element of $\Fhat$. Inequality \eqref{eq:exp-prebound-split} implies
\begin{equation*}
\left|\Psi_\mathscr{C}(\lambda)\right|_\uhat \leq \kappa_2^{|\mathscr{C}|} \exp\left(-\frac{\tau}{f_\what}\sum_{\omega\in\mathscr{C}}\sum_{\alpha\in R_\omega} \left|\log |\alpha(a)|_\what \right| \right)
\end{equation*}

How many times each root of $\mathbf{H}$ appears in the sum above? Because $\mathfrak{G}$ is 2-transitive it is easy to see that all the roots appear in the sum with the same multiplicity. There are $n(n-1)$ roots in total. Each $\omega\in\mathscr{C}$ contributes the same number of roots to the sum as $\sigma$, which contributes exactly $n$ of them as it has \emph{no fixed points}. Hence each root appears in the sum with multiplicity 
\begin{equation*}
\frac{n|\mathscr{C}|}{n(n-1)}=\frac{|\mathscr{C}|}{n-1}
\end{equation*}

In the expression for $2 h_{\mathrm{Haar}}(a)$ each root appears with multiplicity $1$. 
We conclude that
\begin{equation}\label{eq:bound-exp}
\left|\Psi_\mathscr{C}(\lambda)\right|_\uhat \leq \kappa_2^{|\mathscr{C}|} 
\exp\left(- 2\frac{h_\mathrm{Haar}(a)}{n-1} |\mathscr{C}| \tau\right)
\end{equation} 

\paragraph{Triviality of $\Psi_\mathscr{C}(\lambda)$}
For each place of $u\neq\uhat$ we apply either bound \eqref{eq:bound-arch} or \eqref{eq:bound-nonarch} appropriately. For the place $u=\uhat$ we apply the exponential bound \eqref{eq:bound-exp}. Combining these we have 
\begin{align*}
{\Nr}_\mathbb{A}\Psi_\mathscr{C}(\lambda)
&\coloneqq \prod_{u\in\mathcal{V}_\mathbb{F}}|\Psi_\mathscr{C}(\lambda)|_u\\
&\leq \left(\kappa_1^{\left[\mathbb{F}\colon\mathbb{Q}\right]} \kappa_2\right)^{|\mathscr{C}|}{D_\uhat}^{-{|\mathscr{C}|}} D^{|\mathscr{C}|} D_\mathrm{ram}^{n/2|\mathscr{C}|}  
\exp\left(- 2\frac{h_\mathrm{Haar}(a)}{n-1} |\mathscr{C}| \tau\right)
\end{align*}
The local discriminant $D_\uhat$ depends only on $\mathbf{H}$. Define $\kappa_3=\kappa_1^{\left[\mathbb{F}\colon\mathbb{Q}\right]} \kappa_2 D_\uhat^{-1}$, it depends only on $B$, $\mathbb{F}$ and $\mathbf{H}$.

We see that if
\begin{align*}
\kappa_3^{|\mathscr{C}|}  D^{|\mathscr{C}|}D_\mathrm{ram}^{n/2|\mathscr{C}|} \exp\left(- 2\frac{h_\mathrm{Haar}(a)}{n-1} |\mathscr{C}| \tau\right) &<1 \\
\Longleftrightarrow 2\frac{h_\mathrm{Haar}(a)}{n-1} \tau &> \log\left(DD_\mathrm{ram}^{n/2}\right)+ \log \kappa_3
\end{align*}
Then ${\Nr}_\mathbb{A}\Psi_\mathscr{C}(\lambda)<1$. But $\Psi_\mathscr{C}(\lambda)\in\mathbb{F}$ and from the product formula for the number field $\mathbb{F}$ we deduce that $\Psi_\mathscr{C}(\lambda)=0$.

\paragraph{Finishing the proof.}
We are now ready to set $\kappa=\log\kappa_3$.
We have shown that if inequality \eqref{eq:tau-ineq} holds then $\Psi_\mathscr{C}(\lambda)=0$. This in turn implies that there exists $\tau\in\mathfrak{G}$ such that $\Psi^\mathbf{T}_{\tau\sigma\tau^{-1}}(\lambda)=0$.

By Corollary \ref{cor:psi-galois-2} the projection of $\lambda\in\mathbf{G}(\mathbb{F})$ to $\dfaktor{\mathbf{T}}{\mathbf{G}}{\mathbf{T}}$ is $\mathbf{T}e\mathbf{T}$. Now we can use Proposition \ref{prop:identity-fiber} to conclude that there exits $t\in\mathbf{T}(\mathbb{F})$ such that 
$\lambda=t$ and the claim is proved.
\end{proof}
\subsection{Proof of the Entropy Lower Bound}
We now prove our lower bound on the asymptotic entropy. The reader should remember that our final argument is an extension of a rank $1$ argument which is weaker then Linnik's powerful method. Indeed, in rank $1$ Linnik's results are akin to an optimal bound on the entropy \cite{LinnikBook,ELMVPGL2}. A generalization of Linnik's basic lemma and a finer analysis of integral points on the variety $\dfaktor{\mathbf{T}}{\mathbf{G}}{\mathbf{T}}$ in conjugation with the Galois action remain open.

In addition we have no contribution to the question of escape of mass in the non-cocompact case. We assume a priori that mass does not escape. This is unnecessary if $\mathbf{B}$ is a division algebra as congruence lattices are co-compact in that case.

\begin{theorem-quote}[\ref{thm: entropy-bound}]
Suppose we have a sequence of $H$-invariant homogeneous toral sets of maximal type  $Y_i=\lfaktor{\mathbf{G}(\mathbb{F})}{\mathbf{T}_i(\mathbb{A})g_i}$ with $\mathbf{T}_i$ a torus defined and anisotropic over $\mathbb{F}$ and $g_i\in G_\mathbb{A}$. Let $\mathbb{L}_i/\mathbb{F}$ be the splitting field of $\mathbf{T}_i$. We assume for all $i$ that $\Gal(\mathbb{L}_i/\mathbb{F})$ is 2-transitive.

Denote by $D_i$ the global discriminant of $Y_i$. Assume $D_i\to_{i\to\infty}\infty$. Let $D_{\mathrm{ram},i}$ be the ramified discriminant of $Y_i$ which is the product of local discriminant at finite places where $\mathbf{B}$ is ramified. Suppose $D_{\mathrm{ram},i}=D_i^{o(1)}$

Let $\mu_{i}$ be the probability measure on $\lfaktor{\Gamma}{G_S}$ induced by the probability measure on the homogeneous toral set $Y_i$ using \eqref{eq:adelic-to-S}. If $\mu_i$ converges in the weak-$*$ topology to a probability measure $\mu$ on $\lfaktor{\Gamma}{G_S}$ then for any $a\in H$ we have
\begin{equation*}
h_\mu(a)\geq \frac{h_\mathrm{Haar}(a)}{2(n-1)}
\end{equation*}
\end{theorem-quote}

\begin{remark} \label{rem:volume}
The volume of the packet is known to be equal to ${D_i}^{1/2+o(1)}$, see \cite[Theorem 4.8]{ELMVCubic}. Although the proof in \cite{ELMVCubic} is written for the case $\mathbf{B}=\mathbf{M}_n$ it applies verbatim to any central simple algebra $\mathbf{B}$. Notice that any torus in $\mathbf{B}$ over a local field can be embedded in $\mathbf{M}_n$ with the same local order\footnote{In the archimedean case with $\mathcal{R}_u$ comparable as convex sets.} $\mathcal{R}_u$. The computation of the volume depends only on this data if the norms are chosen consistently.
\end{remark}

\begin{proof}[Proof of Theorem \ref{thm: entropy-bound}]
The theorem follows from combining Theorem \ref{thm:small-bowen} with Proposition \ref{prop:bowen-entropy}  (\cite[Proposition 3.2]{ELMVPeriodic}).

\paragraph{The necessary bound on Bowen balls.}
We define
\begin{equation*}
\eta=\frac{h_\mathrm{Haar}(a)}{n-1}
\end{equation*}
Choose $B$ as in Theorem \ref{thm:small-bowen} and such that $B_\uhat^{(-\tau,\tau)}\subseteq B_\uhat$ for all\footnote{This can be done by choosing in the Lie algebra $\Lie(\mathbf{G})(\Fhat)$ a small enough identity neighborhood fulfilling an analogues condition and taking its image under the exponent.}  $\tau\geq 0$. Proposition \ref{prop:bowen-entropy} implies that the theorem would follow if for each $\varepsilon>0$ we can find a sequence of positive integers $\tau_i\to_{i\to\infty}\infty$ such that for any compact $F\subset\lfaktor{\Gamma}{G_S}$
\begin{equation}\label{eq:required-bound}
\int_F \mu_i\left(x B^{(-\tau_i,\tau_i)}\right)\dif\mu_i(x) \ll_{F,\mathbb{F},B,\mathbf{H},\varepsilon}
\exp(-(\eta-\varepsilon) \tau_i)
\end{equation}

Let $Z_i\subset \lfaktor{\Gamma}{G_S}$ be the projection of $Y_i$. By choosing $\tau_i$ appropriately, we are going to show that inequality \eqref{eq:required-bound} holds individually for any Bowen ball $x B^{(-\tau_i,\tau_i)}$ with $x\in Z_i$. This will imply the theorem.

\paragraph{Volume of a small Bowen ball.}
Theorem \ref{thm:small-bowen} suggests that we choose 
\begin{equation} \label{eq:tau-condition}
\tau_i=\frac{\log \left(D_i D_{\mathrm{ram},i}^{n/2}\right)}{2\eta} +\kappa+1=\frac{\log \left(D_i\right)(1+o(1))}{2\eta} +\kappa+1
\end{equation}

Fix $x\in Z_i$ and let $q=\delta_0 t_0 g_i k_0\in G(\mathbb{F})\mathbf{T}_i(\mathbb{A})g_iK^S$ be an adelic representative of $x$.
We can now bound the measure of the Bowen ball around $x$ corresponding to $\tau_i$. Recall that the measure $\mu_i$ is the projection of the adelic measure on the homogeneous toral set to $\lfaktor{\Gamma}{G_S}$. For simplicity we use the notation $\mu_i$ for both measures. Also, we denote by $\mu_{\mathbf{T}_i}$ the Haar measure on $\mathbf{T}_i$ normalized so that $\lfaktor{\mathbf{T}_i(\mathbb{F})}{\mathbf{T}_i(\mathbb{A})}$ has volume 1.

\begin{align}
&\mu_i\left(\Gamma x B^{(-\tau_i,\tau_i)}\right)
=\mu_i\left(\mathbf{G}(\mathbb{F}) \delta_0 t_0 g_i k_s  B^{(-\tau_i,\tau_i)}\times K^S\right)
\nonumber\\
&=\mu_i\left(\mathbf{G}(\mathbb{F}) t_0 g_i B^{(-\tau_i,\tau_i)}\times K^S\right)
\nonumber\\
&=\mu_{\mathbf{T}_i}\left(t\in \mathbf{T}_i(\mathbb{A}) 
\mid \exists\delta\in \mathbf{G}(\mathbb{F})\colon 
\delta t g_i\in \delta_0 t_0 g_i B^{(-\tau_i,\tau_i)}\times K^S  \right)
\nonumber\\
&=\mu_{\mathbf{T}_i}\left(t\in \mathbf{T}_i(\mathbb{A}) 
\mid \exists\delta\in \mathbf{G}(\mathbb{F})\colon 
\delta_0^{-1} \delta t\in t_0 \cdot g_i \left( B^{(-\tau_i,\tau_i)} \times K^S \right) {g_i}^{-1} \right) \label{eq:bowen-volume-adelic}
\end{align}
For any $\delta\in\mathbf{G}(\mathbb{F})$ contributing to the measure of the Bowen ball in \eqref{eq:bowen-volume-adelic} we have that $\lambda=\delta_0^{-1}\delta$ fulfills the condition of Theorem \ref{thm:small-bowen}, so $\lambda\in\mathbf{T}_i(\mathbb{F})$. In particular, it is enough to take in \eqref{eq:bowen-volume-adelic} only $\delta=\delta_0$. This implies
\begin{align}
\mu_i\left(\Gamma x B^{(-\tau_i,\tau_i)}\right)
&=\mu_{\mathbf{T}_i}\left(t\in \mathbf{T}_i(\mathbb{A}) 
\mid t\in t_0\cdot g_i \left( B^{(-\tau_i,\tau_i)} \times K^S\right) {g_i}^{-1} \right)
\nonumber\\
&=\mu_{\mathbf{T}_i}\left(t\in \mathbf{T}_i(\mathbb{A}) 
\mid t\in  g_i \left( B^{(-\tau_i,\tau_i)} \times K^S\right) {g_i}^{-1} \right)
\nonumber\\
&\leq\mu_{\mathbf{T}_i}\left(t\in \mathbf{T}_i(\mathbb{A}) 
\mid t\in  g_i \left( B \times K^S\right) {g_i}^{-1} \right)
\nonumber\\
&\asymp_B \vol(Y_i)^{-1}=D_i^{-1/2-o(1)} \label{eq:single-ball-volume}
\end{align}
Where the last equality is \cite[Theorem 4.8]{ELMVCubic} (see Remark \ref{rem:volume}).

Combining \eqref{eq:tau-condition} and \eqref{eq:single-ball-volume} and setting $\kappa'=2\eta(\kappa+1)$
we see that 
\begin{align*}
\mu_i\left(x B^{(-\tau_i,\tau_i)}\right)
&\leq \exp\left(-\left(\frac{1}{2}+o(1)\right)\log D_i\right)\\
&=\exp\left(-\left(\frac{1}{2}+o(1)\right)(2\eta \tau_i/(1+o(1))-\kappa')\right)\\
&\ll_{\kappa',\varepsilon} \exp\left(-(\eta-\varepsilon) \tau_i\right)
\end{align*}
This concludes the proof of Theorem \ref{thm: entropy-bound}.
\end{proof}

\section{Rigidity of Limit Measures}\label{sec:rigidity}
We now combine Theorem \ref{thm: entropy-bound} with measure rigidity for higher rank torus actions. Specifically, we are going to use the results of \cite{EntropySArithmetic} which generalizes \cite{EKL}.

The following proposition is a simple analogue of \cite[Theorem 5.1]{ELMVPeriodic} combined with the improved entropy bounds of Theorem \ref{thm: entropy-bound}. Its proof is standard and follows the same lines as  \cite[Theorem 5.1]{ELMVPeriodic}.
\begin{proposition}\label{prop:limit-rigidity}

We assume for simplicity $\mathbb{F}=\mathbb{Q}$.
Let $\uhat\in S$ and denote $\Fhat\coloneqq \mathbb{F}_\uhat$. Fix a maximal torus $\mathbf{H}<\mathbf{G}_\Fhat$ defined over $\Fhat$ and set $H\coloneqq \mathbf{H}(\Fhat)$. Assume that $\mathbf{H}$ is \emph{split} over $\Fhat$.

Suppose we have a sequence of $H$-invariant homogeneous toral sets of maximal type  $Y_i=\lfaktor{\mathbf{G}(\mathbb{F})}{\mathbf{T}_i(\mathbb{A})g_i}$ with $\mathbf{T}_i$ a torus defined and anisotropic over $\mathbb{F}$ and $g_i\in G_\mathbb{A}$. Let $\mathbb{L}_i/\mathbb{F}$ be the splitting field of $\mathbf{T}_i$. We assume for all $i$ that $\Gal(\mathbb{L}_i/\mathbb{F})$ is 2-transitive.

Denote by $D_i$ the global discriminant of $Y_i$ and by $D_{\mathrm{ram},i}$ the ramified discriminant of $Y_i$. Assume $D_i\to_{i\to\infty}\infty$ and $D_{\mathrm{ram},i}=D_i^{o(1)}$. 

Let $\mu_{i}$ be the probability measure on $\lfaktor{\Gamma}{G_S}$ induced by the probability measure on the homogeneous toral set $Y_i$. If $\mu_i$ converges in the weak-$*$ topology to a probability measure $\mu$ on $\lfaktor{\Gamma}{G_S}$ then 
\begin{equation*}
\mu\geq \int \nu \dif\tau(\nu)
\end{equation*}
Where $\tau$ is a finite measure on the space of $H$-invariant Borel probability measures on $\lfaktor{\Gamma}{G_S}$ and
\begin{enumerate}
\item For $\tau$-almost every $\nu$ there is a reductive algebraic subgroup $\mathbf{L}<\mathbf{G}$ defined and anisotropic over $\mathbb{Q}$ and some $g\in G_S$ such that $\nu$ is the unique ${g}^{-1} \mathbf{L}(\mathbb{Q}_S) g$ invariant measure supported on the periodic orbit $\lfaktor{\Gamma}{\mathbf{L}(\mathbb{Q}_S)g}$. Moreover, $H<{g_{\uhat}}^{-1} \mathbf{L}(\Fhat) g_{\uhat}$ and $\nu$ is $H$ ergodic.
\item For any $a\in H$ we have
\begin{equation*}
\int h_{\nu}(a) \dif\tau(\nu)=h_\mu(a)\geq  \frac{h_\mathrm{Haar}(a)}{2(n-1)}
\end{equation*}
In particular, $\int  \dif\tau(\nu)\geq \frac{1}{2(n-1)}$. 
\end{enumerate}
\end{proposition}
\begin{remark}\label{rem:BdS}
The reductive subgroups $\mathbf{L}<\mathbf{G}$ appearing in the decomposition above include over $\Fhat$ a maximal split torus of $\mathbf{G}_\Fhat$. In particular, $\mathbf{L}_\Fhat$ is a reductive subgroup of maximal rank of the split group $\mathbf{G}_\Fhat$ which has type $\mathrm{A}_n$. By the Borel-de Siebenthal algorithm \cite{BdS} adapted to algebraic groups, $\mathbf{L}_\Fhat$ is an almost-direct product of a torus and groups of type $\mathrm{A}_k$ with $k\leq n$.
\end{remark}

\begin{proof}
Decompose $\mu$ into $H$-invariant and ergodic probability measures
\begin{equation*}
\mu=\int \nu \dif \widetilde{\tau}(\nu)
\end{equation*}
where $\widetilde{\tau}$ is a probability measure on the space of $H$-invariant Borel probability measures on $\lfaktor{\Gamma}{G_S}$ supported on the $H$ \emph{ergodic} measures. Let $\mathcal{P}_{>0}$ be the Borel set\footnote{
To see that this is a Borel set in the space of probability measures fix a countable dense subset $\left\{a_i\right\}_{i=1}$ of the separable space $H$. It holds that $\mathcal{P}_{>0}=\bigcup_{i=1}^\infty \left\{\nu \mid h_\nu(a_i)>0\right\}$. Each sets in this union is Borel because we can compute the entropy using a countable set of finite partitions of $\lfaktor{\Gamma}{G_S}$ that generate together the whole $\sigma$-algebra.} of $H$-invariant probability measures $\nu$ such that $h_\nu(a)>0$ for \emph{some} $a\in H$.  Define 
\begin{equation*}
\tau=\widetilde{\tau}\restriction_{\mathcal{P}_{>0}} 
\end{equation*}
Obviously, $\mu\geq \int \nu \dif\tau(\nu)$ and (2) follows from the linearity of the entropy with respect to a fixed $a\in H$ when considered as a function on the convex set of Borel probability measures on $\lfaktor{\Gamma}{G_S}$.

Assertion (1) is a direct consequences of \cite[Theorem 1.1 and Remark after Corollary 1.2]{EntropySArithmetic} combined with the fact that, by definition of $\tau$, $h_\nu(a)>0$ for some $a\in H$ for $\tau$-almost all $\nu$.
\end{proof}

The following corollary is a modest qualitative restriction on possible limit measures arising from the entropy bound.
\begin{corollary}\label{cor:no-single-periodic}
Fix $R\in\mathbb{N}$. Let $\mu\geq\int \nu \dif\tau(\nu)$ be as in Proposition \ref{prop:limit-rigidity}. Let $\mathbf{L}<\mathbf{G}$ be a reductive subgroup of maximal rank, such that the rank of all its simple parts (of type $A_k$) is less then $R$.

Set $N_R=(R+2)(R+1)R-1$. If the absolute rank of $\mathbf{G}$ is greater then $N_R-1$, then $\tau(\{\nu\})<1$ for each $\nu$ supported on a periodic orbit of a conjugate of  $\mathbf{L}(\mathbb{Q}_S)$, namely, $\mu$ is not supported on a \emph{single} periodic orbit of type $\mathbf{L}$.
\end{corollary}
\begin{proof}
Assume $\nu$ is an $H$-invariant and ergodic probability measure on $\lfaktor{\Gamma}{G}$ corresponding to the periodic orbit $\lfaktor{\Gamma}{\mathbf{L}(\mathbb{Q}_S)g}$. Set $\widetilde{\mathbf{L}}\coloneqq g^{-1}\mathbf{L}_\Fhat g$ -- a reductive split algebraic group over $\Fhat$. Recall that $\mathbf{H}<\widetilde{\mathbf{L}}$ is a maximal split torus.

Let $\jmath\colon \mathbf{SL}_{k_1,\Fhat}\times\ldots\times\mathbf{SL}_{k_l,\Fhat}\times \mathbf{T}_0\to \widetilde{\mathbf{L}}$, where $\mathbf{T}_0$ is a torus defined over $\Fhat$, be the obvious isogeny from the simply connected cover to $\widetilde{\mathbf{L}}$, see Remark \ref{rem:BdS}. Let $\mathbf{H}_i<\mathbf{SL}_{k_i,\Fhat}$ a maximal torus whose image under $\jmath$ is contained in $\mathbf{H}$. By replacing $\jmath$ with a composition of $\jmath$ with an inner automorphism of $\widetilde{\mathbf{L}}$ we can assure that $\mathbf{H}=\jmath\left(\prod_{i} \mathbf{H}_i \times \mathbf{T}_0\right)$.

Denote by $q$ the cardinality of the residue field of $\Fhat$ if $\uhat$ is non-archimedean or $q\coloneqq \exp(1)$ if $\uhat$ is archimedean. Then by choosing an element $a_i\in \mathbf{H}_i(\Fhat)$ and $t_0\in\mathbf{T}_0(\Fhat)$ similarly to Remark \ref{rem:elmv-bound}\footnote{Instead of $\exp((n-j)/2)$ one puts $\varpi^{(n-j)/2}$ where $\varpi$ is a uniformizer for $\Fhat$ if $\uhat$ is non-archimedean and $\varpi=\exp(1)$ otherwise.} for $a\coloneqq\jmath(a_1,\ldots,a_l,t_0)\in H$ we have
\begin{align*}
&h_\mathrm{Haar}(a)= \frac{(n+1)n(n-1)}{6}\log q\\
&h_\nu(a)=\sum_{i=1}^l \frac{(k_i+1)k_i(k_i-1)}{6}\log q\leq \frac{n}{2} \frac{(R+2)(R+1)R}{6}\log q 
\end{align*}
and
\begin{align}
&h_\mu(a)/h_\mathrm{Haar}(a)\geq\frac{1}{2(n-1)} \label{eq:h-bound1}\\
&h_\nu(a)/h_\mathrm{Haar}(a)\leq \frac{(R+2)(R+1)R}{2(n+1)(n-1)} \label{eq:h-bound2}
\end{align}

Using the inequality
\begin{align*}
h_\mu(a)&=\int h_{\nu'}(a) \dif\tau(\nu')\leq h_\nu(a)\nu(\{\tau\})+h_\mathrm{Haar}(a)(1-\nu(\{\tau\}))\\
&\implies \nu(\{\tau\})\leq \frac{h_\mathrm{Haar}(a)-h_\mu(a)}{h_\mathrm{Haar}(a)-h_\nu(a)}
\end{align*}
together with inequalities \eqref{eq:h-bound1} and \eqref{eq:h-bound2} we conclude that $\nu(\{\tau\})<1$ if 
$n>(R+2)(R+1)R-1$.
\end{proof}

\bibliographystyle{alpha}
\bibliography{bib_torus_lower}

\end{document}